\documentclass[a4paper,10pt,sort&compress]{amsart}

\usepackage[utf8]{inputenc}
\usepackage[english]{babel}
\usepackage{hyperref}
\hypersetup{
	colorlinks={true},
	citecolor={NavyBlue},
	linkcolor={BrickRed}
	}
\usepackage[nameinlink,capitalize]{cleveref}
\usepackage[dvipsnames]{xcolor}
\usepackage{graphics}

\crefname{equation}{}{}

\usepackage{dsfont}
\usepackage{tikz}

\usepackage{booktabs}
\usepackage{amssymb,amsmath,amsfonts}
\usepackage{algorithm}
\usepackage{algorithmicx,algpseudocode}
\algrenewcommand{\algorithmiccomment}[1]{\hspace{-19pt} {\normalsize\color{black}$\triangleright\;$\textit{#1}}}
\usepackage{xypic}

\usepackage[margin=1.2in]{geometry}
\usepackage{enumitem}

\DeclareMathAlphabet{\mathpzc}{OT1}{pzc}{m}{it}
\newcommand{\pzc}[1]{\mathpzc{#1}}

\newcommand{\RR}{\mathbb{R}}

\newcommand{\NN}{\mathbb{N}}

\newcommand{\calD}{\mathcal{D}}

\newcommand{\bbP}{\mathbb{P}}
\newcommand{\bbN}{\mathbb{N}}

\newcommand{\ba}{\mathbf{a}}
\newcommand{\bb}{\mathbf{b}}

\newcommand{\bj}{\mathbf{j}}
\newcommand{\bi}{\mathbf{i}}
\newcommand{\bk}{\mathbf{k}}
\newcommand{\bn}{\mathbf{n}}
\newcommand{\br}{\mathbf{r}}

\newcommand{\tensor}[1]{\pzc{#1}}
\newcommand{\vect}[1]{\mathbf{#1}}
\newcommand{\Var}[1]{\mathcal{#1}}

\newcommand{\rk}{\mathrm{rank}}

\DeclareMathOperator{\bighadamard}{\vcenter{\hbox{\scalebox{2.1}{$\circledast$}}}}
\DeclareMathOperator{\hadamard}{\circledast}
\newcommand{\bhadamard}[2]{\overset{#2}{\underset{#1}{\bighadamard}}\;}

\theoremstyle{plain}
\newtheorem{theorem}{Theorem}[section]
\newtheorem*{theorem*}{Theorem}
\newtheorem{proposition}[theorem]{Proposition}
\newtheorem{corollary}[theorem]{Corollary}
\newtheorem{lemma}[theorem]{Lemma}

\newtheorem{claim}[theorem]{Claim}

\theoremstyle{definition}
\newtheorem{definition}[theorem]{Definition}

\newtheorem{remark}[theorem]{Remark}

\numberwithin{equation}{section}

\title[Hadamard--Hitchcock decompositions]{Hadamard--Hitchcock decompositions:\\ identifiability and computation}
\author{Alessandro Oneto}
\address[A. Oneto]{U. Trento, Department of Mathematics, Via Sommarive 14, 38123 Trento, Italy.}
\email{alessandro.oneto@unitn.it}

\author{Nick Vannieuwenhoven}
\address[N. Vannieuwenhoven]{KU Leuven, Department of Computer Science, Celestijnenlaan 200A, B-3001 Leuven, Belgium. \& Leuven.AI, KU Leuven Institute for AI, B-3000 Leuven, Belgium.}
\email{nick.vannieuwenhoven@kuleuven.be}

\thanks{ORCiDs: 0000-0002-8142-6382 (A. Oneto); 0000-0001-5692-4163 (N. Vannieuwenhoven)}
\thanks{The paper has been accepted for publication on {\it Mathematics of Computation}.}


\begin{document}

\begin{abstract}
A Hadamard--Hitchcock decomposition of a multidimensional array is a decomposition that expresses the latter as a Hadamard product of several tensor rank decompositions. Such decompositions can encode probability distributions that arise from statistical graphical models associated to complete bipartite graphs with one layer of observed random variables and one layer of hidden ones, usually called restricted Boltzmann machines. We establish generic identifiability of Hadamard--Hitchcock decompositions by exploiting the reshaped Kruskal criterion for tensor rank decompositions. An algorithm leveraging existing decomposition algorithms for tensor rank decomposition is introduced for computing a Hadamard--Hitchcock decomposition. Numerical experiments illustrate its computational performance and numerical accuracy {in a noiseless setting}.
\end{abstract}

\maketitle

\section{Introduction}

\subsection*{A motivation from statistics}
{Our interest for the Hadamard--Hitchcock decomposition (HHD) of a multidimensional array arose because of its interesting connection to a classic family of \emph{graphical models} called \textit{discrete restricted Boltzmann machines} (RBM) \cite{smolensky1986information,fischer2012introduction,Mon16}.}
Recall that graphical models are statistical models in which the conditional dependence structure among random variables is modelled by a graph; see, e.g., \cite{lauritzen1996graphical,Bishop2006,sullivant2018algebraic}.
{In an RBM,} the dependence graph is complete and bipartite:
	\begin{equation}
	\begin{tikzpicture}[scale=.9,baseline=(current bounding box.center)]
	\tikzset{cir/.style={circle,draw=black!80,thick,minimum size=0.8cm},y=0.6cm,font=\sffamily}

	\begin{scope}[rotate=90]
	\node[cir] (b1) at (0,0) {$Y_1$};
	\node[cir] (b2) at (0,-1*2) {$Y_2$};
	\node[cir] (b3) at (0,-2*2) {$Y_3$};
	\node[cir] (b4) at (0,-4*2) {$Y_m$};
	\node[cir] (c1) at (-2,1) {$X_1$};
	\node[cir] (c2) at (-2,-1*2+1) {$X_2$};
	\node[cir] (c3) at (-2,-2*2+1) {$X_3$};
	\node[cir] (c4) at (-2,-4*2-1) {$X_d$};

	\foreach \i in {-2,0} {
		\draw (\i,-3*2) node {$\cdots$};
	}

	\foreach \cnto in {1,2,3,4} {
		\foreach \cntt in {1,2,3,4} {
			\draw [-] (b\cnto.south)--(c\cntt.north);
		}
	}
	\end{scope}
	\end{tikzpicture}
	\label{fig:RBM}
\end{equation}
The discrete random variables $X_1,\dots,X_d$ are assumed to be \emph{observable} and {the variables} $Y_1,\ldots,Y_m$ are \emph{hidden}, meaning that the latter are not directly measurable. The observable variables are assumed to be mutually independent conditional on the state of the hidden variables and, vice versa, the hidden variables are assumed to be independent conditional on the state the observable variables.
The probability distribution over the observable variables is given by marginalization over the hidden variables.
As in multi-layer artificial neural networks that chain multiple perceptrons, RBMs can be organized in layers to form practical machine learning models called deep belief networks \cite{hinton2006fast}.

An RBM can be interpreted as a \emph{product of experts}, in particular, a product of \emph{mixtures of independence models}. 
The \emph{independence model} on observed discrete random variables $X_1,\ldots,X_d$ is the set of product distributions, namely, probability distributions whose density functions on $X_1 \times\cdots\times X_d$ are {\it $d$-arrays}, or {\it tensors of order $d$}, of the form $\tensor{P} = \mathbf{p}_1 \otimes \cdots \otimes \mathbf{p}_d$, where $\mathbf{p}_i$ is a probability density function (pdf) on $X_i$ and $\otimes$ denotes the {\it tensor product} \cite{Greub1978}. {In general, any tensor product of nonzero vectors, 
\(
\vect{a}_1 \otimes\dots\otimes \vect{a}_d,
\)
is called a \emph{rank-$1$ tensor}. The pdf of an independence model is thus a rank-$1$ tensor.}
An \emph{$r$-mixture of independence models} is the set of probability distributions that are convex combinations of product distributions, namely, probability distributions whose density functions are
\begin{equation}\label{eqn_cpd}
	\tensor{P}^{(k)} = \sum_{i=1}^{r_k} \pi_i \tensor{P}_{i}^{(k)},
	\quad\text{where }
	\tensor{P}_{i}^{(k)} :=
	\mathbf{p}_{1,i}^{(k)} \otimes \mathbf{p}_{2,i}^{(k)} \otimes\cdots\otimes \mathbf{p}_{d,i}^{(k)}
\end{equation}
and $\sum_{i=1}^{r_k} \pi_i = 1$ with $\pi_i \geq 0$.
This corresponds to the marginalization of the graphical model as in \cref{fig:RBM} with only $m=1$ hidden variable, say $Y$ with $r_k$ states, in which the observed variables are independent conditional on $Y$. More generally, the RBM model is the set of marginal distributions over the visible random variables $X_1,\ldots,X_d$ of probability distributions corresponding to the graphical model as in \cref{fig:RBM}. In this way, the corresponding pdfs are an \emph{elementwise}, \emph{Hadamard}, or \emph{Schur product} of mixtures of independence models, i.e., 
\begin{equation}\label{eqn_basic_hhd}
	\tensor{P} = \lambda \cdot \tensor{P}^{(1)} \hadamard \tensor{P}^{(2)} \hadamard\cdots\hadamard \tensor{P}^{(m)},
\end{equation}
where $\tensor{P}^{(k)}$ is an $r_k$-mixture over $X_1\times\cdots\times X_d$ (where $r_k$ is the number of states of the hidden variable $Y_i$), $\hadamard$ denotes the Hadamard product of tensors and $\lambda$ is a normalization constant so that $\tensor{P} \in \RR^{n_1\times\cdots\times n_d}$ is a pdf. 

For a more elaborate description {of RBMs}, see \cite[Proposition 6]{MM15} or \cite[Section 3]{Mon16}.

\subsection*{{Applied algebraic geometry setup}}
There is a renewed interest in the mathematical properties of RBMs and other graphical models from the perspective of algebraic statistics and applied algebraic geometry; see, e.g., \cite{CMS10,CTY10,MM15,MM17,FOW17,SM18,sullivant2018algebraic}. The reason for this is that the set of probability distributions that can be realized by an RBM with a fixed graph corresponds to a \emph{semi-algebraic subvariety} of a space of \emph{tensors}, i.e., {a set of tensors that satisfy a system of polynomial equalities and inequalities.}

{In this article, we continue along this line of thought, but we will leave the statistical setting and remove} all of the restrictions {on the model} that originate from it: (i) the constraints that $\pi$ and the $\mathbf{p}_{ij}$ in \cref{eqn_cpd} lie on probability simplices and (ii) that $\tensor{P}^{(k)}$ represents a pdf. {Hence, we regard \cref{eqn_basic_hhd} as a polynomial parametrization of a constructible subset of tensors.} One can even relax $\RR$ to an arbitrary field, though we will deal only with fields of characteristic zero. {The motivation for relaxing these constraints is threefold. First, it opens the door to additional tools from applied algebraic geometry, as the semi-algebraic sets describing RBMs will now become \emph{algebraic varieties}, whose theory is more developed \cite{Harris1992}. Second, we will show that the statistical restrictions are not a strict requirement for the identifiability and decomposition algorithms of such models, even though they can be helpful during numerical computations. Third, since our algorithmic and identifiability results hold outside a complex algebraic subvariety of the aforementioned constructible set, they also hold for the generic probability distributions representable by RBMs. The reason is that the latter form a Euclidean open subset of the constructible set, as they are obtained by restricting both the domain and codomain of the polynomial parametrization to full-dimensional subsets.}

The decomposition of a tensor into a sum of {rank-$1$} tensors {as in \cref{eqn_cpd} was introduced by} Hitchcock \cite{Hit27}. Its variants have been rediscovered several times in the literature under various names: \emph{polyadic decomposition}, \emph{canonical decomposition}, \emph{parallel factor analysis}, \textit{CANDECOMP/PARAFAC}, \textit{CP decomposition}, \textit{canonical polyadic decomposition (CPD)}, \emph{tensor rank decomposition}, and \emph{Hitchcock decomposition}. 
Throughout this paper, we will refer to a Hitchcock decomposition \cref{eqn_cpd} of minimum length $r_k$ as a \emph{(rank-$r_k$) CPD} and we call $r_k = \mathrm{rank}(\tensor{P}^{(k)})$ the \emph{rank} of $\tensor{P}^{(k)}$. If we {do not know whether} the decomposition in \cref{eqn_cpd} is of minimum length, then we refer to it as a Hitchcock decomposition.

CPDs are one of the classic tensor decompositions studied in multilinear algebra, and they appear in a wide variety of applications \cite{KB09,Lan12,PFS2016,BC2019}. By contrast, the decomposition \cref{eqn_basic_hhd}, while {connected to the important applications of RBMs,} seems, at present, comparatively poorly understood{, see \cite{Mon16,FOW17,SM18}}. Let us first define and name it appropriately.\footnote{A note about terminology is needed here. Mathematically, a tensor is not simply a box of numbers like the $d$-arrays $\RR^{n_1\times\dots\times n_d}$, but instead it (i) satisfies transformation laws, (ii) represents multilinear maps, and (iii) lives in a tensor product of vector spaces \cite{Lim2021}. An \ref{eqn_hhd} decidedly does not satisfy these properties because of its dependency on the choice of basis. For this reason, it would be imprecise to refer to them as \emph{tensor} decompositions, and one should more correctly refer to them as a decomposition of a tensor in coordinates. 
We stress that in many applications, a natural basis is anyhow forced upon us and that other decompositions of tensors in coordinates, such as nonnegative matrix or tensor decompositions, have in recent years proven their mathematical intricacy and utility in applications \cite{Gillis2020,NMTFBook2009}.}

\begin{definition}\label{def_hhd}
Let $\Bbbk^{n_1 \times \cdots \times n_d} := \Bbbk^{n_1} \otimes \cdots \otimes \Bbbk^{n_d}$ be the vector space of order-$d$ tensors, with coefficients in the field $\Bbbk$, and let $\tensor{T} \in \Bbbk^{n_1 \times \cdots \times n_d}$ be a tensor. For any $\br = (r_1,\ldots,r_m) \in \bbN^m$, an \textit{$\br$-Hadamard--Hitchcock decomposition} (\ref{eqn_hhd}) of $\tensor{T}$ is, if it exists, {a Hadamard product of Hitchcock decompositions:}
	\begin{equation}\tag{$\br$-HHD}\label{eqn_hhd}
		\tensor{T} = \bhadamard{k=1}{m} \sum_{i_k=1}^{r_k} \tensor{A}_{k i_k},
	\end{equation}
where $\tensor{A}_{k i_k}$ is a rank-$1$ tensor for every $k \in [m]$ and $i_k \in [r_k]$ and $\hadamard$ denotes the \textit{Hadamard product}.
We refer to the Hitchcock decomposition of length $r_k$ appearing as $k$th factor in the above Hadamard product as the \emph{$k$th Hadamard factor}.
\end{definition}

In this paper, we address the inverse problem of identifying all rank-$1$ terms that appear in an \ref{eqn_hhd} of a tensor. Specifically, we establish its uniqueness properties {when $d \ge 3$} and propose a direct algorithm for solving it under certain conditions.

{\begin{remark}
In the case of matrices ($d=2$), the foregoing decomposition was also recently considered by Ciaperoni, Gionis, and Mannila \cite{CGM2024}, who called it a Hadamard decomposition. They primarily studied the case with only $2$ hadamard factors ($m=2$), proposing optimization algorithms and observing the non-uniqueness of the \ref{eqn_hhd} model when $d=2$. The question of uniqueness is ultimately left unresolved, concluding that ``there is no obvious way of characterizing the nonuniqueness'' in \cite[Section 6]{CGM2024}. Sufficient conditions for uniqueness were subsequently proposed in \cite{DJL2025}. For $d\ge3$, we fully characterize the generic uniqueness properties for tensors admitting a low-rank \ref{eqn_hhd}. 
\end{remark}}

\subsection*{Outline and main contributions}

In \cref{sec:preliminaries}, we recall standard concepts of tensor decompositions such as the notion of $R$-\textit{identifiability}: the property of admitting a unique rank-$R$ CPD up to symmetries. This property can be certified with a classic criterion due to Kruskal \cite{kruskal1977three}, which will play a main role in our identifiability results for \ref{eqn_hhd}s. We also recall known results on the algebraic geometry of sets of \ref{eqn_hhd}s, which again relate to the question of identifiability.

The other crucial ingredient in our approach to identifiability of \ref{eqn_hhd}s is developed thoroughly in \cref{sec_rk1permpres}. It is the new concept of a \emph{rank-$1$ permutation}, which we believe could be of independent interest. In particular, we present the efficient \cref{alg_rank1_permutation} that finds, if possible, a permutation of a tuple of $r_1\cdots r_m$ numbers into a rank-$1$ tensor in $\Bbbk^{r_1 \times \cdots \times r_m}$.

In \cref{sec:identifiability}, we prove \emph{generic identifiability} of \ref{eqn_hhd}s when $R=r_1\cdots r_m$ is sufficiently small. That is, identifiability holds on a Zariski open subset of \ref{eqn_hhd}s. The main idea is to expand an \ref{eqn_hhd} as a Hitchcock decomposition by using the distributivity between the sum and the Hadamard product.\footnote{In the statistical context, this means that we (temporarily) interpret the hidden variables $Y_1,\ldots,Y_m$ with $Y_k\in[r_k]$ as a unique hidden variable with $R = r_1\cdots r_m$ states.} If the (reshaped) Kruskal criterion guarantees that such a Hitchcock decomposition is $R$-identifiable, then we show, based on the results of \cref{sec_rk1permpres}, that the rank-$1$ terms in the \ref{eqn_hhd} can be recovered {mathematically}, up to certain group actions. {Since $R$-identifiability of CPDs holds only for tensors of order three and higher, it is henceforth assumed that $d \ge 3$.}

The main results, \cref{thm:main_identifiabilityHHD,prop_hhdkruskal}, are summarized as follows.

\begin{theorem*}[Identifiability of generic \ref{eqn_hhd}s]
	Let $\br\in\NN^m$ and $R=r_1\cdots r_m$.
	Let $\tensor{T}$ admit a generic \ref{eqn_hhd}.
	If $\tensor{T}$ is $R$-identifiable, then the rank-$1$ tensors $\tensor{A}_{k i_k}$ in the \ref{eqn_hhd} are uniquely determined up to the following symmetries:
	\begin{enumerate}[nolistsep]
	 \item[(i)] permuting Hadamard factors of equal length, 
	 \item[(ii)] permuting the summands of the Hitchcock decomposition in any Hadamard factor, and
	 \item[(iii)] Hadamard multiplication by rank-$1$ tensors with nonzero entries. 
	 \end{enumerate}
	 Moreover, if $R$ falls within the range of the reshaped Kruskal criterion \cite[Theorem 4.6]{COV17}, then the generic $\tensor{T}$ admitting a decomposition as in \cref{eqn_hhd} is $R$-identifiable and satisfies the previous statement.
\end{theorem*}
The third symmetry follows from the observation that $\rk(\tensor{A}\hadamard\tensor{T}) \leq \rk(\tensor{T})$ if $\tensor{A}$ has rank $1$. Hence, if $\tensor{T} = \tensor{T_1}\hadamard\cdots\hadamard\tensor{T}_m$ is an \ref{eqn_hhd}, then for any rank-$1$ tensor $\tensor{A}$ with nonzero entries, we get a new \ref{eqn_hhd} by Hadamard-multiplying $\tensor{T_i}$ by $\tensor{A}$ and $\tensor{T_j}$ by the Hadamard-inverse~$\tensor{A}^\ominus$ for any $i \ne j$.

As a corollary of the main result, \cref{cor:dim_RBM} will state the dimension of the variety of tensors admitting an \ref{eqn_hhd}, showing that it matches the expected dimension, according to \cite{BCK16}.

In \cref{sec:algorithm}, we propose \cref{alg_hhd} for computing the unique {decomposition of a generic low-rank \ref{eqn_hhd} (noiseless) tensor}, based on the proof strategy of the aforementioned main result. Essentially, it computes a rank-$R$ CPD of the tensor, then computes a rank-$1$ permutation using \cref{alg_rank1_permutation}, and finally obtains the rank-$1$ tensors $\tensor{A}_{k i_k}$ through the computation of several rank-$1$ CPDs. The cost of computing the rank-$R$ CPD generally dominates the asymptotic time complexity of this algorithm.
{While \cref{alg_hhd} is conceived as an algorithm suitable for implementation in floating-point arithmetic, we stress that it is designed only for tensors that admit an exact \cref{eqn_hhd} up to numerical round-off errors. It is not an optimization-based algorithm that can be applied even under model violations.}

\Cref{sec:experiments} presents numerical experiments illustrating the computational performance and accuracy of \cref{alg_hhd} on data admitting an exact \ref{eqn_hhd}, up to numerical round-off errors. {We also investigate to what extent slight model violations can be tolerated by the algorithm.}

Finally, the main conclusions and a future outlook are presented in \cref{sec_conclusions}.

\subsection*{Notation}
We fix some notation that will be used throughout this paper.
	Let $\Bbbk$ be a field of characteristic zero. We use the shorthand notations
\[
[s] := \{1,\ldots,s\} \quad\text{and}\quad
[\mathbf{s}] := [s_1]\times\cdots\times [s_m]
\]
for $s\in\NN$ and $\mathbf{s} \in \NN^m$, respectively.
We reserve the symbols
\[
	\bn := (n_1,\ldots,n_d),\quad
	N := \prod_{i=1}^d n_i, \quad
	\br := (r_1,\ldots,r_m), \quad\text{and}\quad
	R := \prod_{i=1}^m r_i
\]
{for, respectively,} the shape of the tensor to decompose, the dimension of the ambient space $\Bbbk^{n_1\times\cdots\times n_d}$, the rank parameters of the \ref{eqn_hhd}, and the length of the induced Hitchcock decomposition (see \cref{ssec_hhd_cpd} for this interpretation).
We assume throughout the paper that all $r_k \ge 2$, which will be important for the correctness of \cref{alg_rank1_permutation}; as explained in \cref{sec:preliminaries}, if some $r_k=1$, then a strictly simpler HHD can be constructed, so this assumption does not limit generality.

The elementwise multiplication of two tensors $\tensor{R}, \tensor{S} \in \Bbbk^{n_1\times\cdots\times n_d}$ is {defined as}
\[
	\tensor{T} = \tensor{R} \hadamard \tensor{S}\text{, where }\tensor{T}_{i_1,\ldots,i_d}=\tensor{R}_{i_1,\ldots,i_d}\cdot\tensor{S}_{i_1,\ldots,i_d}\text{ for all }(i_1,\ldots,i_d)\in[\bn].
\]
This is known as the \emph{elementwise}, \emph{Hadamard}, or \emph{Schur} product. To streamline some expressions involving Hadamard products of $\tensor{T}\in\Bbbk^{n_1\times\cdots\times n_d}$, we let $\tensor{T}^\ominus$ denote the elementwise inverse.

Whenever we refer to a \textit{generic} property on an algebraic variety $\Var{V}\subset\Bbbk^N$, we mean that the property holds everywhere outside of a strict Zariski closed subset $\Var{Z}\subsetneq\Var{V}$. In particular, this implies that if an element $x \in \Var{V}$ is sampled randomly from a probability distribution absolutely continuous with respect to the Lebesgue measure on $\Var{V}$, then the probability that $x \in\Var{Z}$ is zero. Clearly, if a property $A$ holds generically, and $B$ holds generically as well, then $A \wedge B$ holds generically as well. This fact will be used often without explicitly mentioning it.

Tensors will be typeset in a calligraphic font (e.g., $\tensor{T}$). We consistently use the letter ``a'' to refer to quantities related to rank-$1$ tensors, such as rank-$1$ tensors ($\tensor{A}$), vectors that can be reshaped into a rank-$1$ tensor ($\ba$), the set of such vectors ($\Var{A}$), and vectors whose tensor product is taken $(\ba_1\otimes\cdots\otimes\ba_d$). The rank-$1$ ``unit'' tensor filled with ones is denoted by $\mathds{1} = \mathbf{1}\otimes\cdots\otimes\mathbf{1}$, where ${\bf1}$ is the vector of ones.

\subsection*{Acknowledgments}
{We thank the reviewers for their suggestions and questions that led to significant improvements of this paper, notably resulting in an improved introduction and the addition of \cref{sec_sub_advanced_version,sec_sub_comparisonpba,sec_sub_noise_tolerance,sec_sub_rbm}. We thank the editor for handling the process.}

We thank Benjamin Lovitz for pointing us to Westwick's result.
AO thanks Guido Montufar for useful discussions about RBMs and the corresponding algebraic statistical models.
NV thanks Piotr Zwiernik for an interesting discussion on the interpretation of HHDs {and Dirk Nuyens for a helpful discussion on kernel density estimation and quasi-Monte Carlo sampling}.

This project was initiated when NV visited the TensorDec Lab of the Department of Mathematics of the Universit\`a di Trento from October 23, 2021 to November 20, 2021. The collaboration continued during the \emph{Algebraic Geometry with Applications to Tensors and Secants} (AGATES) thematic semester in Warsaw, Poland {during the Fall of 2022. The revision was concluded when NV visited the TensorDec Lab from May 3, 2025 to July 3, 2025.}

\subsection*{Funding}
AO is member of INdAM-GNSAGA (Italy). {AO has been partially supported by MUR (Ministero dell'Università e della Ricerca, Italy) through the PRIN2022SC project Prot. 2022NBN7TL (CUP: E53C24002320006) “Applied Algebraic Geometry of Tensors”.}

NV was partially supported by KU Leuven Interne Fondsen grants STG/19/002 {and C16/21/002}.

AO and NV were partially supported by the Thematic Research Programme ``Tensors: geometry, complexity and quantum entanglement,'' University of Warsaw, Excellence Initiative -- Research University and the Simons Foundation Award No.~663281 granted to the Institute of Mathematics of the Polish Academy of Sciences for the years 2021--2023. 

\section{Background}\label{sec:preliminaries}
We summarize relevant aspects of the geometry of the spaces of CPDs and HHDs from the literature in this section.

\subsection{CPDs}\label{ssec:identifiability_cpd}
Let $\bn\in\NN^d$.
The space of rank-$1$ tensors is denoted by
\[
\Var{S}_\bn = \{\tensor{A} \in \Bbbk^{n_1\times\dots\times n_d} ~:~ {\rm rank}(\tensor{A}) = 1\}.
\]
The ambient space is usually clear from the context, in which case we drop the subscript $\bn$.
$\Var{S}_\bn$ is a smooth manifold. Its closure $\overline{\Var{S}}=\Var{S}\cup\{0\}$ is an irreducible algebraic variety, called the \textit{Segre variety}, corresponding to the Zariski closure of the  ${\rm GL}(\Bbbk^{n_1})\times\cdots\times {\rm GL}(\Bbbk^{n_d})$-orbit of the unit tensor $\mathds{1}$. We will not distinguish in notation between $\Var{S}$ and $\overline{\Var{S}}$, using the former uniformly.

CPDs have been successfully approached by means of algebraic geometry through \textit{secant varieties}. The \textit{$r$th secant variety $\sigma_r(\Var{S})$} is the (Zariski) closure of the set of tensors of rank at most $r$:
\[
	 \sigma_r(\Var{S}) = \overline{\sigma^\circ_r(\Var{S})},
	 \quad \text{where } \sigma^\circ_r(\Var{S}) = \{\tensor{T} \in \Bbbk^{n_1\times\dots\times n_d} ~:~ {\rm rank}(\tensor{T}) \leq r\}.
\]
For completeness, recall that for $r \geq 2$, the set of {higher-order tensors of rank at most $r$ is in general \textit{not} closed (this does not happen in the case of matrices $d = 2$)}. Indeed, there are examples of tensors {whose rank is strictly larger than $r$ which are nevertheless the limit of a sequence of} rank-$r$ tensors; see, e.g., \cite{dSL2008}.

Consider the following addition map for a fixed $r\in\NN$:
\begin{equation*}
\begin{array}{c c c c}
	+_r : & \Var{S}^{\times r} & \longrightarrow & \sigma_r^\circ(\Var{S}) \subset \Bbbk^{n_1\times\dots\times n_d}, \\ 
	& \left(\tensor{A}_1,\dots,\tensor{A}_r\right)& \longmapsto & \tensor{A}_1+\dots+\tensor{A}_r.
\end{array}
\end{equation*}
Clearly, if $\tensor{T} \in \sigma^\circ_r(\Var{S})$ and $\left(\tensor{A}_1,\dots,\tensor{A}_r\right) \in +_r^{-1}(\tensor{T})$, then the same is also true for any permutation of the $r$-tuple of rank-$1$ tensors. 
Let $\mathfrak{S}_r$ be the group of permutations of $[r]$. We denote by
\[
	\Var{S}^{(r)} := \Var{S}^{\times r} / \mathfrak{S}_r
\]
the symmetric power of the Segre variety, i.e., the set of $r$-tuples of rank-$1$ tensors up to permutations. Given $\tensor{T} \in \sigma^\circ_r(\Var{S})$, we consider the \textit{set of CPDs} of $\tensor{T}$ as the set $+_r^{-1}(\tensor{T})$ up to permutation:
\[
	\calD_r(\tensor{T}) := +_r^{-1}(\tensor{T}) / \mathfrak{S}_r \subset \Var{S}^{(r)}.
\]
Geometrically, if $\{\tensor{T}_1,\ldots,\tensor{T}_r\} \in \calD_r(\tensor{T})$, then $\tensor{T}$ lies in the linear span of the $\tensor{T}_i$'s.

\begin{definition}
	A rank-$r$ tensor $\tensor{T}$ is $r$-\textit{identifiable} if $\calD_r(\tensor{T})$ is a singleton. 
\end{definition}

There is a rich literature studying the dimension and cardinality of $\calD_r(\tensor{T})$; see, e.g., \cite{abo2009induction,CO2012,Gesmundo2013,COV2014,bernardi2018hitchhiker,MM2022,blomenhofer2023nondefectivity}. Subject to a few understood exceptions, the conjecture is that if {the dimension of $\sigma_r(\Var{S})$ is strictly smaller than the dimension $n_1\cdots n_d$ of the ambient space}, then $\calD_r(\tensor{T})$ is a singleton for generic $\tensor{T}\in\sigma_r(\Var{S})$.

Since one aims to recover the parameters of a specific rank-$r$ tensor $\tensor{T}$ uniquely, it is important to have tools to certify identifiability. In pioneering work, Kruskal \cite{kruskal1977three} provided a criterion to guarantee identifiability of order-$3$ tensors which is nowadays still popular and widely used. The method is based on the notion of \textit{Kruskal rank}, which is recalled next.

Let $K = \{k_1,\ldots,k_m\} \subset [d]$. In $\Bbbk^{n_1\times\dots\times n_d}$, consider the set of points $A = \{\tensor{A}_1,\ldots,\tensor{A}_r\} \in \Var{S}^{(r)}$. Let $p_K$ be the forgetful projection {that takes $\vect{a}_1\otimes\cdots\otimes \vect{a}_d \in \Var{S} \subset \Bbbk^{n_1\times\dots\times n_d}$ to the $1$-dimensional linear space $\langle \vect{a}_{k_1}\otimes\cdots\otimes \vect{a}_{k_m} \rangle \subset \Bbbk^{n_{k_1}\times\dots\times n_{k_m}}$}. Then, the Kruskal rank of $p_K(A)$, denoted by $\mathrm{krank}_K(A)$, is defined as the maximum value $j \in \mathbb{N}$ such that for every set $S \subset p_K(A)$ of cardinality $\sharp S=j$, the dimension of the linear span of $S$ is $j$, i.e., $\dim\mathrm{span}(S) = j$.

Kruskal \cite{kruskal1977three} proved a sufficient condition for the identifiability of order-$3$ tensors based on the Kruskal ranks of the projections $p_{\{i\}}(A)$ for $i=1,2,3$. For higher-order tensors $\tensor{T}$, the same criterion can be applied by interpreting $\tensor{T}$ as an order-$3$ tensor through \textit{reshaping} \cite{AMR2009,COV17}.

\begin{lemma}[{\cite[Theorem 4.6]{COV17}}\footnote{The necessary condition that all Kruskal ranks are at least $2$ is missing in the statement of \cite[Theorem~4.6]{COV17}.}]\label{lem_reshaped_kruskal}
Let $\tensor{T}\in\Bbbk^{n_1\times\cdots\times n_d}$, {$d\ge3$,} and $A = \{\tensor{A}_1,\ldots,\tensor{A}_R\} \in \calD_R(\tensor{T})$. If there exists a partition $I \sqcup J \sqcup K = [d]$ such that
\[
 R \le \frac{1}{2}(\mathrm{krank}_I(A) + \mathrm{krank}_J(A) + \mathrm{krank}_K(A) - 2)
\]
and $\mathrm{krank}_X(A) \ge 2$ for all $X \in \{I,J,K\}$, then $\tensor{T}$ is $R$-identifiable. 

Furthermore, this criterion is effective if $\Pi_I \ge \Pi_J \ge \Pi_K \ge 2$, where $\Pi_X := \prod_{x\in X} n_x$, and
\begin{equation}\label{eqn_kruskal_range}
 R \le \Pi_I + \min\left\{\frac{1}{2}\delta, \delta\right\}, \quad\text{where } \delta:= \Pi_J + \Pi_K - \Pi_I - 2;
\end{equation}
that is, a generic rank-$R$ tensor in $\Bbbk^{n_1\times\cdots\times n_d}$ is $R$-identifiable under these assumptions.
\end{lemma}

We will refer to the value on the right hand side of \cref{eqn_kruskal_range} as the range of applicability of the reshaped Kruskal criterion with partition $I\sqcup J\sqcup K=[d]$. A heuristic for determining partitions with a large range of applicability was discussed in \cite[Section 4.2]{COV17}. We refer to the largest range of applicability over all partitions of $[d]$ simply as the \emph{range of applicability} of \cref{lem_reshaped_kruskal}.

\subsection{HHDs}\label{sec:preliminaries_hhd}
The space of tensors admitting an \ref{eqn_hhd} is the image of the map 
\begin{equation}\label{eq:HHD_map}
	\begin{array}{c c c c}
		\hadamard_\br : & \sigma^\circ_{r_1}(\Var{S})\times\cdots\times \sigma^\circ_{r_m}(\Var{S}) & \longrightarrow &  \Bbbk^{n_1\times\dots\times n_d}, \\ 
		& \left(\tensor{T}_1,\ldots,\tensor{T}_m\right)& \longmapsto & \tensor{T}_1 \hadamard \cdots \hadamard \tensor{T}_m.
	\end{array}
\end{equation}
We denote it as 
\[
	\sigma_{\br}^\circ(\Var{S}) := {\rm im}(\hadamard_\br) = \{\tensor{T} \in \Bbbk^{n_1\times\dots\times n_d}~:~\tensor{T} \text{ admits an \ref{eqn_hhd}}\}.
\]
Note that the Hadamard product of rank-$1$ tensors is of rank bounded by $1$ because
\[
	(\ba_1\otimes\cdots\otimes \ba_d) \hadamard (\bb_1\otimes\cdots\otimes \bb_d) = (\ba_1\hadamard \bb_1)\otimes\cdots\otimes(\ba_d\hadamard \bb_d),
\]
for all $\ba_k, \bb_k \in \Bbbk^{n_k}$, so that $\rk(\tensor{T} \hadamard \tensor{A}) \leq \rk(\tensor{T})$ if $\rk(\tensor{A}) = 1$. This implies that we can assume that $r_k \geq 2$ for all $k \in [m]$. Indeed,
if for example $r_m = 1$,
then $\hadamard_\br(\tensor{T}_1,\ldots,\tensor{T}_m) = \hadamard_{(r_1,\ldots,r_{m-1})}(\tensor{T}_1\hadamard\tensor{T}_m,\tensor{T}_2,\ldots,\tensor{T}_{m-1})$. Consequently, by dropping all components of $\br$ equal to $1$, we obtain a strictly simpler decomposition.

After introducing \ref{eqn_hhd}s, it is natural to ask if there exists an $\br\in\NN^m$ such that \textit{every} tensor in $\Bbbk^{n_1\times\dots\times n_d}$ admits an \ref{eqn_hhd}, i.e., if $\hadamard_\br$ is surjective for some $\br$. For $m=1$, so $\br=(R)$, the answer is affirmative because rank-$1$ tensors span $\Bbbk^{n_1\times\dots\times n_d}$.
For $m \geq 2$, if $r_k \geq 2$ for all $k\in [m]$, then every $\tensor{T} \in \Bbbk^{n_1\times\cdots\times n_d}$ admits an \ref{eqn_hhd} if the number of Hadamard factors $m$ is sufficiently large compared to the tensor's order, i.e., $m \gg d$. Indeed, if $\tensor{T}_{:,i_2,\ldots,i_d} = ( \tensor{T}_{1,i_2,\ldots,i_d}, \dots, \tensor{T}_{n_1,i_2,\dots,i_d} ) \in \Bbbk^{n_1}$, then,
\[
	\tensor{T} = \bighadamard_{i_2,\ldots,i_d=1}^{n_2,\ldots,n_d} \left((\tensor{T}_{:,i_2,\ldots,i_d} - \mathbf{1}) \otimes \mathbf{e}_{i_2} \otimes \cdots \otimes \mathbf{e}_{i_d} + \mathds{1}\right),
\]
where $\mathbf{e}_i$ is zero everywhere except at position $i$ where it is one.

The question of identifiability of \ref{eqn_hhd}s is closely related to the dimension of the space of tensors admitting an \ref{eqn_hhd}.
The (Zariski closure of the) space of tensors admitting an \ref{eqn_hhd} can be constructed as a \textit{Hadamard product of secant varieties}, as introduced in \cite{CMS10}.

Since the tensor rank is invariant under scalar multiplication, in the classic algebraic geometry literature, secant varieties are considered \textit{projectively}. Let $\bbP V$ be the projective space of nonzero vectors in the vector space $V$ up to nonzero scalar multiplication. {If $\vect{v}\in V$, we denote by $[\vect{v}] \in \bbP V$ the corresponding equivalence class, i.e., $[\vect{v}]=[\vect{w}]$ if and only if $\vect{v}=\lambda \vect{w}$ for some $\lambda \neq 0$}. With a minor abuse of notation, we refer here to the Segre variety $\Var{S}$ and its secant varieties $\sigma_r(\Var{S})$ as \emph{projective} varieties in $\bbP \Bbbk^{n_1\times\cdots\times n_d}$. As projective variety, the dimension of $\Var{S}$ is $\dim\Var{S}=\sum_{i=1}^d (n_i-1)$.

Fixing coordinates on $\Bbbk^N$, let $(x_1:\ldots:x_N)$ be the projective coordinates of $[(x_1,\ldots,{x_N})] \in \bbP\Bbbk^N$. Following \cite{CMS10}, the \textit{Hadamard product} {defines a map} 
\[
	\hadamard : \bbP\Bbbk^N \times \bbP\Bbbk^N \dashrightarrow \bbP\Bbbk^N, \quad ((a_1:\ldots:a_N),(b_1:\ldots:b_N)) \mapsto (a_1b_1:\ldots:a_Nb_N),
\]
{which is well defined only on a Zariski open subset of $\bbP\Bbbk^N \times \bbP\Bbbk^N$, hence the dashed arrow. In particular, if we consider two nonzero vectors $\vect{a}, \vect{b}\in \Bbbk^N$ whose Hadamard product is nevertheless equal to zero, then $\vect{a}\hadamard\vect{b} = 0$ does not define a point in projective space}. In general, given two projective varieties $X,Y \subset \bbP \Bbbk^N$, their \textit{Hadamard product}, denoted by $X \hadamard Y$, is defined as the Zariski closure of the image of $X \times Y$ under the Hadamard product, i.e.,
\[
	X \hadamard Y = \overline{\{x \hadamard y ~:~ x \in X,\, y \in Y,\, x \hadamard y \text{ is well defined}\}}.
\]
Then, we have that 
\[
	\sigma_{\br}(\Var{S}) = \overline{\sigma^\circ_\br(\Var{S})} = \bhadamard{k=1}{m} \sigma_{r_k}(\Var{S}).
\]
As mentioned above, the Hadamard product of rank-$1$ tensors still has rank $1$, i.e., $\Var{S} \hadamard \Var{S} = \Var{S}$. Geometrically, this is a particular case of the general fact that any \textit{toric variety}, embedded in the sense of \cite{sturmfels1996equations}, is invariant by taking Hadamard product with itself; see e.g. \cite{FOW17}. 

The last observation also gives a heuristic for the fact that, given $X,Y \subset \bbP \Bbbk^N$, the dimension of $X \hadamard Y$ depends on the largest toric variety $H$ contained in $X \cap Y$. In particular, if $H$ is the largest torus acting on both $X$ and $Y$, then
\begin{equation}\label{eq:expdim_hadamard}
	\dim(X \hadamard Y) \leq \min\{\dim(X) + \dim(Y) - \dim(H), N-1\}
\end{equation}
by \cite[Proposition~5.4]{BCK16}.
The right hand side is called the \textit{expected dimension}. There are examples of \textit{defective} Hadamard products of projective varieties for which this inequality is strict; {see \cite[Example 5.5]{BCK16}}.

For Hadamard products of secant varieties of Segre varieties, we find the Segre variety $\Var{S}$, which is a toric variety, contained in all of its secant varieties. Therefore, applying \cref{eq:expdim_hadamard}, we get the following inequality among projective dimensions:
\begin{equation}\label{eqn:upperbound_dim_HHD}
	\dim \sigma_\br(\Var{S}) \leq
	{\exp}.\dim \sigma_\br(\Var{S}) := \min\left\{ \sum_{j=1}^m \dim\sigma_{r_j}(\Var{S}) - (m-1)\dim \Var{S}, N-1 \right\};
\end{equation}
see also \cite[Proposition 14]{MM15}.
The right-hand side of \cref{eqn:upperbound_dim_HHD} is called the \textit{expected dimension} of $\sigma_\br(\Var{S})$.
It is expected that the expected dimension is the actual dimension of $\sigma_\br(\Var{S})$, see \cite{MM15,FOW17}. 
For more on Hadamard products of algebraic varieties we refer to \cite{CTY10,BCK16,FOW17,CCFL20,BC22,AB+22}.

\subsection{From HHD to CPD and back}\label{ssec_hhd_cpd}
An \ref{eqn_hhd} with $\br=(R)$ is evidently a rank-$R$ CPD. More generally, by distributivity of Hadamard multiplication with the sum, it follows that every \ref{eqn_hhd} of a tensor $\tensor{T} \in \Bbbk^{n_1\times\cdots\times n_d}$ induces a special Hitchcock decomposition of $\tensor{T}$:
\begin{equation*} 
	\tensor{T}
	= \bhadamard{k=1}{m} \sum_{i_k=1}^{r_k} \tensor{A}_{k i_k}
	= \sum_{\bi\in[\br]} \bhadamard{k=1}{m} \tensor{A}_{k i_k}.
\end{equation*}
Note that the induced Hitchcock decomposition on the right side does not necessarily correspond to a CPD of $\tensor{T}$, namely it is not guaranteed in general that $\rk(\tensor{T}) = R$. In fact, we see that
\begin{equation}\label{inclusion}
	\sigma_\br(\Var{S}) \subset \sigma_R(\Var{S}),
\end{equation}
but the inclusion is expected to be of high codimension, implying that the induced Hitchcock decomposition of length $R$ could have a smaller tensor rank.
Indeed, by a parameter count, the dimension of the $R$th secant variety satisfies \cite{Harris1992,bernardi2018hitchhiker}:
\begin{equation}\label{eqn:upperbound_dim_CPD}
\dim \sigma_R(\Var{S}) \le {\exp}.\dim \sigma_R(\Var{S}) := \min\left\{R \cdot (1+\dim\Var{S}) - 1, N-1 \right\}.
\end{equation}
While there are \textit{defective} cases in which the previous inequality is strict, for most secant varieties of the Segre variety the inequality is an equality; see \cite{abo2009induction,COV2014,bernardi2018hitchhiker,blomenhofer2023nondefectivity}.
On the other hand, the expected dimension of $\sigma_\br(\Var{S})$ satisfies \cref{eqn:upperbound_dim_HHD}.
By comparing the two expected dimensions in \cref{eqn:upperbound_dim_HHD,eqn:upperbound_dim_CPD}, it is clear that \cref{inclusion} is expected to be strict.
This implies that \textit{a priori} there is no reason to expect that generic properties of CPDs would also apply to the decompositions induced from an \ref{eqn_hhd}. Nonetheless, for identifiability this favorable result is precisely what the main result of this paper will establish for sufficiently small ranks.

\section{Rank-1 permutations} \label{sec_rk1permpres}

The central hypothesis of this paper is that an \ref{eqn_hhd} can be recovered from its induced Hitchcock decomposition under appropriate assumptions. However, one may anticipate that the indeterminate ordering of the rank-$1$ summands in the induced Hitchcock decomposition poses an obstacle.
We propose \emph{rank-$1$ permutations} as a way out of this conundrum.

\subsection{Preservers and permutations}

We start by defining rank-$1$ permutations and stating some of their useful properties that will be used in the next sections.

\begin{definition}\label{def_admissible}
	{Let $\rho : [R] \to [\br]$ be a bijection. Given $\ba\in\Bbbk^R$, we denote by $\rho(\ba)\in \Bbbk^{r_1\times\cdots\times r_m}$ the tensor whose entries are given by $\rho(\ba)_{i_1,\ldots,i_m} := \ba_{\rho^{-1}(i_1,\ldots,i_m)}$ for all $(i_1,\ldots,i_m)\in[\br]$. We call such a map} a \emph{rank-$1$ permutation} of $\ba \in \Bbbk^R$ if $\rho(\ba) \in \Bbbk^{r_1\times\cdots\times r_m}$ is a rank-$1$ tensor. For any factoring $\br$ of $R$, and bijection $\rho : [R] \to [\br]$, let
	\[
		\Var{A}^\rho = \{\ba \in \Bbbk^R ~:~ \rho \text{ is a rank-$1$ permutation of }\ba\} = \rho^{-1}(\Var{S}_\br)
	\]
be the set of \emph{admissible} vectors for $\rho$.
\end{definition}

Note that $\rho$ induces an isomorphism between $\Bbbk^R$ and $\Bbbk^{r_1\times\cdots\times r_m}$ which shows that $\Var{A}^\rho$ is an irreducible algebraic variety isomorphic to the Segre variety. 

Before continuing, we recognize that the terminology ``rank-$1$ permutation'' may be confused with the \emph{tensor product of permutations}, introduced in \cite[Section 2]{BC2021}.
If we let $\mathfrak{S}([\br])$ denote the set of permutations of the set of multi-indices $[\br] = \{ (i_1,\ldots,i_m) : i_1\in[r_1],\ldots,i_m\in[r_m]\}$, then the tensor product of the permutations $\pi_k \in \mathfrak{S}_{r_k}$, $k\in[m]$, is the permutation $\pi=\pi_1\otimes\cdots\otimes\pi_m \in \mathfrak{S}([\br])$ that acts like
\[
\pi(i_1,\ldots,i_m)
= (\pi_1\otimes\cdots\otimes\pi_m)(i_1,\ldots,i_m)
= (\pi_1(i_1), \ldots, \pi_m(i_m)).
\]
The notation $\pi_1\otimes\cdots\otimes\pi_m$ and name are appropriate in light of the next basic result, stated in \cite[Section 2]{BC2021}. Recall that a \emph{permutation matrix} is a matrix whose columns are a permutation of the identity matrix.

\begin{lemma}\label{lem_tenprodperm}
Let the permutation $\pi_k \in \mathfrak{S}_{r_k}$ be represented by the permutation matrix $\Pi_k \in \mathbb{R}^{r_k \times r_k}$, for $k\in[m]$. Then, the permutation $\pi = \pi_1\otimes\cdots\otimes\pi_m \in \mathfrak{S}([\br])$ is represented by the permutation matrix given by the Kronecker product $\Pi_1 \otimes\cdots\otimes \Pi_m \in \Bbbk^{R \times R}$.
\end{lemma}
\begin{proof}
A matrix is a permutation matrix if and only if it is orthogonal and nonnegative \cite[Lemma 5]{ZP2008}.
The Kronecker product preserves orthogonality and non-negativity. This concludes the proof.
\end{proof}

The tensor product of permutations is closely related to the nonuniqueness of rank-$1$ permutations of an admissible vector.
Indeed, rank-$1$ permutations are not unique; for example, if $\mathbf{a}=(1,1,1,1)$ with $\br=(2,2)$, then all $4!$ permutations are rank-$1$ permutations.
To understand which rank-$1$ permutations exist for an admissible vector, we need to recall the concept of a \emph{rank-$1$ preserver} from Westwick \cite{Westwick1967}. This is a linear map $L : \Bbbk^{r_1\times\cdots\times r_m} \to \Bbbk^{r'_1\times\cdots\times r'_n}$ with the property that it maps rank-$1$ tensors into rank-$1$ tensors.
We are interested in rank-$1$ preservers that are also permutations.

\begin{definition}\label{def:rank1_preserver}
	Given a bijection $\sigma : [\br] \rightarrow [\br']$, we denote the induced isomorphism $\Bbbk^{r_1\times\cdots\times r_m} \to \Bbbk^{r'_1\times\cdots\times r'_n}$ also by $\sigma$. We say that $\sigma$ is a \textit{rank-$1$ preserver permutation} if it maps rank-$1$ tensors to rank-$1$ tensors, i.e., $\sigma(\Var{S}_\br) \subset \Var{S}_{\br'}$.
\end{definition}

Before we proceed, we need to recall some notation. 
For any partition $I_1 \sqcup \cdots \sqcup I_s = [m]$, the $(I_1,\dots,I_s)$\textit{-flattening} is the bijective linear map
\[
	{\rm flatt}_{(I_1,\dots,I_s)} : \Bbbk^{r_1\times\cdots\times r_m} \longrightarrow \bigotimes_{j \in [s]} \bigotimes_{k \in I_j} \Bbbk^{r_k}
\]
such that over rank-$1$ tensors $\tensor{A}=\ba_1\otimes\cdots\otimes\ba_m$ it is given by
\[
	\tensor{A}_{(I_1,\dots,I_s)} := {\rm flatt}_{(I_1,\dots,I_s)}(\tensor{A}) := \bigotimes_{j \in [s]} \bigotimes_{k \in I_j} \,\ba_k.
\]
The foregoing completely specifies their action because rank-$1$ tensors form a basis of the domain.
If $I = \{h\}$ and $J = [m] \smallsetminus \{h\}$, we call $\tensor{T}_{(I,J)}$ the \emph{$h$-flattening} of $\tensor{T}$ and we denote it by $\tensor{T}_{(h)}$. For any partition $I_1 \sqcup \cdots \sqcup I_s = [m]$, we also call the bijection $\varphi_{(I_1,\ldots,I_s)} : [\br] \rightarrow [\br']$ induced by the flattening ${\rm flatt}_{(I_1,\dots,I_s)}$ on the multi-indices a \emph{flattening}; that is,
\[
	[\tensor{T}_{(I_1,\dots,I_s)}]_{\varphi_{(I_1,\ldots,I_s)}(\bj)} = \tensor{T}_\bj, \quad \forall\bj \in [\br].
\]
Every flattening of a rank-$1$ tensor is a rank-$1$ tensor, by definition. It is well known that this is a complete characterization of rank-$1$ tensors; see \cref{lem:char_rank1_flattening}.

Then, from \cite{Westwick1967}, we obtain the following characterization of rank-$1$ preserver permutations that relates them to tensor products of permutations.

\begin{proposition}\label{prop_rk1presperm}
	Let $R = r_1\cdots r_m = r'_1\cdots r'_n$.
	Let $\rho : [R] \rightarrow [\br]$ and $\rho' : [R] \rightarrow [\br']$ be bijections and let $\sigma = \rho' \circ \rho^{-1} : [\br] \rightarrow [\br']$. Then, $\sigma$ is a rank-$1$ preserver permutation if and only if there exist permutations $\pi_k \in \mathfrak{S}_{r'_k}$, $k\in[n]$, and a flattening $\varphi : [\br] \rightarrow [\br']$ such that
	\(
 		\sigma = (\pi_1 \otimes\cdots\otimes \pi_n) \circ \varphi.
	\)
\end{proposition}
\begin{proof}
	Throughout the proof let $\Phi, \Sigma \in \RR^{R\times R}$, and $\Pi_k \in \RR^{r'_k\times r'_k}$ denote the permutation matrices corresponding to $\varphi$, $\sigma$, and the $\pi_k$'s, respectively. Note that we use the multi-indexing induced by the lexicographic ordering on $[\br]$ and $[\br']$ to label rows and columns of the matrix $\Phi$ and $\Sigma$.
	
	If $\sigma$ has the stated form, then by \cref{lem_tenprodperm}, the basic properties of tensor products of linear maps \cite{Greub1978}, and the fact that we can express $\ba'_1 \otimes \cdots \otimes \ba'_n = \Phi(\ba_1\otimes \cdots \otimes \ba_m)$ for any rank-$1$ tensor $\ba_1\otimes\cdots\otimes\ba_m$, it follows that
\[
 \Sigma(\ba_1\otimes\cdots\otimes\ba_m)
 = (\Pi_1 \ba'_1)\otimes\cdots\otimes (\Pi_n \ba'_n),
\]
which is again of rank equal to $1$. This proves the first direction.

Conversely, \cite[Theorem 3.4]{Westwick1967} entails that every invertible rank-$1$ preserver, hence including every $\Sigma$ representing $\sigma$, must be of the form
\[
	\Sigma = (F_1 \otimes \cdots \otimes F_n) \Phi, \quad\text{or, equivalently,}\quad
	\Sigma' := \Sigma \Phi^T = F_1 \otimes \cdots \otimes F_n,
\]
where each $F_k \in \Bbbk^{r_k' \times r_k'}$ is invertible and $\Phi$ represents a flattening from $\Bbbk^{r_1\times\cdots\times r_m} \to \Bbbk^{r'_1\times\cdots\times r'_n}$.

We will prove the following claim shortly.

\begin{claim}\label{claim:rank1_preservers_1}
	For every $k\in [n]$, the matrix $F_k$ is a scaled permutation matrix, i.e., there exists a permutation matrix $\Pi_k$ and a diagonal matrix $D_k$ so that $F_k = \Pi_kD_k$. 
\end{claim}

By \cref{claim:rank1_preservers_1}, we get
	\[
	 \Sigma' = (\Pi_1 D_1) \otimes\cdots\otimes (\Pi_n D_n) = (\Pi_1\otimes\cdots\otimes\Pi_n)(D_1\otimes\cdots\otimes D_n).
	\]
	The left hand side is a permutation matrix, while the right hand side is a permutation matrix whose columns are multiplied by nonzero scalars. If any of these scalars is different from $1$, then $\Sigma'$ cannot be a permutation matrix. Hence, $D_1\otimes\cdots\otimes D_{n} = I$. This concludes the proof.
\end{proof}

\begin{proof}[Proof of \cref{claim:rank1_preservers_1}]
	It is enough to show that
	\begin{enumerate}[nolistsep]
		\item $F_k$ has orthogonal rows, and
		\item every row of $F_k$ contains only one nonzero entry.
	\end{enumerate}
Since $\Sigma'$ is a permutation matrix, then it is orthogonal and it satisfies
\[
 	{\rm Id} = \Sigma' \Sigma'^{T} = (F_1 F_1^T)\otimes\cdots\otimes(F_{n} F_{n}^T).
\]
This implies that each $H^k := F_kF_k^T$ is a diagonal matrix. Indeed, for $\bi\in[\br']$, we have $1 = I_{\bi,\bi} = H^1_{i_1,i_1}\cdots H^{n}_{i_{n},i_{n}}$, so that the diagonal of each $H^k$ is nonzero. In turn, this implies, for arbitrary $k$ and $i_k \ne j_k$, that we have
\[
0 = {\rm Id}_{(1,\ldots,1,i_k,1,\ldots,1),(1,\ldots,1,j_k,1,\ldots,1)} = \underbrace{H^1_{1,1}\cdots H^{k-1}_{1,1}}_{\ne0} H^k_{i_k,j_k} \underbrace{H^{k+1}_{1,1}\cdots H^{{n}}_{1,1}}_{\ne0},
\]
so $H^k$ is nonzero only on the diagonal. Hence, each $F_k$ has orthogonal  rows.

	Next, consider an arbitrary $\bi,\bj \in [\br']$ such that
	\[
		\Sigma_{\bi,\bj}' = (F_1)_{i_1,j_1}\cdots(F_n)_{i_n,j_n} = 1.
	\]
	This implies that $(F_k)_{i_k,j_k} \neq 0$ for $k\in[n]$.
	Moreover, for every $\ell \ne j_k$, we have $(F_k)_{i_k,\ell} = 0$. Indeed, let $\bj_\ell=(j_1,\ldots,j_{k-1},\ell,j_{k+1},\ldots,j_n)$, and then we have
	\[
		0 = \Sigma_{\bi,\bj_\ell}' = \underbrace{(F_1)_{i_1,j_1}\cdots(F_{k-1})_{i_{k-1},j_{k-1}}}_{\ne 0} (F_k)_{i_k,\ell} \underbrace{(F_{k+1})_{i_{k+1},j_{k+1}} \cdots (F_{n})_{i_n,j_n}}_{\ne 0}.
	\]
	This shows that every row of $F_k$ contains only one nonzero element.
\end{proof}

We can now relate the various rank-$1$ permutations of an admissible vector.

\begin{proposition}\label{prop_allrk1perms}
Let $\rho : [R] \rightarrow [\br]$ be a bijective map.
Then, $\rho' : [R] \rightarrow [\br']$ is a rank-$1$ permutation for all $\ba \in \Var{A}^\rho$ if and only if there exists a rank-$1$ preserver permutation $\sigma : [\br] \to [\br']$ such that $\rho' = \sigma \circ \rho$.
\end{proposition}
\begin{proof}
	If there exists such a $\rho' = \sigma \circ \rho$, then, by \cref{def:rank1_preserver,def_admissible}, we get
	\[
		\Var{A}^\rho \subset \Var{A}^{\sigma \circ \rho} = \Var{A}^{\rho'}.
	\]
	Hence, $\rho'$ is a rank-$1$ permutation for all $\ba \in \Var{A}^\rho$.

	Conversely, consider $\sigma = \rho' \circ \rho^{-1}$. Since $\rho$ and $\rho'$ are rank-$1$ permutations by assumption, $\sigma$ is a rank-$1$ preserver permutation. This concludes the proof.
\end{proof}

	These results tell us that for a generic admissible vector $\ba \in \Var{A}^\rho$, there exists a unique \textit{maximal} rank-$1$ permutation. We say that $\rho$ is a maximal rank-$1$ permutation for $\ba \in \Var{A}^\rho$ if, for every $\rho' \neq \rho$ such that $\Var{A}^{\rho'} \subset \Var{A}^\rho$, then $\ba \not\in \Var{A}^{\rho'}$. By \cref{prop_rk1presperm,prop_allrk1perms}, this means that there is no rank-$1$ permutation $\rho'$ of $\ba$ such that $\rho$ is, up to rank-$1$ preservers, a flattening of $\rho'$. Note that, for every $\rho$, the set of $\ba$'s for which $\rho$ is a maximal rank-$1$ permutation is a Zariksi open subset of $\Var{A}^\rho$. Moreover, by \cref{prop_rk1presperm,prop_allrk1perms}, if, up to rank-$1$ preservers, $\rho$ and $\rho'$ are not related by a flattening, then $\Var{A}^\rho \cap \Var{A}^{\rho'}$ is a strict subvariety of $\Var{A}^{\rho}$. Hence, also the set of $\ba$'s for which $\rho$ is the \emph{unique} maximal rank-$1$ permutation is a Zariksi open subset of $\Var{A}^\rho$.

	In computing rank-$1$ permutations of a given vector $\ba \in \Bbbk^R$, we will always assume that $\ba$ is generic in the sense that it admits a unique maximal rank-$1$ permutation and our goal will be to retrieve this permutation.
	In this way, we can always assume that the size of the maximal rank-$1$ permutation is fixed and we will employ \cref{prop_allrk1perms} under the assumption that $[\br] = [\br']$, namely, we will use the fact that $\Var{A}^\rho = \Var{A}^{\rho'}$ if and only if $\rho' = \sigma \circ \rho$ for some $\sigma \in \mathfrak{S}([\br])$.

	When we say that $\ba \in \Bbbk^R$ is a \textit{generic admissible vector}, then we specifically mean {that it has} a maximal rank-$1$ permutation that is unique up to rank-$1$ preservers.

\subsection{Computing a rank-1 permutation}\label{sec_rk1perm}

Conceptually, computing a rank-$1$ permutation of a given $\ba$, or certifying that one does not exist, is trivial: try out each of the $R!$ permutations and check whether $\pi(\ba)$ is a rank-$1$ tensor.
This is feasible when $R \le 10$, but due to the superexponential growth quickly becomes infeasible for the next few values. In the context of \ref{eqn_hhd}s, this means that the six rank parameters $\br=(2,2),(2,3),(2,4),(2,5),(3,3),(2,2,2)$ are the only feasible ones with this approach. By exploiting the symmetries from \cref{prop_allrk1perms}, a clever enumeration approach might reduce the worst-case complexity to about $\frac{R!}{r_1! \cdots r_m!}$. Unfortunately, this still grows superexponentially with $R=r_1\cdots r_m$. A different approach is required.

As we are unaware of other algorithms for computing rank-$1$ permutations, we present a new polynomial-time algorithm in the next subsection. We have to make one concession, however: the proposed algorithm applies only to \textit{generic} inputs. That is, it will output a correct solution on a Zariski open dense subset of inputs, but can fail on the other inputs.

Before proceeding, note that we can make some simplifying assumptions. If $\rho$ is a rank-$1$ permutation of an admissible $\ba \in \Bbbk^R$, then $\sigma \circ \rho$ is still a rank-$1$ permutation for every rank-$1$ preserver $\sigma \in \mathfrak{S}([\br])$, by \cref{prop_allrk1perms}. Moreover, if $\sigma \in \mathfrak{S}_{R}$ is any permutation, then $\rho \circ \sigma$ is a rank-$1$ permutation of $\sigma^{-1}(\ba)$.
Hence, we will assume without loss of generality in the remainder of this section that $\ba \in \Bbbk^R$ is sorted by decreasing magnitude, i.e., $|a_1|\geq |a_2| \geq \dots \geq |a_R|$, and that the element $a_1$ is mapped by the rank-$1$ permutation $\rho$ to position ${\bf1} := (1,\ldots,1)$, i.e., $\rho(1) = \bf 1$. Since the case $\ba = 0$ is trivial, we may assume that at least one entry of $\ba$ is nonzero and we can assume that such nonzero entry appears in $\rho(\ba)$ in position ${\bf1}$.

We identify the maximal rank-$1$ permutation of a generic admissible vector by exploiting the following two characterizations of rank-$1$ tensors. 

\begin{lemma}[{\cite[Proposition 8.3.7]{BC2019}}]\label{lem:char_rank1_flattening}
	Let $\tensor{A} \neq 0 \in \Bbbk^{r_1\times\cdots\times r_m}$. The following are equivalent:
	\begin{itemize}[nolistsep]
		\item $\tensor{A}$ has rank $1$;
		\item $\tensor{A}_{(I_1,\dots,I_s)}$ has rank $1$ for every partition $I_1 \sqcup \cdots \sqcup I_s = [m]$;
		\item $\tensor{A}_{(h)}$ has rank $1$ for every $h \in [m]$.
	\end{itemize}
\end{lemma}

\begin{lemma}[\cite{EGH2009}]\label{rmk:char_rank1}
Let $\tensor{A} \neq 0 \in \Bbbk^{r_1\times\cdots\times r_m}$. Then,
$\tensor{A}$ is of rank $1$ if and only if it has an \textit{interpolatory decomposition}, i.e., if $a := \tensor{A}_{i_1,\ldots,i_m} \neq 0$, then \[
 \tensor{A} = a \cdot
 \begin{bmatrix}
  \tensor{A}_{(1,i_2,\ldots,i_m)} a^{-1} \\
  \vdots\\
  \tensor{A}_{(r_1,i_2,\ldots,i_m)} a^{-1}
 \end{bmatrix} \otimes \cdots \otimes
 \begin{bmatrix}
  \tensor{A}_{(i_1,\ldots,i_{m-1},1)} a^{-1}\\
  \vdots \\
  \tensor{A}_{(i_1,\ldots,i_{m-1},r_m)} a^{-1}
 \end{bmatrix}.
\]
\end{lemma}

For the following definition, recall that the \textit{Hamming distance} of two ordered vectors of the same length is the number of entries that are different.

\begin{definition}\label{def_cross}
	Let $\tensor{T} \in \Bbbk^{r_1\times\cdots\times r_m}$. For any $\mathbf{i} = (i_1,\ldots,i_m) \in [\bn]$, the \emph{$\mathbf{i}$-cross} of $\tensor{T}$ is union of all entries of $\tensor{T}$ whose index is at Hamming distance at most $1$ from ${(i_1,\ldots,i_m)}$, i.e.,
	\[
		\Var{C}_{\mathbf{i}}(\tensor{T})
		= \underbrace{\{ \tensor{T}_{j,i_2,i_3,\ldots,i_m} : j \in [r_1] \}}_{\Var{C}^1_{\mathbf{i}}(\tensor{T})}
		\cup \dots \cup \underbrace{\{ \tensor{T}_{i_1,\ldots,i_{m-1},j} : j \in [r_m] \}}_{\Var{C}^m_{\mathbf{i}}(\tensor{T})}.
	\]
	We call the subset $\Var{C}^h_{\mathbf{i}}(\tensor{T})$ the \emph{$h$th direction} of the $\mathbf{i}$-cross. Note that the $h$th direction of the ${\bf1}$-cross comprises the elements of the first column of $\tensor{A}_{(h)}$.
\end{definition}

	To determine the rank-1 permutations of an admissible vector $\ba \in \Bbbk^R$ such that $\rho(\ba)_{1,\ldots,1} \neq 0$, \cref{rmk:char_rank1} states that it suffices to determine the elements of $\ba$ that are mapped into the ${\bf1}$-cross.
	Indeed, the entire rank-$1$ tensor is a function of these elements.
	The next result shows how \cref{lem:char_rank1_flattening} offers a tool to identify which elements of $\ba$ are mapped to the $\bf 1$-cross by looking at the set $\Var{M}(\ba)$ of all vanishing $2 \times 2$ \textit{top-left corner minors}, i.e., the minors
	\begin{equation}\label{eq:vanishing_minors}
	\Var{M}(\ba) := \left\{m_{i,j,k} :=
	\begin{bmatrix}
  	a_1 & a_j \\
  	a_i & a_k
	\end{bmatrix} : \mathrm{det}(m_{i,j,k}) = 0 \text{ and } \sharp\{1,i,j,k\}=4 \right\}.
	\end{equation}
	Note that all the $2 \times 2$ top-left corner minors of the $h$-flattenings of $\tensor{A}$ are like this.

\begin{lemma}\label{lemma:rank1_a}
	Let $\ba \in \Bbbk^R$ be a generic admissible vector and let $\rho$ be any of its maximal rank-$1$ permutations with $\rho(1) = {\bf 1} =(1,\ldots,1)$. Let $\tensor{A} = \rho(\ba) \in \Bbbk^{r_1\times\cdots\times r_m}$. Then, the complement of the $\bf 1$-cross of $\tensor{A}$ is given by the set of elements appearing in the $(2,2)$-entry of a $2 \times 2$ top-left corner vanishing minor of $\ba$ defined in \cref{eq:vanishing_minors}, i.e.,
	\[
		\Var{C}_{\bf1}(\tensor{A}) = \{a_i : i \in [R]\} \setminus\{a_k : \exists i,j \in [R], ~m_{i,j,k}\in\Var{M}(\ba) \}.
	\]
\end{lemma}
\begin{proof}
	By \cite[Section 6.4]{BC2019}, $2 \times 2$ matrices of the form
	\[
 	\begin{bmatrix}
  		\tensor{A}_{\mathbf{l}} & \tensor{A}_{\mathbf{j}} \\
  		\tensor{A}_\mathbf{i} & \tensor{A}_\mathbf{k}
 	\end{bmatrix},
	\]
with distinct $\mathbf{i},\mathbf{j},\mathbf{k},\mathbf{l} \in [\br]$, whose determinant vanishes on a generic $\tensor{A}\in\Var{S}$ must originate as a submatrix of some flattening of $\tensor{A}$. 
Therefore, by maximality of $\rho$, those in $\Var{M}(\ba)$ can only originate as minors of flattenings of $\rho(\ba) = \tensor{A}$. Now, note that the Hamming distance of the index of any element of $\tensor{A}$ to ${\bf 1}$ can only decrease by flattenings.
In particular, if an element $a_k$ belongs to a top-left corner vanishing minor as in \cref{eq:vanishing_minors}, then it means that, in $\rho(\ba)$, it has Hamming distance at least two from $\tensor{A}_{\bf 1}$. Hence, it belongs to the complement of the $\bf 1$-cross of $\tensor{A}$.

Conversely, if an entry $a_k$ is not in the $\bf 1$-cross of $\tensor{A}$, then there exists a flattening in which $a_k$ is not mapped to the $(1,1)$-cross of the flattening. Hence, it appears in the $(2,2)$-entry of a top-left corner vanishing minor of $\ba$ defined in \cref{eq:vanishing_minors}.
\end{proof}

Observe that \cref{lemma:rank1_a} states that for \emph{every} maximal rank-$1$ permutation $\rho$ of a generic admissible vector $\ba\in\Bbbk^R$ with $\rho({\bf1})=1$, the ${\bf1}$-cross of $\rho(\ba)$ is the same set of numbers.

Now that we {characterized} the entries of a generic admissible $\ba \in \Bbbk^R$ that are mapped into the $\bf 1$-cross by all of its maximal rank-$1$ permutations, we need to understand the partition of such elements into the directions of the cross to reconstruct the whole rank-$1$ tensor {via} \cref{rmk:char_rank1}. The next result provides a criterion to determine which elements are mapped to the same direction in the cross.

\begin{lemma}\label{lemma:rank1_b}
	Let $\ba \in \Bbbk^R$ be a generic admissible vector and let $\rho$ be any of its maximal rank-$1$ permutations with $\rho(1)=(1,\ldots,1)$. Let $\tensor{A} = \rho(\ba) \in \Bbbk^{r_1\times\cdots\times r_m}$. Then, for any $2 \le i \neq j \le R$, the elements $a_i$ and $a_j$ belong to the same direction $h$ of the $\mathbf{1}$-cross of $\tensor{A}$ if and only if the minors $m_{i,j,k} \in \Var{M}(\ba)$ do not vanish for any $k \neq i,j$. Equivalently, for $a_i,a_j\in\Var{C}_{\bf1}(\tensor{A})\setminus\{a_1\}$ we have
	\[
		\{a_i,a_j\} \subset \Var{C}^h_{\mathbf{1}}(\tensor{A}) \setminus\{a_1\} \quad \text{ if and only if } \quad \forall k\neq i,j:  m_{i,j,k} \not\in \Var{M}(\ba).
	\]
\end{lemma}
\begin{proof}
	We prove the negation of the statement. 

	Let $a_i, a_j \in \Var{C}_{\bf1}(\tensor{A})\setminus\{a_1\}$ and assume $a_i$ is mapped to the $h$th direction and $a_j$ to the $h'$th direction. If $h \ne h'$, then by taking the $h$-flattening, $a_i$ is mapped to the first column of $\tensor{A}_{(h)}$, while $a_j$ is mapped to the first row.
	Since $a_i, a_j \ne a_1$, neither of them is mapped to position $(1,1)$ in $\tensor{A}_{(h)}$. Consequently, since $\tensor{A}_{(h)}$ is a rank-$1$ matrix by \cref{lem:char_rank1_flattening} there is a vanishing minor $m_{i,j,k} \in \Var{M}(\ba)$ for some $k$.

	Conversely, let $k\ne i,j$ be such that $m_{i,j,k} \in \Var{M}(\ba)$ vanishes. Assume there exists a direction $h$ such that both $a_i, a_j \in \Var{C}_{\bf1}^h(\tensor{A})\setminus\{a_1\}$. Then, in the $(I,J)$-flattening of $\rho(\ba)$, for arbitrary $I\sqcup J=[m]$, they are either both mapped to the first column if $h \in I$ or both to the first row if $h\in J$. However, by genericity of $\ba$ in $\Var{A}^\rho$, the only vanishing minors in $\Var{M}(\ba)$ are the ones originating from the flattenings by \cref{lem:char_rank1_flattening}. This means there is no vanishing minor $m_{i,j,k} \in \Var{M}(\ba)$, concluding the proof.
\end{proof}

By combining \cref{lemma:rank1_a,rmk:char_rank1,lemma:rank1_b}, we arrive at \cref{alg_basic_minors,alg_rank1_permutation} for computing a rank-$1$ permutation of a generic admissible vector $\ba\in\Bbbk^R$. 

\begin{algorithm}[tb]\small
\caption{{Compute the vanishing top-left corner minors of a generic admissible vector.}}
\begin{algorithmic}[1]
\Require $\ba$ is a generic admissible vector in $\Bbbk^R$.
\vspace{4pt}
\State $\Var{M} \leftarrow []$; \label{line_step1start}
\For{$2 \le i < j \le R$}
	\State Let $k$ be the first index such that $a_k = a_i a_1^{-1} a_j \in \ba$;
	\If{$k \not\in \{\infty, i, j\}$}\label{line_binsearch}
 		\State $\Var{M} \leftarrow \Var{M} \cup \{ (i,j,k), (j,i,k) \}$;
	\EndIf
\EndFor
\State \Return $\Var{M}$; 
\end{algorithmic}
\label{alg_basic_minors}
\end{algorithm}

{\Cref{alg_basic_minors} is a basic algorithm for identifying the set of vanishing minors $\Var{M}(\ba)$ in \cref{eq:vanishing_minors} from $\ba$. We do not actually store the $2\times2$ matrices $m_{i,j,k}$ themselves, but rather only keep track of the three index positions $i,j,k \in [R]$.}
{If binary search is exploited to look up values in $\ba$, after sorting it, then its time complexity is only $\mathcal{O}(R^2 \log_2 R)$ operations.}
\Cref{alg_basic_minors} is presented as a symbolical algorithm, but it also can be executed in floating-point arithmetic. {The discussion of how to compute a rank-$1$ permutation robustly in the presence of round-off errors is delayed until \cref{sec_sub_advanced_version}.}

\begin{algorithm}[tb]\small
\caption{Compute a maximal rank-$1$ permutation $\rho : [R]\to[\br]$ of a generic admissible vector.}
\begin{algorithmic}[1]
\Require $\ba$ is a generic admissible vector in $\Bbbk^R$.
{\Require $\Var{M}$ contains the vanishing minors of $\ba$.}
\vspace{4pt}
\Statex\Comment{Extract the ${\bf1}$-cross $\Var{C}_{\bf1}$.}
\State $\Var{C}_{\bf1} \leftarrow \{2,3,\ldots,R\}\setminus\{ k : (i,j,k) \in \Var{M} \}$; \label{line_step2}
\vspace{4pt}
\Statex \Comment{Partition the ${\bf1}$-cross $\Var{C}_{\bf1} = \Var{C}^1_{\bf1} \sqcup \Var{C}^2_{\bf1} \sqcup \cdots \sqcup \Var{C}^m_{\bf1}$.}
\State $\Var{C}^\ell_{\bf1} \leftarrow \emptyset$ for $\ell=1,\ldots,m$; \label{line_step3start}
\For{$s = 1, 2, \ldots, m$}
	\State Let $i$ be the first element from $\Var{C}_{\bf1}$;
	\State $H \leftarrow \Var{C}_{\bf1} \setminus \{ j : (i,j,k) \in \Var{C}_{\bf1} \}$;\label{line_sort1}
	\State Let $1 \le \ell \le m$ be the first number such that $\Var{C}^\ell_{\bf1} = \emptyset$ and $r_\ell = \sharp H + 1$;
	\State $\Var{C}_{\bf1}^\ell = (\Var{C}_{\bf1}^\ell(2),\ldots,\Var{C}_{\bf1}^\ell(r_\ell)) \leftarrow H$;
	\State $\Var{C}_{\bf1} \leftarrow \Var{C}_{\bf1} \setminus \Var{C}^\ell_{\bf1}$;\label{line_sort2}
\EndFor \label{line_step3stop}
\vspace{4pt}
\Statex\Comment{Compute the tensor $\tensor{A}$ from its interpolatory decomposition.}
\State $\mathbf{x}_\ell \leftarrow [1, a_1^{-1} a_{\Var{C}^\ell_{\bf1}(2)}, \ldots, a_1^{-1} a_{\Var{C}^\ell_{\bf1}(r_\ell)}]$ for $\ell=1,\ldots,m$; \label{line_step4start}
\State $\tensor{A} \leftarrow a_1 \cdot \mathbf{x}_1 \otimes\cdots\otimes \mathbf{x}_m$; \label{line_step4stop}
\vspace{4pt}
\Statex\Comment{Determine the rank-$1$ permutation.}
\State Sort $\mathrm{vec}(\tensor{A})$ and let $\pi_1$ denote the permutation that accomplishes this; \label{line_step5start}
\State {Sort $\ba$ and let $\pi_2$ denote the permutation that accomplishes this;} \label{line_alsosort}
\State $\rho \leftarrow \mathrm{vec}^{-1} \circ {(\pi_1^{-1} \circ \pi_2)}$; \label{line_step5stop}
\State \Return $\rho$;
\end{algorithmic}
 \label{alg_rank1_permutation}
\end{algorithm}

{\Cref{alg_rank1_permutation} takes as input a generic admissible vector $\ba$ and its vanishing minors $\Var{M}$, computed for example by \cref{alg_basic_minors}, and computes the rank-$1$ permutation from them. Its} key steps are as follows:
\begin{itemize}[nolistsep]
	\item line \ref{line_step2}: we extract the ${\bf1}$-cross as suggested by \cref{lemma:rank1_a};
	\item lines \ref{line_step3start}--\ref{line_step3stop}: we partition the ${\bf1}$-cross as explained in \cref{lemma:rank1_b};
	\item lines \ref{line_step4start}--\ref{line_step4stop}: we exploit the interpolatory decomposition from \cref{rmk:char_rank1} to compute the values of $\tensor{A}\in\Bbbk^{r_1\times\cdots\times r_m}$;
	\item lines \ref{line_step5start}--\ref{line_step5stop}: we finally efficiently determine the permutation $\rho$ from $\ba$ to $\tensor{A} \in \Bbbk^R$ by sorting the values of $\ba$ and the vectorization $\mathrm{vec}(\tensor{A})$ of $\tensor{A}$.
\end{itemize}

The time complexity of \cref{alg_rank1_permutation} is straightforward to determine. We assume the following implementation choices are made:
\begin{itemize}[nolistsep]
 \item In lines \ref{line_step2}, \ref{line_sort1}, and \ref{line_sort2} a sorting algorithm with time complexity $\mathcal{O}(n \log_2 n)$ for vectors of length $n$ is used to compute the set differences.
 \item Lines \ref{line_step5start} and \ref{line_alsosort} use an $\mathcal{O}(n \log_2 n)$ sorting algorithm.
\end{itemize}
Then, we have
\begin{equation*}
 \underbrace{\mathcal{O}(R^2 \log_2 R)}_\text{line \ref{line_step2}} +
 \underbrace{\mathcal{O}(m S \log_2 S)}_\text{lines \ref{line_step3start}--\ref{line_step3stop}} +
 \underbrace{\mathcal{O}(m R)}_\text{lines \ref{line_step4start}--\ref{line_step4stop}} +
 \underbrace{\mathcal{O}(R \log_2 R)}_\text{lines \ref{line_step5start}--\ref{line_step5stop}} = \mathcal{O}( R^2 \log_2 R ),
\end{equation*}
where $S = r_1 + \cdots + r_m$ and $R = r_1 \cdots r_m$. Since $m \le S \le R$, \cref{alg_rank1_permutation} has a quasi-quadratic time complexity in $R$, the product of the HHD ranks.

\section{Identifiability of low-rank HHDs}\label{sec:identifiability}

The naive idea underlying our approach is the observation in \cref{sec:preliminaries_hhd} that by distributivity of the Hadamard multiplication with the sum, every \ref{eqn_hhd} of a tensor $\tensor{T} \in \Bbbk^{n_1\times\cdots\times n_d}$ induces a special Hitchcock decomposition of $\tensor{T}$:
\begin{equation*} 
	\tensor{T} 
	= \bhadamard{k=1}{m} \sum_{i_k=1}^{r_k} \tensor{A}_{k i_k} 
	= \sum_{\bi\in[\br]} \bhadamard{k=1}{m} \tensor{A}_{k i_k}.
\end{equation*}
Here, we show that if the induced Hitchcock decomposition is $R$-identifiable, then the \ref{eqn_hhd} is also identifiable, up to the symmetries that are described next.

It follows immediately from the definition of an \ref{eqn_hhd} as a Hadamard product of Hitchcock decompositions that an \ref{eqn_hhd} is invariant under the action of $\mathfrak{S}_{r_k}$ on the summands of the $k$th Hadamard factor for every $k \in [m]$. {Moreover,} Hadamard factors whose Hitchcock decompositions have the same length can be permuted without consequence. This is encoded by the action of the subgroup of $\mathfrak{S}_m$ that preserves the multi-index $\br \in \NN^m$:
\begin{align}\label{eqn_prespermgroup}
	\mathfrak{P}(\br) = \{ \sigma \in \mathfrak{S}_m ~:~ \sigma(\br)=\br \}.
\end{align}
That is, we naturally regard an \ref{eqn_hhd} as an element of
\begin{equation}\label{eq:upto_permutation}
	\Var{S}^{(\br)} := (\Var{S}^{(r_1)}\times\cdots\times\Var{S}^{(r_m)}) / \mathfrak{P}({\br}).
\end{equation}
Similarly as we did for CPDs, we consider the \textit{set of HHDs} of a given tensor $\tensor{T} \in \Bbbk^{n_1\times\cdots\times n_d}$. If we fix the rank parameters $\br$, then the set of \ref{eqn_hhd}s of $\tensor{T}$ is
\[	
	\Var{H}_\br(\tensor{T}) = \Bigl\{ [(\{ \tensor{A}_{k1},\ldots,\tensor{A}_{k r_k}\})_{k=1}^{m}] \in \Var{S}^{(\br)} : \tensor{T} = \bhadamard{k=1}{m} \sum_{i_k=1}^{r_k} \tensor{A}_{k i_k}\Bigr\}.
\]
{Herein, the square brackets denote the equivalence class under the quotient by $\mathfrak{P}(\br)$ in \eqref{eq:upto_permutation}. In other words, if $(\tensor{T}_1,\dots,\tensor{T}_m)\in\Var{S}^{(r_1)}\times\dots\times\Var{S}^{(r_m)}$, then $[(\tensor{T}_1,\dots,\tensor{T}_m)] \in \Var{S}^{(\br)}$ is the equivalence class of all tuples that can be obtained from $(\tensor{T}_1,\dots,\tensor{T}_m)$ by permuting all $\tensor{T}_i$ and $\tensor{T}_j$ for which $r_i=r_j$.}

Contrary to the CPD case, there is a nondiscrete action preserving the set $\Var{H}_\br(\tensor{T})$. Consider the group $(\Var{S}^{*},\hadamard)$ of rank-$1$ tensors with nonzero entries in $\Bbbk^{n_1\times\cdots\times n_d}$ equipped with Hadamard multiplication. It is a group with unit $\mathds{1}={\bf1}\otimes\cdots\otimes{\bf1}\in\Var{S}^*$ because the Hadamard inverse of a rank-$1$ tensor with all entries nonzero is still of rank $1$ (see \cref{sec:preliminaries_hhd}), i.e., $\tensor{A} \in \Var{S}^{*}$ if and only if $\tensor{A}^{\ominus}\in \Var{S}^{*}$, and the Hadamard product of rank-$1$ tensors with nonzero entries is again a rank-$1$ tensor, i.e., if $\tensor{A}, \tensor{B} \in \Var{S}^*$ then $\tensor{A}\hadamard\tensor{B} \in \Var{S}^*$. Now, we observe that if $H := [(\{ \tensor{A}_{k1},\ldots,\tensor{A}_{k r_k} \})_{k=1}^m] \in \Var{H}_\br(\tensor{T})$, then we also have the following set of nontrivial alternative \ref{eqn_hhd}s of $\tensor{T}$:
\begin{equation}\label{eq:equivalent_HHDs}
[H]_{\Var{S}^*} :=
\Bigl\{ [(\{ \tensor{D}_k \hadamard \tensor{A}_{k1},\ldots,\tensor{D}_k \hadamard\tensor{A}_{k r_k} \})_{k=1}^m] \in \Var{H}_\br(\tensor{T}) : \hspace{-5pt}\begin{array}{c}\tensor{D}_1 \hadamard\cdots\hadamard \tensor{D}_m = \mathds{1} \\ \text{ with } \tensor{D}_k\in\Var{S}^*, ~\forall k \in [m]\end{array}\Bigr\}.
\end{equation}
In conclusion, we can identify the constituent rank-$1$ tensors in the Hadamard factors of an \ref{eqn_hhd} up to the actions of
\begin{enumerate}[nolistsep]
 \item[(i)] $\mathfrak{S}_{r_k}$, which reorders the rank-$1$ tensors in the $k$th Hadamard factor,
 \item[(ii)] $\mathfrak{P}(\br)$, which reorders the Hadamard factors of equal length, and
 \item[(iii)] $\Var{S}^*$, which Hadamard multiplies the Hadamard factors by rank-$1$ tensors that factor $\mathds{1}$.
\end{enumerate}
{The actions (i) and (ii) were taken in consideration in the definition of $\Var{H}_\br(\tensor{T})$. In order to mod out also action (iii), we let} 
\[
\Var{H}^{\Var{S}^{*}}_\br(\tensor{T}) 
:= \Var{H}_\br(\tensor{T}) / \Var{S}^* 
= \{ [H]_{\Var{S}^{*}} : H \in \Var{H}_\br(\tensor{T}) \}
\]
denote the set of \textit{essential} \ref{eqn_hhd}s of $\tensor{T}$. 

An \ref{eqn_hhd} is called $\br$-identifiable if there is but one essential decompostion.

\begin{definition}
	The tensor $\tensor{T}$ is \textit{$\br$-identifiable} if $\Var{H}^{\Var{S}^{*}}_\br(\tensor{T})$ is a singleton.
\end{definition}

\begin{remark}
It follows from this definition that an \ref{eqn_hhd} is not $\br$-identifiable if there exists a factor $k\in[m]$ such that the $k$th Hadamard factor is not an $r_k$-identifiable CPD. Another way $\br$-identifiability may fail is if $\tensor{T}$ factors as
\[
 \tensor{T}= \tensor{T}_1 \hadamard\cdots\hadamard \tensor{T}_m = \tensor{T}_1'\hadamard\cdots\hadamard\tensor{T}_m'
\]
in such a way that for every $\psi \in \mathfrak{P}(\br)$ there exists an $i \in [m]$ such that $\tensor{T}_i^\ominus \hadamard \tensor{T}'_{\psi(i)}$ is not of rank-$1$. 
\end{remark}

The key to distinguishing the essential \ref{eqn_hhd}s is an auxiliary tensor with a special \textit{double rank-$1$ structure}. For brevity, let $\Var{S}^{\times\br}=\Var{S}^{\times r_1}\times\cdots\times\Var{S}^{\times r_m}$ denote the space of \emph{ordered \ref{eqn_hhd}s} in the next definition and proposition.

\begin{definition}\label{def_stacked}
Let $H= ( \tensor{A}_{k1}, \ldots, \tensor{A}_{k i_k} )_{k=1}^m \in \Var{S}^{\times \br}$ be an ordered \ref{eqn_hhd}.
The \emph{Hadamard--Hitchcock tensor} (HHT) $\tensor{H}=h(H)$ is the $(m+d)$-tensor of size $r_1\times\cdots\times r_m \times n_1\times\cdots\times n_d$ that satisfies
\begin{enumerate}[nolistsep]
 \item $\tensor{H}_{\bi,\bj} = (\tensor{A}_{1 i_1})_\bj \cdots (\tensor{A}_{m i_m})_\bj$, for all $\bi\in[\br]$ and $\bj\in[\bn]$;
 \item $\tensor{H}_{\bi,:} = \tensor{A}_{1 i_1} \hadamard\cdots\hadamard \tensor{A}_{m i_m} \in \Var{S}_{\bn}$, for all $\bi\in[\br]$; and 
 \item $\tensor{H}_{:,\bj} = A^1_\bj \otimes\cdots\otimes A^m_\bj \in \Var{S}_{\br}$, for all $\bj \in [\bn]$,
\end{enumerate}
where $A^k_\bj = \bigl( (\tensor{A}_{k1})_\bj, \ldots, (\tensor{A}_{k r_k})_\bj \bigr) \in \Bbbk^{r_k}$ for all $k\in[m]$ and the subscript colon denotes the set of all other indices.
\end{definition}

Note that the three items are all equivalent; {using elementary computations one} verifies that $(1) \Leftrightarrow (2)$ and $(1) \Leftrightarrow (3)$. Consider the matrix $\overline{\tensor{H}} \in \Bbbk^{R \times N}$ obtained as the $(I,J)$-flattening of the tensor $\tensor{H}$ for $I=(1,\ldots,m)$ and $J=(m+1,\ldots,m+d)$, which we call the \emph{natural flattening} of the HHT. Then, item $(2)$ of \Cref{def_stacked} essentially states that all the rank-$1$ tensors $\tensor{A}_{1 i_1}\hadamard\cdots\hadamard \tensor{A}_{m i_m} \in \Bbbk^{n_1\times\cdots\times n_d}\simeq\Bbbk^N$ are placed as rows of $\overline{\tensor{H}} \in \Bbbk^{R \times N}$ in a suitable order so that the columns of this matrix can be interpreted (by identifying $\Bbbk^R \simeq \Bbbk^{r_1\times\cdots\times r_m}$) as rank-$1$ tensors as well.
The HHT $\tensor{H}=h(H)$ has the nice property that it is invariant to the action of $\Var{S}^*$ on generic $H$.

With a little abuse of notation, we denote by $[H]_{\Var{S}^*}$ a set of equivalent \ref{eqn_hhd}s as in \Cref{eq:equivalent_HHDs}, but considering them as ordered \ref{eqn_hhd}s so that we can define the HHT for any $H' \in [H]_{\Var{S}^*}$.

\begin{proposition}\label{prop_stacked}
Let $H=(\tensor{A}_{k1},\ldots,\tensor{A}_{k r_k})_{k=1}^m \in \Var{S}^{\times\br}$ be an ordered \ref{eqn_hhd}. 
Let $\tensor{H}=h(H)$ be the HHT.
If there exists an index $\bj\in[\br]$ such that $\tensor{A}_{1 j_1} \hadamard\cdots\hadamard \tensor{A}_{m j_m} \in \Var{S}^*$, then $h^{-1}(\tensor{H}) = [H]_{\Var{S}^*}$.
\end{proposition}
\begin{proof}
	The inclusion $[H]_{\Var{S}^*} \subset h^{-1}(\tensor{H})$ follows from \cref{def_stacked}$(2)$ and the associativity of the Hadamard product: for any $H' \in [H]_{\Var{S}^*}$ and $\bi\in[\br]$, we have
\[
 h(H')_{\bi,:}
 = (\tensor{D}_1 \hadamard \tensor{A}_{1 i_1})\hadamard\cdots\hadamard(\tensor{D}_m \hadamard \tensor{A}_{m i_m}) 
 = (\tensor{D}_1 \hadamard\cdots\hadamard\tensor{D}_m) \hadamard \tensor{H}_{\bi,:}
 = \tensor{H}_{\bi,:},
\]
This inclusion holds even without the assumption on the existence of the special $\bj$ in the claim.

Conversely, let $H' = (\tensor{B}_{k1},\ldots,\tensor{B}_{k r_k})_{k=1}^m \in h^{-1}(\tensor{H})$ be an arbitrary element in the fiber. We show that $H' \in [H]_{\Var{S}^*}$.
Without loss of generality, up to reordering the rank-$1$ summands of each Hadamard factor by $\mathfrak{S}_{r_1}\times\cdots\times\mathfrak{S}_{r_m}$, we can assume $\bj={\bf1}$, simplifying the notation.

The assumption implies that $\tensor{A}_{k1} \in \Var{S}^*$, i.e., it has nonzero elements, for every $k\in[m]$.
Since $h(H') = \tensor{H}$, we have from \cref{def_stacked}(2) that
\[
 \tensor{H}_{{\bf1},:} =  \tensor{A}_{1 1}\hadamard\cdots\hadamard\tensor{A}_{m 1} = \tensor{B}_{1 1}\hadamard\cdots\hadamard \tensor{B}_{m 1}.
\]
Since $\tensor{H}_{{\bf1},:} \in \Var{S}^*$, we have $\tensor{B}_{k1}\in\Var{S}^*$ for every $k\in[m]$ as well.
Then, for all $i_k \in [r_k]$, we see that
\begin{align*}
\tensor{H}_{(1,\ldots,1,i_k,1,\ldots,1),:} 
 &= \tensor{A}_{1,1}\hadamard\cdots\hadamard\tensor{A}_{k-1,1} \hadamard\tensor{A}_{k,i_k}\hadamard \tensor{A}_{k+1,1}  \hadamard\cdots\hadamard\tensor{A}_{m 1} \\
 &= \tensor{B}_{1,1}\hadamard\cdots\hadamard \tensor{B}_{k-1,1} \hadamard\tensor{B}_{k,i_k}\hadamard \tensor{B}_{k+1,1}   \hadamard\cdots\hadamard\tensor{B}_{m 1}.
\end{align*}
Dividing this elementwise by $\tensor{H}_{{\bf1},:}$ results in nontrivial relations for $i_k\ge2$, namely
\[
\tensor{H}_{(1,\ldots,1,i_k,1,\ldots,1),:} \hadamard \tensor{H}_{{\bf1},:}^{\ominus}
= \tensor{A}_{k i_k} \hadamard \tensor{A}_{k1}^{\ominus} 
= \tensor{B}_{k i_k} \hadamard \tensor{B}_{k1}^{\ominus}.
\]
This implies that for all $k$ and all $i_k$ we have 
\[
 \tensor{B}_{k i_k} = \underbrace{\tensor{B}_{k 1} \hadamard\tensor{A}_{k 1}^{\ominus}}_{\tensor{D}_k \in \Var{S}^*} \hadamard \tensor{A}_{k i_k}.
\]
Since we also have 
\begin{align*}
 \tensor{D}_1 \hadamard\cdots\hadamard\tensor{D}_m 
 &= (\tensor{B}_{11}\hadamard\tensor{A}_{11}^{\ominus}) \hadamard\cdots\hadamard(\tensor{B}_{m1}\hadamard\tensor{A}_{m1}^{\ominus})\\
 &= (\tensor{A}_{11}\hadamard\cdots\hadamard\tensor{A}_{m1})^{\ominus} \hadamard(\tensor{B}_{11}\hadamard\cdots\hadamard\tensor{B}_{m1})
 = \tensor{H}_{{\bf1},:}^{\ominus}\hadamard \tensor{H}_{{\bf1},:} = \mathds{1},
\end{align*}
it follows that $H' \in [H]_{\Var{S}^*}$ as required. This concludes the proof.
\end{proof}

We can now prove the first main result about the identifiability of \ref{eqn_hhd}s.
	
\begin{theorem}\label{thm:main_identifiabilityHHD}
Let $\br\in\NN^m$ and $R=r_1\cdots r_m$. Let $\tensor{T}\in\sigma_\br(\Var{S}) \subset\Bbbk^{n_1\times\cdots\times n_d}$ be generic. If $\tensor{T}$ is $R$-identifiable, then $\tensor{T}$ is $\br$-identifiable.
\end{theorem}
\begin{proof}
	Assume that we have two \ref{eqn_hhd}s of the tensor 
	\[
		\tensor{T} = \bhadamard{k=1}{m} \sum_{i_k=1}^{r_k} \tensor{A}_{k i_k} = \bhadamard{k=1}{m} \sum_{i_k=1}^{r_k} \tensor{A}_{k i_k}',
	\]
	which induce two Hitchcock decompositions as in \cref{ssec_hhd_cpd}, namely
	\begin{equation*}
		\tensor{T} = \sum_{\bj \in [\br]} \tensor{A}_{\bj} =  \sum_{\bj \in [\br]}\tensor{A}_{\bj}',
	\end{equation*}
	where $\tensor{A}_\bj = \tensor{A}_{1j_1}\hadamard\cdots\hadamard\tensor{A}_{mj_m}$, for all $\bj \in [\br]$, and similarly for $\tensor{A}_\bj'$. By $R$-iden\-ti\-fi\-a\-bi\-li\-ty, there is a permutation $\sigma \in \mathfrak{S}([\br])$ such that $\tensor{A}_{\bj} = \tensor{A}_{\sigma(\bj)}'$.
	
	Consider the HHTs $\tensor{H}$ and $\tensor{H}'$ {that correspond to} the ordered \ref{eqn_hhd}s $H=(\tensor{A}_{k1},\ldots,\tensor{A}_{k r_k})_{k=1}^m$ and $H'=(\tensor{A}_{k1}',\ldots,\tensor{A}_{k r_k}')_{k=1}^m$, respectively. In their natural flattenings, the existence of $\sigma$ implies that
	\[
	 \overline{\tensor{H}}' = \Sigma \overline{\tensor{H}},
	\]
	where $\Sigma\in \RR^{R\times R}$ is the permutation matrix of $\sigma$ in the lexicographic order on $[\br]$.

 	Through the lexicographic identification $\rho_{{\rm lex}} : [R] \simeq [r_1]\times\cdots\times[r_m]$, whose inverse is the flattening map $\mathrm{flatt}_{I}$ with $I=(1,\ldots,m)$, we can interpret $\Sigma$ as a rank-$1$ permutation $\rho$ of the columns of $\overline{\tensor{H}}$, i.e., $\tensor{H}_{:,\bj} \in \Var{A}^\rho$ for all $\bj\in[\bn]$.
	To be precise, we have the following commuting diagram:
	\[
		\xymatrix{
			\Bbbk^R \ar[r]^{\Sigma} \ar[d]^{\rho_{{\rm lex}}} \ar[rd]^\rho & \Bbbk^{R} \ar[d]^{\rho_{{\rm lex}}} \\
			\Bbbk^{r_1\times\cdots\times r_m}  \ar[r]_\sigma & \Bbbk^{r_1\times\cdots\times r_m}
		}
	\]
	Hence, we have the two irreducible varieties $\Var{A}^{\rho_{\rm lex}}$ and $\Var{A}^\rho$. By construction, all columns of $\overline{\tensor{H}}$ belong to the intersection $\Var{A}^{\rho_{\rm lex}}$ and $\Var{A}^\rho$. By the genericity of $\tensor{T}$, we have that the $\tensor{H}_{:,\bj}$'s are generic in $\Var{A}^{\rho_{\rm lex}}$. Hence, by irreducibility, $\Var{A}^{\rho_{\rm lex}} \subset \Var{A}^\rho$. \Cref{prop_allrk1perms} then tells us that $\sigma = \rho \circ \rho_{\rm lex}^{-1}$ is a rank-$1$ preserver permutation on $\Bbbk^{r_1\times\cdots\times r_m}$. Consequently, by \cref{prop_rk1presperm}, $\sigma = (\pi_1 \otimes\cdots\otimes \pi_m) \circ \psi$ for some flattening $\psi : [\br] \to [\br]$ and $\pi_k \in \mathfrak{S}_{r_k}$.
	Such a flattening $\psi$ corresponds to an element of $\mathfrak{P}(\br)$. Indeed, every $\varsigma \in \mathfrak{S}_m$ induces an $(I_1,\ldots,I_m)$-flattening in which each $I_i$ is a singleton and $\varsigma \in \mathfrak{P}(\br)$ if and only if it defines an automorphism on $\Bbbk^{r_1\times\cdots\times r_m}$.
	
	In conclusion, $\sigma$ is not just any permutation taking $\tensor{A}_\bj$ to $\tensor{A}_{\sigma(\bj)}'$, but rather it acts by (i) reordering the Hadamard factors of equal rank via $\psi$, and (ii) by reordering the rank-$1$ summands in the Hitchcock decomposition of each Hadamard factor via the $\pi_i$'s. Since essential \ref{eqn_hhd}s are invariant under these actions, we can assume without loss of generality that the tensors $\tensor{A}_{k i_k}'$ were ordered so that $\sigma$ is the identity. 

	By genericity of $\tensor{T}$, we can assume that $H$ is a generic ordered \ref{eqn_hhd} in the sense that it satisfies the assumption of \cref{prop_stacked}, i.e., there exists a $\bj \in [\br]$ such that $\tensor{A}_\bj$ has nonzero entries. Then, $H' \in [H]_{\Var{S}^*}$. That is, all of the ordered \ref{eqn_hhd}s of $\tensor{T}$ are equal up to the actions of $\mathfrak{P}(\br)$, the $\mathfrak{S}_{r_k}$'s, and $\Var{S}^*$. Hence, $\Var{H}_\br^{\Var{S}^*}(\tensor{T})$ is a singleton, concluding the proof.
\end{proof}

While this result shows that $\br$-identifiability is implied by $R$-identifiability for generic \ref{eqn_hhd}s, we recall from \cref{ssec_hhd_cpd} that there is no pressing reason why the very special Hitchcock decompositions induced from \ref{eqn_hhd}s should be identifiable. The next result shows that generic \ref{eqn_hhd}s are $R$-identifiable when $R$ is small.

\begin{proposition}\label{prop_hhdkruskal}
 Let $\br\in\NN^m$ and $R=r_1\cdots r_m$. 
 Let $\tensor{T}\in\sigma_\br(\Var{S})\subset\Bbbk^{n_1\times\cdots\times n_d}$ be generic. If there exists a partition $I \sqcup J \sqcup K = [d]$ such that $R$ satisfies \cref{eqn_kruskal_range}, then $\tensor{T}$ is $R$-identifiable.
\end{proposition}
\begin{proof}
By \cite[Lemma 4.4]{COV17}, a generic tuple of rank-$1$ tensors $A \in \Var{S}^{(R)}$ has maximal Kruskal ranks, i.e., 
\begin{equation}\label{eq:general_Kranks}
\mathrm{krank}_X(A) = \min\left\{R, \prod_{x\in X} n_x\right\}, \quad \forall X\subset [d].
\end{equation}
Namely, there exists a closed algebraic subvariety $\Var{Z} \subsetneq \Var{S}^{(R)}$ such that \cref{eq:general_Kranks} is satisfied for all $A \not\in \Var{Z}$. In particular, \cref{eq:general_Kranks} is satisfied for $X\in\{I,J,K\}$ for all $A \not\in \Var{Z}$. Note that the proof of \cite[Lemma 4.4]{COV17} holds for arbitrary fields of characteristic zero (dropping the statement about Euclidean denseness) as the polynomials have coefficients in $\mathbb{Q}$.

Since Kruskal's theorem holds for arbitrary fields, as proved by Rhodes \cite{Rhodes2010}, it follows from \cite[Theorem 4.6]{COV17} that a tensor $\tensor{T}$ is $R$-identifiable if \cref{eqn_kruskal_range} holds and it admits a CPD $A \subset \Var{S}^{(R)}$ with maximal Kruskal ranks $\mathrm{krank}_X(A)$ for $X\in\{I,J,K\}$.

Therefore, since the numerical condition given by \cref{eqn_kruskal_range} is satisfied by assumption, we just need to show that the CPD induced by a generic \ref{eqn_hhd}, say $\tensor{T} = \hadamard_{k=1}^m \sum_{i_k=1}^{r_k} \tensor{A}_{k i_k}$ and
\[
 A = \{ \tensor{A}_{1 i_1}\hadamard\cdots\hadamard\tensor{A}_{m i_m} : \bi \in [\br] \},
\]
is generic in the sense that $A \not\in \Var{Z}$. In other words, if $\Var{U} \subset \sigma_R(\Var{S})$ is the non-empty Zariksi open subset of rank-$R$ tensors with a CPD $A \not\in \Var{Z}$, then we need to show that $\Var{U} \cap \sigma_\br(\Var{S})\neq\emptyset$.

Next, we use an idea inspired by \cite[Proof of Lemma 13]{AMR2009}.
Let $X \subset [d]$ be any nonempty subset.
Without loss of generality, we can relabel the indices so that $X = \{1, \ldots, e\}$.
Let $p_{k,i_k} \in \NN$ denote the $(i_k + \sum_{\ell=0}^{e-1} r_\ell)$th prime number, so that 
\[
	\left\{ p_{k,i_k} : k\in[e], i_k\in[r_k] \right\}
\]
is the set of the first $r_1+\cdots+r_e$ consecutive prime numbers.
Then, we define the following rank-$1$ tensors
\[
\tensor{A}_{k,i_k}^X
:= 
\begin{bmatrix}
1 \\ 
p_{k,i_k} \\
p_{k,i_k}^2 \\ 
\vdots \\
p_{k,i_k}^{n_1-1}
\end{bmatrix}
\otimes 
\begin{bmatrix}
1 \\
p_{k,i_k}^{n_1}\\
p_{k,i_k}^{2 n_1}\\
\vdots \\
p_{k,i_k}^{(n_2-1)n_1}
\end{bmatrix}
\otimes 
\cdots 
\otimes 
\begin{bmatrix}
 1 \\
 p_{k,i_k}^{n_1+\cdots+n_{e-1}} \\
 p_{k,i_k}^{2(n_1+\cdots+n_{e-1})} \\
 \vdots \\
 p_{k,i_k}^{(n_k-1)(n_1+\cdots+n_{e-1})}
\end{bmatrix}.
\]
These rank-$1$ tensors are constructed so that the entries of $\tensor{A}_{k,i_k}^X$ are all consecutive powers of $p_{k,i_k}$. That is, if $E=(1,\ldots,e)$, then the $E$-flattening satisfies
\[
 (\tensor{A}_{k,i_k}^X)_{E}  = (1, p_{k,i_k}, p_{k,i_k}^2, \ldots, p_{k,i_k}^{N_e-1})^T,
\]
where $N_e=\prod_{i=1}^e n_i$. Moreover, as each $\tensor{A}_{k,i_k}^X$ is built from its own unique prime number $p_{k,i_k}$, we have 
\[
 ( \tensor{A}_{1,i_1}^X \hadamard\cdots\hadamard\tensor{A}_{m,i_m}^X )_E = 
 ( \tensor{A}_{1,i_1}^X )_E \hadamard\cdots\hadamard  (\tensor{A}_{m,i_m}^X )_E = 
 \begin{bmatrix}
 1 \\ 
 p_{1,i_1} p_{2,i_2} \cdots p_{m,i_m} \\
 p_{1,i_1}^2 p_{2,i_2}^2 \cdots p_{m,i_m}^2 \\ 
 \vdots \\
  p_{1,i_1}^{N_e-1} p_{2,i_2}^{N_e-1} \cdots p_{m,i_m}^{N_e-1}
 \end{bmatrix}.
\]
In other words, the rank-$1$ tensor $\tensor{A}^X_{i_1,\ldots,i_m} := \tensor{A}_{1,i_1}^X\hadamard\cdots\hadamard\tensor{A}_{m,i_m}^X \in \NN^{n_1\times\cdots\times n_e}$ has the consecutive powers of $q_{i_1,\ldots,i_m} = p_{1,i_1}p_{2,i_2}\cdots p_{m,i_m}$ as entries. By construction, all $q_{i_1,\ldots,i_m}$ are distinct numbers as their prime factorizations differ. Hence, it follows from a basic property of Vandermonde matrices that the set $A^X = \{ \tensor{A}_{i_1,\ldots,i_m}^X : \bi\in[\br] \}$ has maximal Kruskal rank equal to $\min\{R, N_e\}$, indeed, for every subset $S \subset A^X$ of cardinality $j$ we have that $\dim {\rm span}(S) = \min\{j,N_e\}$.

Finally, consider the rank-$1$ tensors
\[
 \tensor{A}_{k,i_k} := \tensor{A}_{k,i_k}^I \otimes \tensor{A}_{k,i_k}^J \otimes \tensor{A}_{k,i_k}^K.
\]
It follows from the previous argument that $A$ has maximal Kruskal ranks $\mathrm{krank}_X(A)$ for $X=I,J,K$. As $R$ lies in the range of applicability of \cref{lem_reshaped_kruskal}, the CPD induced from the \ref{eqn_hhd} $[(\{\tensor{A}_{k1},\ldots,\tensor{A}_{k r_k}\})_{k=1}^m]$ is $R$-identifiable. Hence, $A \in \Var{U}$.
\end{proof}

Combining the foregoing two results yields the next result.
\begin{corollary}\label{cor:identifiability}
Let $\br \in \NN^m$. If $R=r_1\cdots r_m$ lies in the range of applicability of \cref{lem_reshaped_kruskal}, then {the generic $\tensor{T}\in {\sigma_\br(\Var{S})} \subset \Bbbk^{n_1\times\cdots\times n_d}$ is $\br$-identifiable}.
\end{corollary}

As byproduct, we obtain a non-defectivity result on varieties parameterized by \ref{eqn_hhd}s. 

\begin{corollary}\label{cor:dim_RBM}
	Let $\br \in \NN^m$ and $\bn \in \NN^d$ such that $R$ lies in the range of applicability of \Cref{lem_reshaped_kruskal}. Then, $\sigma_\br(\Var{S}_\bn)$ has the expected dimension, i.e., \cref{eqn:upperbound_dim_HHD} is an equality.
\end{corollary}
\begin{proof}
	For a generic $\tensor{T} \in \sigma^\circ_\br(\Var{S})$ and an equivalence class $[(\{\tensor{A}_{k,1},\ldots,\tensor{A}_{k,r_k}\})_{k=1}^m]_{\Var{S}^*} \in \Var{H}_\br^{\Var{S}^*}(\tensor{T})$ of \ref{eqn_hhd}s, we consider as representative the \ref{eqn_hhd}
\[
	(\{\mathds{1},\tensor{A}_{1,2}',\ldots,\tensor{A}'_{1,r_k}\},\cdots,\{\mathds{1},\tensor{A}'_{m-1,2},\ldots,\tensor{A}'_{m-1,r_{m-1}}\},\{\tensor{A}'_{m,1},\ldots,\tensor{A}'_{m,r_{m-1}}\}),
\] 
obtained by Hadamard multiplying the first $m-1$ Hadamard factors by $\tensor{A}_{k,1}^\ominus$, for $k \in [m-1]$, and the last Hadamard factor by $\tensor{A}_{1,1}\hadamard\cdots\hadamard\tensor{A}_{m-1,1}$. As $\br$ lies in the range of applicability of \cref{lem_reshaped_kruskal} by assumption, \cref{cor:identifiability} tells us the map
\begin{equation*}
	\begin{array}{c c c c}
		 & \sigma^\circ_{r_1-1}(\Var{S}) \times\cdots\times \sigma^\circ_{r_{m-1}-1}(\Var{S}) \times  \sigma_{r_m}^\circ(\Var{S}) & \longrightarrow &  \Bbbk^{n_1\times\ldots\times n_d} \\
		 & \left(\tensor{T}'_1,\ldots,\tensor{T}'_{m-1},\tensor{T}_m\right)& \longmapsto & (\mathds{1}+\tensor{T}'_1) \hadamard \cdots \hadamard (\mathds{1}+\tensor{T}'_{m-1}) \hadamard \tensor{T}_m
	\end{array}
\end{equation*}
is $1$-to-$1$. This concludes the proof.
\end{proof}

\section{A numerical algorithm for computing low-rank HHDs}\label{sec:algorithm}

Having studied $\br$-identifiability in the previous section, we continue by proposing an algorithm for computing the essential \ref{eqn_hhd} of an $\br$-identifiable tensor $\tensor{T}$.
This algorithm fuses the three key ideas of this paper: (i) the relation between an ordered \ref{eqn_hhd} $H$ and its induced CPD, (ii) rank-$1$ permutations, and (iii) the HHT.

\subsection{{Basic version}}
Assume we are given a generic, $\br$-identifiable, $R$-identifiable tensor $\tensor{T}\in\Bbbk^{n_1\times\cdots\times n_d}$.
That is, $\tensor{T}$ has the decompositions
\begin{align}\label{eqn_hhd_cpd}
 \tensor{T} = \overset{m}{\underset{k=1}{\bighadamard}} \sum_{i_k=1}^{r_k} \tensor{A}_{k i_k} = \sum_{i=1}^R \tensor{B}_i,
\end{align}
where $H:= (H_1,\ldots,H_m) :=(\tensor{A}_{k1},\ldots,\tensor{A}_{k r_k})_{k=1}^m$ is an ordered \ref{eqn_hhd} and the set $\{\tensor{B}_1,\ldots,\tensor{B}_R\}$ is the CPD of $\tensor{T}$. Our goal is to deduce the $\tensor{A}_{k i_k}$'s by knowing the $\tensor{B}_i$'s. By $R$-identifiability, the $\tensor{B}_i$'s correspond, up to reordering, to the CPD induced by the \ref{eqn_hhd} in \eqref{eqn_hhd_cpd}. The correct reordering, unknown a priori, is obtained by using an appropriate rank-$1$ permutation, as shown shortly in \cref{lem_hhtfromrk1perm}. In this way, we construct the HHT of one of $\tensor{T}$'s ordered \ref{eqn_hhd}s from the CPD. From there, we can reconstruct the \ref{eqn_hhd} as explained after \cref{lem_hhtfromrk1perm}.

The idea of \cref{lem_hhtfromrk1perm} is the following. Fixing an entry $\bi \in [\bn]$ of the tensor, we look at the vector $\mathbf{b}_\bi = ( (\tensor{B}_1)_\bi, \ldots, (\tensor{B}_R)_\bi ) \in \Bbbk^R$. By assumption, there exists a permutation $\sigma \in \mathfrak{S}_R$ of the vector $\mathbf{a}_\bi = (\ldots, (\tensor{A}_{1j_1})_\bi \cdots (\tensor{A}_{kj_k})_\bi,\ldots)_{\bj \in [\br]} \in \Bbbk^R$ such that $\sigma^{-1}(\mathbf{a}_\bi) = \mathbf{b}_\bi$. Note that ${\bf a}_\bi \in \Var{A}^{\rho_{\rm lex}}$ by construction, where $\rho_{\rm lex} : [R] \to [\br]$ is the lexicographic identification $[R] \simeq [\br]$. Hence, $\mathbf{b}_\bi \in \Var{A}^{\rho_{\rm lex} \circ \sigma}$. In the following, we show that stacking the $\tensor{B}_i$'s following the order induced by the permutation $\rho$ gives us the HHT of one of $\tensor{T}$'s ordered \ref{eqn_hhd}s.

\begin{lemma}\label{lem_hhtfromrk1perm}
Let $\tensor{T}$ be a generic tensor as in \cref{eqn_hhd_cpd}.
Fix any $\bk \in [\bn]$.
Let $\rho : [R] \to [\br]$ be the unique (up to rank-$1$ preserver permutations) maximal rank-$1$ permutation of $\mathbf{b}_\bk = ( (\tensor{B}_1)_\bk, \ldots, (\tensor{B}_R)_\bk )$. Then, up to rank-$1$ preserver permutations, $\rho$ is the unique, maximal, simultaneous rank-$1$ permutation of $\mathbf{b}_\bi$ for all $\bi\in[\bn]$. Moreover, the tensor $\tensor{H} \in \Bbbk^{r_1\times\cdots\times r_m\times n_1\times\cdots\times n_d}$ given by
\[
 \tensor{H}_{\rho(i),:} = \tensor{B}_i, \quad i\in[R],
\]
is the HHT of one of the ordered \ref{eqn_hhd}s of $\tensor{T}$.
\end{lemma}
\begin{proof}
Let $H'=(\tensor{A}_{k 1}, \ldots, \tensor{A}_{k r_k})_{k=1}^m$ denote an ordered \ref{eqn_hhd} of $\tensor{T}$, and let $\tensor{H}'$ denote its corresponding HHT.

By $R$-identifiability, there exists a bijection $\rho' : [R] \to [\br]$ such that we have $\tensor{B}_i =  \tensor{A}_{1 i_1} \hadamard\cdots\hadamard \tensor{A}_{m i_m}$ where $\bi = \rho'(i)$, for all $i\in[R]$. In other words, $\mathbf{b}_\bj \in \Var{A}^{\rho'}$, for \emph{all} $\bj \in [\bn]$. In particular, $\mathbf{b}_\bk \in \Var{A}^{\rho'} \cap \Var{A}^{\rho}$. This intersection cannot be a strict subvariety of $\Var{A}^\rho$, by the genericity assumption of uniqueness and maximality of $\rho$ for $\mathbf{b}_\bk$. Hence, we deduce that $\Var{A}^{\rho'} \cap \Var{A}^{\rho} = \Var{A}^\rho$, i.e., $\Var{A}^{\rho} \subset \Var{A}^{\rho'}$. By irreducibility and equidimensionality, we have that $\Var{A}^{\rho} = \Var{A}^{\rho'}$. Then, by \cref{prop_allrk1perms}, up to a rank-$1$ preserver permutation $\sigma$, $\rho'$ corresponds to the rank-$1$ permutation $\rho$. That is, $\rho' = \sigma \circ \rho$. By symmetry of the proof, we can conclude that $\rho'$ is, up to rank-$1$ preserver permutations, the maximal and unique \emph{simultaneous} rank-$1$ permutation of all $\mathbf{b}_\bi$'s.

As used in the proof of \cref{thm:main_identifiabilityHHD}, every rank-$1$ preserver permutation in $\mathfrak{S}([\br])$ has a natural action on an ordered \ref{eqn_hhd}. Indeed, by \cref{prop_rk1presperm}, $\sigma = (\pi_1 \otimes \cdots \otimes \pi_m) \circ \psi$ where $\pi_i \in \mathfrak{S}_{r_i}$ and $\psi \in \mathfrak{P}(\br)$. Therefore, $\sigma$ can act on ordered \ref{eqn_hhd}s by (i) reordering Hadamard factors with equal ranks by $\psi$ and (ii) reordering the rank-$1$ summands of each Hadamard factor by the $\pi_i$'s. Let $H'' = \sigma(H')$ denote the resulting ordered \ref{eqn_hhd} and $\tensor{H}''$ its HHT.

By \cref{def_stacked}(2), we then have for all $\bi\in[\br]$ that 
\(
 \tensor{H}_{\bi,:}''= \tensor{H}_{\sigma(\bi),:}'.
\)
Hence, for all $i\in[R]$, we find
	\[
	 \tensor{B}_i = \tensor{A}_{1 i_1}\hadamard\cdots\hadamard\tensor{A}_{m i_m} = \tensor{H}_{\rho'(i),:}' = \tensor{H}_{\sigma(\rho(i)),:}' = \tensor{H}_{\rho(i),:}''.
	\]
	This shows that $\tensor{H}''=\tensor{H}$, and, hence, the latter is the HHT of $H''$.
\end{proof}

From the last lemma, we know how to construct the HHT of an \ref{eqn_hhd} from a CPD of a generic tensor. Now, we want to recover the rank-$1$ elements of the \ref{eqn_hhd} from this HHT. By \cref{def_stacked}(3), for all $\bj\in[\bn]$, we have that $\tensor{H}_{:,\bj}$ is a rank-$1$ tensor realized as a tensor product of vectors containing the $\bj$th element of the unknown rank-$1$ tensors $\tensor{A}_{k i_k}$ forming the \ref{eqn_hhd}:
\[
\tensor{H}_{:,\bj} =
\begin{bmatrix}
 (\tensor{A}_{11})_\bj\\
 \vdots\\
 (\tensor{A}_{1 r_1})_\bj
\end{bmatrix}
\otimes\cdots\otimes
\begin{bmatrix}
 (\tensor{A}_{m1})_\bj\\
 \vdots\\
 (\tensor{A}_{m r_m})_\bj
\end{bmatrix} \in \Var{S}.
\]
The vectors in a rank-$1$ tensor are identifiable only up to scale. Nevertheless, by genericity, we can assume there is an index $\bi\in[\br]$ such that $\tensor{A}_{1 i_1}\hadamard\cdots\hadamard\tensor{A}_{m i_m} \in \Var{S}^*$, so that each $\tensor{A}_{k i_k} \in \Var{S}^*$.
For notational convenience, we can assume without loss of generality that $\bi={\bf1}$.

We denote $\tensor{A}\hadamard(\tensor{T}_1,\ldots,\tensor{T}_k) := (\tensor{A}\hadamard\tensor{T}_1, \ldots, \tensor{A}\hadamard\tensor{T}_k)$.
Then, by \cref{prop_stacked}, $\tensor{H}$ is also the HHT of
\begin{align*}
H' &= ( (\tensor{A}_{21}\hadamard\cdots\hadamard\tensor{A}_{m1}) \hadamard H_1, \tensor{A}_{21}^\ominus \hadamard H_2, \ldots, \tensor{A}_{m1}^\ominus \hadamard H_m) \\
&= (\tensor{A}_{k1}', \ldots, \tensor{A}_{k r_k}')_{k=1}^m,
\end{align*}
{which lives in $[H]_{\Var{S}^*}$.}
It follows that we have, for all $\bj\in[\br]$,
\begin{equation}\label{eq:interpolatory_HHT}
\tensor{H}_{:,\bj} =
\begin{bmatrix}
 (\tensor{A}_{11})_\bj\\
 \vdots\\
 (\tensor{A}_{1 r_1})_\bj
\end{bmatrix}
\otimes\cdots\otimes
\begin{bmatrix}
 (\tensor{A}_{m1})_\bj\\
 \vdots\\
 (\tensor{A}_{m r_m})_\bj
\end{bmatrix}
=
\begin{bmatrix}
 (\tensor{A}_{11}')_\bj\\
 \vdots\\
 (\tensor{A}_{1 r_1}')_\bj
\end{bmatrix}
\otimes
\begin{bmatrix}
 1\\
 (\tensor{A}_{2 2}')_\bj\\
 \vdots\\
 (\tensor{A}_{2 r_1}')_\bj
\end{bmatrix}
\otimes\cdots\otimes
\begin{bmatrix}
1\\
(\tensor{A}_{m 2}')_\bj\\
 \vdots\\
 (\tensor{A}_{m r_m}')_\bj
\end{bmatrix}.
\end{equation}
Given $\mathbf{a} \in \Bbbk^{n-1}$, for any $i \in [n]$, denote $\iota_i(\mathbf{a}) \in \Bbbk^n$ the vector 
\[
	\iota_i(\mathbf{a}) := (a_1,\ldots,a_{i-1},1,a_i,\ldots,a_{n-1}).
\]
It is well known and straightforward to verify that, for every $\bi \in [(n_2,\ldots,n_d)]$,
\begin{equation*}
 \phi_{\bi} : \Bbbk^{n_1} \times \Bbbk^{n_2-1} \times\cdots\times \Bbbk^{n_d-1} \to \Var{S}_\bn,\, (\ba_1, \ba_2, \ldots, \ba_d) \mapsto \ba_1\otimes \iota_{i_1}(\ba_2) \otimes\cdots\otimes \iota_{i_d}(\ba_d)
\end{equation*}
is a bijective polynomial map onto its image $\mathrm{Im}(\phi_\bi)\subsetneq\Var{S}_\bn$, and that the latter is a Zariski open subset of $\Var{S}_\bn$.
Note that $\phi_{\bi}^{-1}(\tensor{A})$ essentially computes the interpolatory decomposition of a rank-$1$ tensor $\tensor{A} \in \mathrm{Im}(\phi_\bi)$ as in \cref{rmk:char_rank1}.

In our case, we compute $\phi_{{\mathbf{1}}}^{-1}(\tensor{H}_{:,\bj})$ to recover the elements of the $\bj$th entry of the rank-$1$ elements forming the \ref{eqn_hhd} that defines the HHT $\tensor{H}$.

\begin{algorithm}[tb]\small
\begin{algorithmic}[1]
\Require $\tensor{T} \in \Bbbk^{n_1\times\cdots\times n_d}$ is a generic tensor admitting an \ref{eqn_hhd}.
\vspace{4pt}
\Statex\Comment{Perform a CPD of $\tensor{T}$.}
\State $\{ \tensor{B}_1, \ldots, \tensor{B}_R \} \leftarrow \texttt{CPD}(\tensor{T}, R)$; \label{line_cpdBasic}
\vspace{4pt}
\Statex\Comment{Determine a rank-$1$ permutation of the ${\bf1}$-entries.}
\State $\mathbf{b}_{\bf1} \leftarrow \bigl( (\tensor{B}_1)_{\bf1}, \ldots, (\tensor{B}_R)_{\bf1} \bigr)$;\label{line_rk1permstart}
\State $\Var{M} \leftarrow \text{\cref{alg_basic_minors}}(\mathbf{b}_{\bf1})$;\label{line_computeminors}
\State $\rho \leftarrow \text{\cref{alg_rank1_permutation}}(\mathbf{b}_{\bf1}, \Var{M})$;\label{line_rk1permstop}
\vspace{4pt}
\Statex \Comment{Arrange the rank-$1$ tensors in an HHT and extract $\tensor{A}_{k i_k}$'s coordinates from it.}
\State $\tensor{H}_{\rho(i),:} \leftarrow \tensor{B}_{i}, \quad i\in[R]$;\label{line_extractstart}
\For {$\bj \in [\bn]$}\label{line_extractstart2}
\State $(\ba_1, \ldots, \ba_m) \leftarrow \phi_{\bf1}^{-1}( \tensor{H}_{:,\bj} )$;
\State $( (\tensor{A}_{11})_\bj,\ldots,(\tensor{A}_{1r_1})_\bj ) \leftarrow \ba_1$;
\State $( (\tensor{A}_{k2})_\bj,\ldots,(\tensor{A}_{kr_k})_\bj ) \leftarrow \ba_k, \quad k=2,\ldots,m$;
\EndFor
\State $\tensor{A}_{k1} \leftarrow \mathds{1}, \quad k=2,\ldots,m$;\label{line_extractstop}
\State \Return $H = ( \tensor{A}_{k1},\ldots,\tensor{A}_{k r_k} )_{k=1}^m$;
\end{algorithmic}
 \caption{{Theoretical computation} of the component CPDs of a generic tensor admitting an \ref{eqn_hhd}, {not robust in floating-point arithmetic}.}
 \label{alg_hhd_basic}
\end{algorithm}

The foregoing discussion leads to \cref{alg_hhd_basic}. The main steps are as follows:
\begin{itemize}[nolistsep]
	\item line \ref{line_cpdBasic}: compute the rank-$R$ CPD of the input tensor $\tensor{T}$;
	\item lines \ref{line_rk1permstart}--\ref{line_rk1permstop}: from the ${\bf1}$-entries of the summands of the CPD, a rank-$1$ permutation is computed using \cref{alg_basic_minors,alg_rank1_permutation}; 
	\item line \ref{line_extractstart}: build the HHT as suggested by \cref{lem_hhtfromrk1perm};
	\item lines \ref{line_extractstart2}--\ref{line_extractstop}: extract the coordinates of the rank-$1$ tensors in the  \ref{eqn_hhd}.
\end{itemize}
It applies to \emph{generic} tensors that are $\br$-identifiable and $R$-identifiable. By carefully inspecting where genericity is invoked, we can be more precise.

\begin{proposition}
 If a tensor $\tensor{T}$
 \begin{enumerate}[nolistsep]
  \item is $\br$-identifiable,
  \item is $R$-identifiable,
  \item has a $\tensor{B}_i \in \Var{S}^*$ in its CPD $\{\tensor{B}_1,\ldots,\tensor{B}_R\}$, and
  \item has an $\bi\in[\bn]$ such that the admissible vector $\mathbf{b}_\bi$ has distinct entries and a unique maximal rank-$1$ permutation,
 \end{enumerate}
 then \cref{alg_hhd_basic} recovers its essential \ref{eqn_hhd}.
\end{proposition}

\subsection{Advanced version}\label{sec_sub_advanced_version}
{Two modifications of \cref{alg_hhd_basic} are recommended for efficient, numerical implementations.}

{The first modification concerns} the observation that one is typically not interested in the coordinates of $\tensor{A}_{k i_k} \in \Var{S}$ with respect to a basis of $\Bbbk^{n_1\times\cdots\times n_d}$ {because} the dimension of $\Var{S}_{\bn}$ is {only} $\dim\Var{S}_{\bn} \approx n_1+\cdots+n_d \ll n_1\cdots n_d = N$. In a practical implementation, one would simply return any set of vectors $\ba_{k \ell_k}^1, \ldots, \ba_{k \ell_k}^d$ such that $\tensor{A}_{k \ell_k} = \ba_{k \ell_k}^1 \otimes \cdots \otimes \ba_{k \ell_k}^d$.
Note that if we know the values of a rank-$1$ tensor $\tensor{A}=\ba_1\otimes\cdots\otimes\ba_d$ at the $\bi$-cross
\(
 \Var{C}_\bi(\tensor{A}) = \Var{C}_{\bi}^1(\tensor{A}) \cup\cdots\cup \Var{C}_\bi^d(\tensor{A})
\)
from \cref{def_cross},
then we can take $\ba_1 = \Var{C}_\bi^1(\tensor{A})$, $\ba_k = \tensor{A}_{\bi}^{-1} \Var{C}_\bi^k(\tensor{A})$ for $k=2,\ldots,d$, by \cref{rmk:char_rank1}, as factorization of $\tensor{A}$.
Let $\Var{C}_{\bf1}(\br)$ denote the set of indices in the ${\bf1}$-cross of an $r_1\times\cdots\times r_m$ tensor, and likewise for $\Var{C}_{\bf1}(\bn)$.
Then, based on the foregoing insight, we recommend to replace lines \ref{line_extractstart}--\ref{line_extractstop} of the basic \cref{alg_hhd_basic} by the equivalent but much more efficient computation in \cref{alg_fast}.

\begin{algorithm}[tb]\small
\begin{algorithmic}[1]
\Require $\tensor{B}_i = \mathbf{b}_i^1\otimes\cdots\otimes\mathbf{b}_i^d$ are rank-$1$ tensors, given in a data-sparse format such as $(\mathbf{b}_i^1,\ldots,\mathbf{b}_i^d)$.
\Require $\rho : [R]\to[\br]$ is a bijection.
\For{ $\bj \in \Var{C}_{\bf1}(\bn)$ }
\State $\Var{C}_{\bf1}(\tensor{H}_{:,\bj}) \leftarrow \bigl( (\tensor{B}_{\rho^{-1}(\bi)})_\bj ~:~ \bi \in \Var{C}_{\bf1}(\br) \bigr)$;
\State $(\ba_1, \ldots,\ba_m) \leftarrow \bigl( \Var{C}_{\bf1}^1(\tensor{H}_{:,\bj}), \tensor{H}_{{\bf1},\bj}^{-1} \Var{C}_{\bf1}^2(\tensor{H}_{:,\bj}), \ldots, \tensor{H}_{{\bf1},\bj}^{-1} \Var{C}_{\bf1}^m(\tensor{H}_{:,\bj}) \bigr) $;
\State $( (\tensor{A}_{k1})_\bj,\ldots,(\tensor{A}_{kr_k})_\bj ) \leftarrow \ba_k, \quad k=1,\ldots,m$;
\EndFor
\For{ $k=1,\ldots,m$ }
\For{ $i_k=1,\ldots,r_k$}
\State $(\ba_{k i_k}^1, \ldots, \ba_{k i_k}^d) \leftarrow \bigl( \Var{C}_{\bf1}^1(\tensor{A}_{k i_k}), (\tensor{A}_{k i_k})_{\bf1}^{-1} \Var{C}_{\bf1}^2(\tensor{A}_{k i_k}), \ldots, (\tensor{A}_{k i_k})_{\bf1}^{-1} \Var{C}_{\bf1}^d(\tensor{A}_{k i_k}) \bigr)$;
\EndFor
\EndFor
\State\Return $(\ba_{k i_k}^1,\ldots,\ba_{k i_k}^d)$ for all $k\in[m]$ and $i_k\in[r_k]$.
\end{algorithmic}
\caption{Fast decoupling of the CPD into the HHD components.}\label{alg_fast}
\end{algorithm}

{The second modification mitigates the impact of numerical round-off errors and noise. \Cref{alg_hhd_basic} is designed to recover the constituent rank-$1$ components of a tensor admitting an exact Hadamard--Hitchcock decomposition as in \cref{eqn_hhd} using floating-point arithmetic. While it can tolerate small round-off errors, it is not specifically designed to deal with \emph{approximation problems} in which the model \cref{eqn_hhd} holds only approximately. Nevertheless, \cref{alg_hhd_basic} can be made more robust to small model violations by further exploiting \cref{lem_hhtfromrk1perm}. It showed that for every choice of index $\vect{k}\in [\vect{n}]$, the unique maximal rank-$1$ permutation $\rho$ of $\mathbf{b}_{\bf k} = \bigl( (\tensor{B}_1)_{\bf k}, \ldots, (\tensor{B}_R)_{\bf k} \bigr)$ simultaneously permutes $\vect{b}_{\vect{i}}$ for \emph{all} $\vect{i}\in[\vect{n}]$ into rank-$1$ tensors. This means that in \cref{alg_basic_minors}, for a generic $\tensor{T}\in\sigma_{\vect{r}}(\Var{S})$, $(i,j,k) \in \Var{M}(\vect{b}_\vect{k})$ is a vanishing top-left corner minor for a fixed $\vect{b}_\vect{k}$ if and only if it is a vanishing minor for \emph{all} $\vect{b}_\vect{i}$'s. While it theoretically suffices to use a single $\vect{b}_\vect{k}$, as in step 3 of \cref{alg_hhd_basic}, numerically it will be much more robust to use multiple indices $\vect{k}_1, \dots, \vect{k}_\ell \in [\vect{n}]$ and construct an approximate simultaneous rank-$1$ permutation from the corresponding vectors $\{ \vect{b}_{\vect{k}_1}, \dots \vect{b}_{\vect{k}_\ell} \}$.}

{The following approach exploits the above observation.
We seek top-left corner minors $m_{i,j,k}$ such that $(\vect{b}_{\vect{k}_l})_1 (\vect{b}_{\vect{k}_l})_k \approx (\vect{b}_{\vect{k}_l})_i (\vect{b}_{\vect{k}_l})_j$ simultaneously holds for all $1 \le l \le \ell$. Let $B$ be the matrix that has the vectors $\vect{b}_{\vect{k}_i}$ as columns.
Let $C \in \Bbbk^{R \times \ell}$ be the matrix whose $k$th row is 
\[
C_{k,:} 
= \bigl( (\vect{b}_{\vect{k}_1})_1 (\vect{b}_{\vect{k}_1})_k, \dots, (\vect{b}_{\vect{k}_\ell})_1 (\vect{b}_{\vect{k}_\ell})_k \bigr)
= B_{1,:} \hadamard B_{k,:}.
\]
If $2 \le i < j \le R$ are fixed, then finding a top-left corner minor $m_{i,j,k}$ that approximately vanishes amounts to solving 
\begin{align}\label{eqn_lin_search}
\underset{2 \le k \le R}{\arg\min} \| C_{k,:} - B_{i,:}\hadamard B_{j,:} \|
\end{align}
and verifying that it is sufficiently close to zero; herein, $\|\cdot\|$ denotes a norm, such as the Euclidean norm.
A linear search solves this problem using $\mathcal{O}( R\ell )$ operations for fixed $i,j$, leading to a total complexity of $\mathcal{O}(\ell R^3)$.}

{We propose to accelerate the solution of this minimization problem further with a \emph{locality sensitive hashing} scheme that performs an \emph{exact nearest neighbor search} when $\Bbbk\in\{\mathbb{Q},\mathbb{R},\mathbb{C}\}$. The approach is inspired by random projection forests \cite{RDF1995}. Consider a vector $\vect{n} \in \Bbbk^\ell$ of unit Euclidean norm, which defines the projection 
\[
\Pi_{\vect{n}} : \Bbbk^\ell \longrightarrow \mathbb{R},\quad \vect{x} \longmapsto \mathfrak{R}(\vect{x}^T \vect{n}),
\]
where $\mathfrak{R}(\cdot)$ takes the real part of its argument. Note that $\Pi_{\vect{n}}$ is \emph{contractive}, as $\| \Pi_{\vect{n}}(\vect{x}-\vect{y}) \| \le \|\vect{n}\| \|\vect{x}-\vect{y}\| = \|\vect{x}-\vect{y}\|$, where $\|\cdot\|$ is the Euclidean norm.
This implies in particular that if $C_{k,:}$ is close to $B_{i,:}\hadamard B_{j,:}$, then they will also be close after applying the random projection $\Pi_{\vect{n}}$. Since the row vectors of $C$ over which we search is a fixed set, we can precompute their projections $s_k := \Pi_{\vect{n}}(C_{k,:})$ and sort these real numbers. This \emph{$\vect{n}$-search index} is built only once. Then, the $\binom{R}{2}$ searches over all pairs $(i,j)$ can exploit this index by using a binary search to find all the consecutive indices $\vect{i}_{\vect{n}} \subset [R]$ for which the relative error is sufficiently small, i.e.,
\[
 \left| \frac{\Pi_{\vect{n}}(B_{i,:}\hadamard B_{j,:}) - s_k}{\Pi_{\vect{n}}(B_{i,:}\hadamard B_{j,:}) } \right| \le  \epsilon \sqrt{\ell},
\]
where $\epsilon\approx0$ is the numerical threshold. In this way, the domain of the optimization problem \cref{eqn_lin_search} can be reduced to search only over the indices $\vect{i}_{\vect{n}}$.
Multiple $\vect{n}_i$-search indices can be built for (orthogonal) unit vectors $\vect{n}_i$ to yield multiple reduced index sets $\vect{i}_{\vect{n}_i}$. Thusly, we can reduce the search in \cref{eqn_lin_search} to $\cap_i \vect{i}_{\vect{n}_i}$.}

\begin{algorithm}[tb]\small
\begin{algorithmic}[1]
\Require $B \in \Bbbk^{R \times \ell}$ is a matrix whose columns are generic admissible vectors in $\Bbbk^R$ for an unknown rank-$1$ permutation $\rho : R \to [\vect{r}]$.
\Require $\epsilon$ is the relative tolerance to identify vanishing minors, e.g., $\epsilon = 5\cdot10^{-3}$.
\vspace{4pt}
\State $C \leftarrow B \hadamard (\vect{1} B_{1,:})$; \label{line_rmstep1}
\vspace{4pt}\Statex\Comment{Sketch $C$ down to $\RR^{R\times2}$.}
\State Sample a random matrix $N=\begin{bmatrix}\vect{n}_1 & \vect{n}_2\end{bmatrix}\in\Bbbk^{\ell\times 2}$ with orthonormal columns;
\State $S \leftarrow \mathfrak{R}(C N)$;
\State Let $\vect{p}_i$ be a permutation so the $i$th column of $S$ is sorted increasingly, for $i=1,2$;
\vspace{4pt}\Statex\Comment{Determine the set of vanishing minors $\Var{M}$.}
\State $\Var{M} \leftarrow []$;
\For{$2 \le i < j \le R$}
	\State $\vect{b} \leftarrow B_{i,:} \hadamard B_{j,:}$;
	\State $(\beta_1, \beta_2) = \mathfrak{R}(N^T \vect{b})$;
	\State Use binary search on $S_{\vect{p}_i,i}$ to find all indices $\vect{i}_{\vect{n}_i}$ of $S_{:,i}$ close to $\beta_i$, for $i=1,2$;
	\State Let $\vect{i} = \vect{i}_{\vect{n}_1} \cap \vect{i}_{\vect{n}_2}$ be the common indices;
	\State Let
	\(
	 k \in \underset{k \in \vect{i}}{\arg\min} \| \vect{b} - C_{k,:} \|;
	\)
	\If{$\sharp \vect{i} > 0$ \and $k \not\in \{1, i, j\}$ \and $\| \vect{b} - C_{k,:} \| < \epsilon \|\vect{b}\|$}
 		\State $\Var{M} \leftarrow \Var{M} \cup \{ (i,j,k), (j,i,k) \}$;
	\EndIf
\EndFor
\State \Return $\Var{M}$;
\end{algorithmic}
 \caption{Numerical identification of approximately vanishing top-left corner minors from multiple generic admissible vectors, robust in floating-point arithmetic. }
 \label{alg_robust_minors}
\end{algorithm}

{The above discussion is summarized as \cref{alg_robust_minors}. Based our experiments, we suggest to take $\epsilon = 5\cdot 10^{-3}$ as default parameter choice. We also found that using $2$ search indices suffices for $R \le 30\,000$.}

{If $h$ denotes the average length of the identified indices $\vect{i}_{\vect{n}_i}$, $i=1,2$, then the cost of the above scheme amounts to $\mathcal{O}(2R \log_2 R)$ for constructing the two search indices, and $\mathcal{O}(\ell + \log_2 R + h)$ operations to project, search, and find the closest vector for a fixed $(i,j)$. The total cost of finding the approximately vanishing minors is thus reduced to $\mathcal{O}(R^2 (\ell + h + \log_2 R))$ operations (plus lower-order terms for the projection and sorting). As $h \le R$, the time complexity is no worse than $\mathcal{O}(R^3)$.

Our experiments indicated that by using only $2$ search indices, the minimization problem \cref{eqn_lin_search} is solved about one order of magnitude faster than all of the exact tree search methods (k-d trees, ball trees, and standard linear search) from \texttt{NearestNeighbors.jl} for $R \approx 1000$ and $\ell = 30$. For $R \le 1000$, the average number of indices in $\vect{i}_{\vect{n}_1} \cap \vect{i}_{\vect{n}_2}$ was less than $2$. This value does seem to slowly increase with $R$ from about $1$ for very small values of $R$, up to about $10$ for $R \approx 30\,000$. A full theoretical analysis of the performance is outside of the scope of this paper.}

\begin{algorithm}[tbp]\small
{\begin{algorithmic}[1]
\Require $\tensor{T} \in \Bbbk^{n_1\times\cdots\times n_d}$ is a generic tensor admitting an \ref{eqn_hhd}.
\Require $1 \le \ell \le N = n_1 \dots n_d$ is the number of samples from which to compute a rank-$1$ permutation.
\vspace{4pt}
\Statex\Comment{Perform a CPD of $\tensor{T}$.}
\State $\{ \tensor{B}_1, \ldots, \tensor{B}_R \} \leftarrow \texttt{CPD}(\tensor{T}, R)$; \label{line_cpd}
\vspace{4pt}
\Statex\Comment{Robustly determine a rank-$1$ permutation.}
\State Sample $\ell$ indices $\vect{k}_1, \dots, \vect{k}_\ell \in [\vect{n}]$ with uniform probability and without replacement;
\State $\vect{b}_{\vect{k}_i} \leftarrow \bigl( (\tensor{B}_1)_{\vect{k}_i}, \dots, (\tensor{B}_R)_{\vect{k}_i}  \bigr)$ for $i=1,2,\dots,\ell$;
\State $B \leftarrow \begin{bmatrix} \vect{b}_{\vect{k}_1} & \vect{b}_{\vect{k}_2} & \dots & \vect{b}_{\vect{k}_\ell} \end{bmatrix}$;
\State $\Var{M} \leftarrow \Cref{alg_robust_minors}(B, \epsilon)$;
\State $\rho \leftarrow \Cref{alg_rank1_permutation}(B_{:,1}, \Var{M})$;
\vspace{4pt}
\Statex \Comment{Arrange the rank-$1$ tensors in an HHT and extract $\tensor{A}_{k i_k}$'s coordinates from it.}
\State $A \leftarrow \Cref{alg_fast}((\tensor{B}_1, \dots, \tensor{B}_R), \rho)$;
\State \Return $A$;
\end{algorithmic}}
 \caption{Efficient numerical computation of the rank-$1$ components of a generic tensor admitting an \ref{eqn_hhd}.}
 \label{alg_hhd}
\end{algorithm}

{Based on the above discussion of the two modifications, we propose \cref{alg_hhd} as a practical numerical algorithm for computing \ref{eqn_hhd}s of generic tensors in which slight violations of the model \cref{eqn_hhd} can be tolerated.
Our experiments suggest that $\ell = \min\{50, N\}$ leads to a good tradeoff between robustness to noise and computational performance.}

\subsection{Time complexity}
The asymptotic time complexity of \cref{alg_hhd} is the sum of (i) the time complexity for computing a rank-$R$ CPD of an $n_1\times\cdots\times n_d$ tensor, (ii) {$\mathcal{O}( R^2 (h + \log_2 R) )$ for computing the minors} with \cref{alg_robust_minors}, {(iii) $\mathcal{O}(R^2 \log_2 R)$ for extracting the rank-$1$ permutation from the minors with \cref{alg_rank1_permutation}, and (iv)} the complexity of extracting the rank-$1$ tensors $\tensor{A}_{k i_k}$ from the CPD with \cref{alg_fast}, namely
\[
 \mathcal{O}( (n_1+\cdots+n_d)( r_1+\cdots+r_m )m )
\]
elementary operations.
Note the factor $m$ originates from evaluating the $\bj$-entries of the rank-$1$ tensors $\tensor{B}_i$ in line $2$ of \cref{alg_fast}, assuming they are also represented as tuples of vectors rather than in ambient coordinates (as required by \cref{alg_fast}).
The cost for computing a CPD dominates all other costs. For example, basic pencil-based algorithms \cite{Harshman1970,LRA1993} {reduce the input tensor to two $R \times R$ matrices and then compute their generalized eigendecomposition to extract one of the factor matrices. The standard QZ algorithm for this subproblem already has a time complexity of $\mathcal{O}(R^3)$ \cite[Section 7.7]{GvL4}.} 

\section{Numerical experiments}\label{sec:experiments}
In this section, we present numerical experiments to evaluate the performance of \cref{alg_rank1_permutation,alg_hhd,alg_robust_minors}. First, the {implementation details and} experimental setup are described in \cref{sec_sub_impl_det,sec_sub_generalexperiments} respectively. Then, in \cref{sec_sub_rk1permexp}, the computational performance of \cref{alg_rank1_permutation,alg_robust_minors} is investigated. Thereafter, we present several experiments in \cref{sec_sub_hhdexper,sec_sub_highorder,sec_sub_comparisonpba,sec_sub_noise_tolerance} evaluating the computational performance and numerical accuracy of \cref{alg_hhd} for computing \ref{eqn_hhd}s, {looking at the influence of the ranks, the order of the tensor, the CPD decomposition algorithm, and additive noise}. {Finally, in \cref{sec_sub_rbm}, we give an example application of \cref{alg_hhd} to a synthetic pdf of an underlying RBM model, and discuss the limitations of it in this high-noise setting.} 

\subsection{Implementation details}\label{sec_sub_impl_det}
All algorithms were implemented in Julia, {version 1.11.5 (2025-04-14), using Julia's packaged OpenBLAS library} as BLAS implementation.
We implemented the algorithms following the clarifications in the text on how certain operations can be implemented efficiently.
The main code relies on {the \texttt{Distributions.jl} (version 0.25.119) package \cite{JSSv098i16,Distributions.jl-2019}}.
Moreover, we compute the CPD in our implementation using \texttt{cpd\_hnf}, which is an {accurate} non-iterative algorithm based on normal form computations rather than numerical optimization, described in \cite{TV2022}. It additionally requires the packages \texttt{Arpack.jl} {(version 0.5.3)} and \texttt{DynamicPolynomials.jl} {(version 0.6.1)}.
Our Julia implementation, including a module \texttt{HHDExperiments.jl} to reproduce the experiments, is available at
\url{https://gitlab.kuleuven.be/u0072863/hhd-from-cpd/}.

All experiments were performed on the desktop \emph{aerie} which contains one AMD Ryzen 7 5800X3D 8-core, {16-thread} processor (each core clocked at 3.4GHz) with 96MB of shared L3 cache memory, $4$ Corsair 32GB DDR4-3600 memory modules, and running Xubuntu {24.04.2 LTS} on the {linux 6.8.0-57-generic} (x86\_64) kernel. Julia was restricted to use only the 8 physical cores with exactly one thread per core using \texttt{numactl --physcpubind=0-7 julia} to run Julia. To mitigate Julia's just in time compilation overhead, a small-scale test run is always executed before executing and timing our main experiments. No other computationally significant processes ran during the experiments.

\subsection{Experimental setup}\label{sec_sub_generalexperiments}
We generate random CPDs {and Hadamard multiply them} to create random HHD problem instances. A random rank-$r$ tensor $\tensor{T}$ of size $n_1 \times\cdots\times n_d$ is always sampled as follows. We sample \emph{factor matrices} $A^i \in \RR^{n_i \times r}$ by independently sampling its entries from the standard normal distribution. Then, we set $\tensor{T} = \sum_{j=1}^r \mathbf{a}^1_j\otimes\cdots\otimes\mathbf{a}^d_j$, where $\mathbf{a}^i_j$ is the $j$th column of $A^i$.

It should be noted that ``this procedure {[for sampling random rank-$r$ tensors]} generates tensors that are heavily biased toward being numerically well-conditioned'' for the problem of computing a CPD from numerical data, as discussed in \cite[Section 1.2]{BBV2022}. We do not know if Hadamard multiplying such randomly generated CPDs creates HHD problems that are similarly more well-conditioned than sampling random inputs from a reasonable probability distribution.

\subsection{Rank-1 permutation}\label{sec_sub_rk1permexp}
{The first experiment investigates the computational performance to identify a rank-$1$ permutation from an admissible vector using \cref{alg_rank1_permutation} preceded by the robust minor identification of \cref{alg_robust_minors}.}

The setup is described next. {For each $i = 25, 26,\ldots, 100$, we compute $R_i = 1 + \lfloor 10^{4.5i/100} \rfloor$} and compute its prime factorization $R_i = r_{i1} \cdots r_{im}$ (with possibly repeated factors). If $R_i$ would chance to be a prime number, then we take $R_i+1$ instead. {Thirty} random rank-$1$ CPDs of size $r_{i1}\times\cdots\times r_{im}$ are generated and the tensors {$\tensor{T}_{i}^{(1)},\dots,\tensor{T}_i^{(30)}$ they} represent are computed. Finally, we {compute the $B \in \mathbb{R}^{R_i \times 30}$ whose $j$th column is $B_{:,j} = \mathrm{vec}(\tensor{T}_i^{(j)})$. This matrix, admitting a simultaneous rank-$1$ permutation, is} the input data for \cref{alg_rank1_permutation,alg_robust_minors}. {The relative tolerance for accepting minors with \cref{alg_robust_minors} is $5\cdot10^{-3}$.} The wall clock time for computing the rank-$1$ permutation is recorded.

\begin{figure}[tb]
\centering 
\includegraphics[width=.98\textwidth]{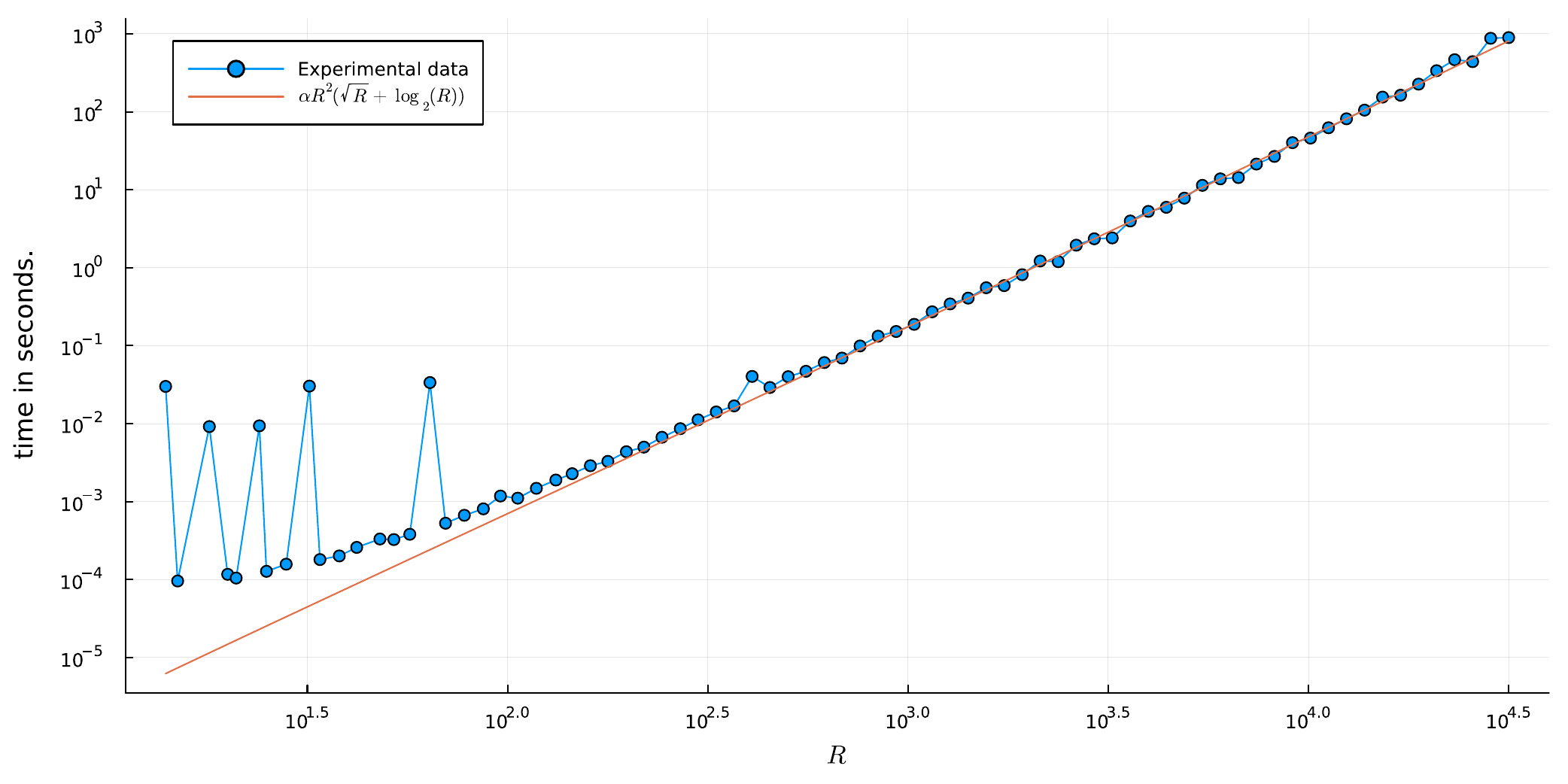}
\caption{The computational time and estimated asymptotic time complexity for computing rank-$1$ permutations of tensors containing $R$ elements as explained in \cref{sec_sub_rk1permexp}.}
\label{fig_times_rk1}
\end{figure}

In \cref{fig_times_rk1}, the time required by our implementation to compute a rank-$1$ permutation with \cref{alg_rank1_permutation,alg_robust_minors} is visualized on a logarithmic scale, along with a curve plotting the estimated asymptotic time complexity.
A few outliers are visible, but these occur for very small problem sizes where the absolute execution time is less than $0.1$ seconds. We attribute these deviations to inconsequential effects, such as Julia's garbage collection overhead.

{We observed in these experiments that the average $h$ from \cref{sec_sub_advanced_version} appears to grow only slightly faster than $\sqrt{R}$. For this reason, we fitted the experimental data to the time complexity estimate $f(R) = \alpha R^2 (\sqrt{R} + \log_2( R ))$.}
The scalar $\alpha$ was determined as the geometric mean of the running times of the largest $30$ problem instances; it was $\alpha \approx 4.2 \cdot 10^{-9}$. The alignment between experimental data and the {time complexity} estimate is {visually good}.

\subsection{Hadamard--Hitchcock decomposition} \label{sec_sub_hhdexper}
Next, we experimentally investigate the computational performance and accuracy of \cref{alg_hhd}. The setup is described next. We consider fifth-order tensors of size $40\times30\times20\times10\times5$, which contain $1.2$ million elements (about $9$MB of data in double-precision floating-point numbers), admitting a random \ref{eqn_hhd}. Such a random tensor is constructed by sampling random CPDs and Hadamard multiplying them. As rank parameters, we consider $2$ Hadamard factors $\br=(r_1,r_2)$ with $1\le r_1 \le r_2 \le 25$. The inequality $r_1 \le r_2$ is no restriction because of the invariance of \ref{eqn_hhd}s under swapping the Hadamard factors. Only $2$ Hadamard factors were chosen in this experiment for ease of visualization. {The random tensors are sampled as described in \cref{sec_sub_generalexperiments}.}

\begin{figure}[tb]
\includegraphics[width=.95\textwidth]{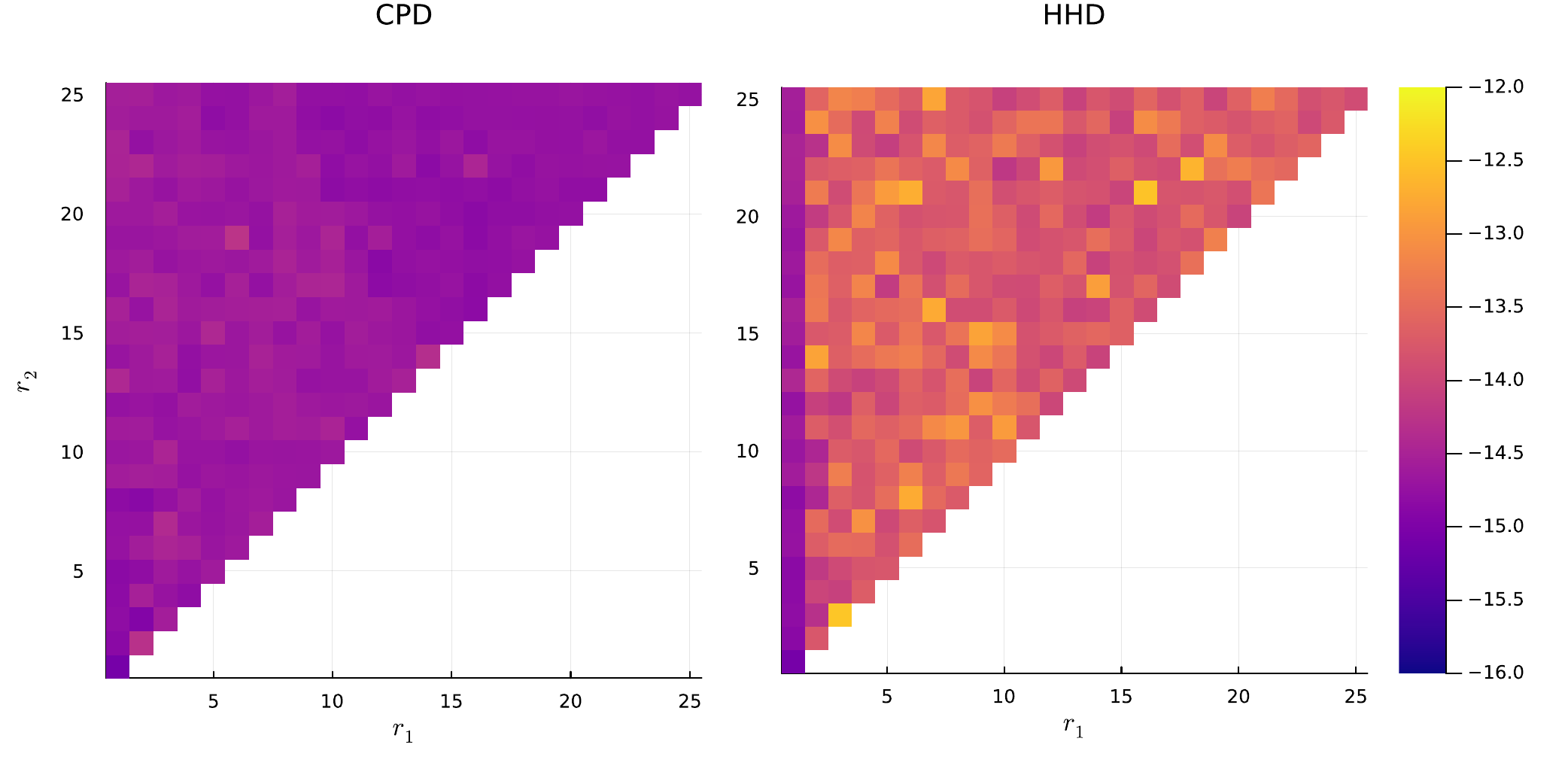}
\caption{Relative backward error of the CPD and HHD for random decomposition problems into $(r_1,r_2)$-HHDs from \cref{sec_sub_hhdexper}.}
\label{fig_accuracy}
\end{figure}

The results of our experiments are summarized in \cref{fig_accuracy}. As measure of accuracy we consider the \emph{backward error}, i.e., the distance between the input tensor and the tensor represented by the numerically computed decomposition. This is a standard way of measuring algorithm performance in numerical analysis \cite{Higham}. Moreover, note that evaluating the \emph{forward error}, as in the case of CPDs, is complicated by the indeterminacy caused by the group action of $(\Var{S},\hadamard)$ on the Hadamard factors. In \cref{fig_accuracy}, the left (respectively right) panel shows the CPD (respectively HHD) relative backward error.
\emph{CPD relative backward error} refers to the relative difference $\frac{\|\tensor{T}-\tilde{\tensor{T}}\|_F}{\|\tensor{T}\|_F}$ in the Euclidean or Frobenius norm between the input tensor $\tensor{T}$ and the numerically computed induced CPD $\tilde{\tensor{T}}=\tilde{\tensor{T}}_1+\cdots+\tilde{\tensor{T}}_R$. The \emph{HHD relative backward error} refers to the same quantity but with respect to the computed \ref{eqn_hhd} $\hat{\tensor{T}}=\hadamard_{k=1}^m \sum_{i_k=1}^{r_k} \hat{\tensor{T}}_{k,i_k}$. Since the \ref{eqn_hhd} is computed from the induced CPD in \cref{alg_hhd}, one should realistically expect the HHD relative error to be at least as large as the CPD error. As can be seen in \cref{fig_accuracy}, the CPD error is generally of the order $10^{-15}$, i.e., within a factor of $10$ of the unit round-off. The HHD error, on the other hand, seems to be systematically a few orders of magnitude worse than the CPD error, though in all cases the relative error is bounded by $10^{-12}$. We believe these are reasonably satisfactory results for a method that is not based on numerical optimization, and taking into account the discrete nature of \cref{alg_robust_minors,alg_rank1_permutation} for determining a rank-$1$ permutation.

\begin{figure}[tb]
\includegraphics[width=.95\textwidth]{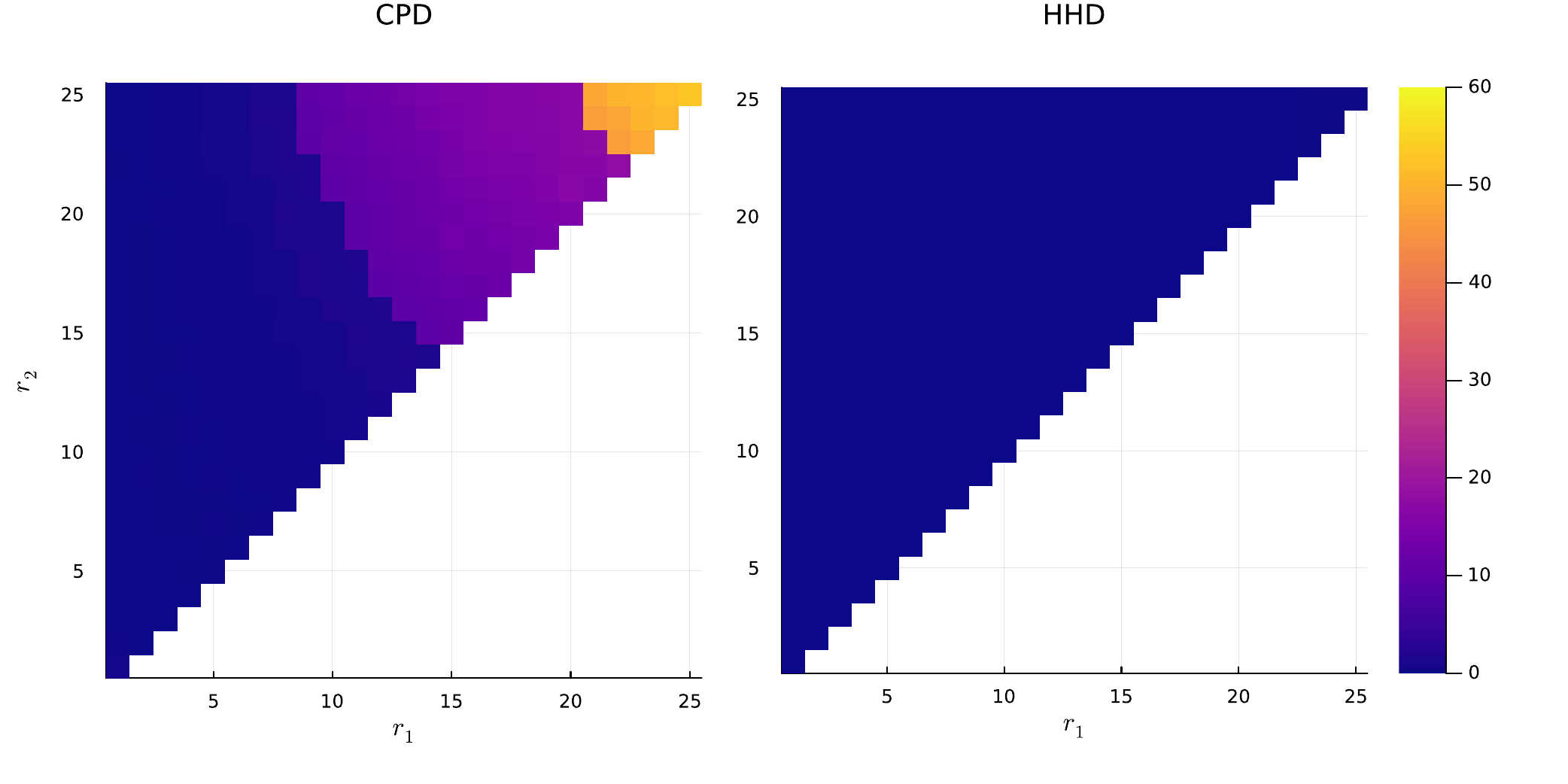}
\caption{The execution time (s) of \cref{alg_hhd} spent in the CPD phase (line \ref{line_cpd}) and the HHD phase (all other lines) for random decomposition problems into rank-$(r_1,r_2)$ HHDs from \cref{sec_sub_hhdexper}.}
\label{fig_times}
\end{figure}

Finally, the individual, non-cumulative computation times for computing the induced CPD and then decomposing that CPD into the HHD is shown in \cref{fig_times}. As can be seen, the time to compute the induced CPD dominates the cost for successive decomposition into the HHD factors; the latter takes less than {$0.5$s} even for the largest problem with $\br=(25,25)$. In the left panel, one notes the discrete jumps between about {$1.5$s}, about $15$s, and approximately {$52$s.} The reason for this is that the normal form algorithm needs to reshape the input tensor in different ways in these regimes. In the first regime it reshapes to {$24000 \times 10 \times 5$, $12000 \times 20 \times 5$, $8000\times 30\times5$, $6000\times40\times5$, $4000\times30\times10$, $3000\times40\times10$ tensors depending on $r_1 r_2$}; in the second to a $1200 \times 200 \times 5$ tensors; and finally in the third regime into a $800 \times 300 \times 5$ tensor. This results in eigendecomposition problems of orders no more than $2200$, $3000$, and $4500$, respectively.

\subsection{More Hadamard factors}\label{sec_sub_highorder}
We show an example output of our Julia implementation, when applied to a computationally challenging scenario. Herein, we constructed a random $50 \times 45 \times 40 \times 35 \times 30$ tensor, containing $N=94.5 \cdot 10^{6}$ elements (about $721$MB in double-precision floating-point data), admitting an \ref{eqn_hhd} with ranks $\br=(3,5,7,9)$. The induced Hitchcock decomposition thus has length $945$. The output, cosmetically formatted, is shown in \cref{tab_output}.

\begin{table}[tb]\footnotesize
\caption{Example output of our implementation of \cref{alg_hhd}.}
\label{tab_output}
\begin{verbatim}
Performing CPD decomposition:
Grouping [[2, 4, 5], [1], [3]] and reshaped to (47250, 50, 40) in 0.594630238s.
1. Performed ST-HOSVD compression to size (945, 50, 40) in 8.362082862s.
2. Constructed kernel of A_1 of size (945, 2000) in 0.543751736s.
3. Constructed resultant map of size (41000, 42200) in 8.298436845s.
4. Constructed res res' in 186.29313676 seconds
5. Computed cokernel of size (945, 41000) in 1217.982204774s.
6. Constructed multiplication matrices in 2.461266143s.
7. Diagonalized multiplication matrices and extracted solution in 3.368446869s.
8. Refined factor matrices Y and Z in 34.995917954s.
9. Recovered factor matrix X in 0.089853506s.
10. Recovered the full factor matrices in 2.165462028s.
  >> Computed a rank-945 CPD of a (50, 45, 40, 35, 30) tensor in 1465.172525932s.

Performing CPD to HHD decomposition:
1. Computed fixed entry of all rank-1 tensors in 0.000501764s.
2. Found all tensor 1-cross 2x2 minors in 0.175331436s.
3. Reconstructed rank-1 permutation in 0.000355509s.
4. Decoupled 945 rank-1 tensors into [3, 5, 7, 9] rank-1 tensors in 0.001734963s.
  >> Decomposed rank-945 CPD into a rank-[3, 5, 7, 9] HHD in 0.177941911s.

Information:
- CPD relative backward error: 2.6311768894610696e-15
- HHD relative backward error: 7.733489380297621e-13
- Extra lop: 2.4682254475559433
- Total computation time: 2203.069830472s.
\end{verbatim}
\end{table}

The first paragraph contains the output of the normal form CPD algorithm, which is described in detail in \cite{TV2022}. Note in particular that the computation of the cokernel of a $41,000 \times 41,000$ matrix is computationally the most demanding step, taking about $83$\% of the total time.

The second paragraph shows information about how the induced CPD is converted into the \ref{eqn_hhd}. Note in particular how fast this phase is compared to computing the CPD. The reason is that the time complexity for this phase is linear in the ambient dimensions $\bn = (50,45,40,35,30)$ and at most cubic in $R=945$, whereas computing the CPD using the normal form algorithm requires the eigendecomposition of a square matrix of order $n_1 \binom{n_3+1}{2}$. This implies a complexity of about $n_1^3 n_3^6$. Hence, one expects that the first phase is about $n_1^3 n_3^6 \ell^{-1} R^{-2.5} \approx 4 \cdot 10^5$ times more expensive. Not exactly equal to the factor of about $8\cdot10^3$ that we observe in practice, but this is readily explained by two complementary effects. First, the eigendecomposition is computed in \texttt{Arpack.jl}, which efficiently exploits \emph{aerie}'s $8$ physical cores while the decomposition from CPD to HHD was fully implemented in Julia without explicit parallelization. Second, \texttt{Arpack.jl}'s code calls the highly optimized \texttt{arpack-ng} library that is built on top of BLAS/LAPACK. The underlying OpenBLAS library is well recognized for its strong computational throughput \cite{Goto2008}. On the other hand, the code for converting the CPD into HHD, especially the search for the minors, should not be expected to attain high throughput because of the high degree of branching. 

Finally, in the third paragraph the code reports useful information on its accuracy and computational performance. The CPD and HHD relative backward errors were explained in \cref{sec_sub_hhdexper}. Both are small, finding a tensor within a relative error of no more than $1\cdot 10^{-12}$. The ``lop'' is short for ``loss of precision'' and equals the base-$10$ logarithm between the CPD and HHD relative backward errors. It measures how many significant digits are additionally lost when recovering the \ref{eqn_hhd} from the induced CPD. Lower values are better. The total computational time for computing the \ref{eqn_hhd}, including the information, is reported in the last line. This includes the time to compute the relative backward errors (which are computed and reported afterwards only for information and do not influence any of the algorithms), which explains the discrepancy of about $738$s between the total time and the sum of the reported times.

\subsection{Other CPD algorithms}\label{sec_sub_comparisonpba}
{A reviewer asked why we selected the \texttt{cpd\_hnf} algorithm \cite{TV2022} to compute the CPD of a tensor. The main focus of this article is decomposing tensors that admit an exact \ref{eqn_hhd}. In this setting, we prefer linear algebra-based algorithms over optimization-based algorithms because the former are conceptually simpler, require no specification of initial point, tend to have less variance in accuracy and computation time, and require no delicate tuning of parameters such as stopping criteria.
We also preferred accuracy over computational performance, as we wish to highlight that \ref{eqn_hhd}s can be computed accurately from an accurate CPD. The main class of linear algebra-based algorithms for computing CPDs are all instances of \emph{pencil-based algorithms} \cite{BBV2019}. However, that article proved that pencil-based algorithms are numerically unstable and exhibit a notable loss of accuracy on practical instances. The numerical experiments from \cite{TV2022} clearly highlighted the superior accuracy of \texttt{cpd\_hnf} as compared to state-of-the-art linear algebra-based algorithms.}

{We repeated the experiment from \cref{sec_sub_hhdexper} for third-order tensors of shape $225 \times 225 \times 20$ with 2 Hadamard factors with ranks $1 \le r_1 \le r_2 \le 15$, where we compute the CPD with either \texttt{cpd\_hnf} \cite{TV2022} or a standard implementation of a pencil-based algorithm, referred to as \texttt{cpd\_pba3}. The implementation we have chosen takes a third-order rank-$r$ tensor $\tensor{T} \in \mathbb{R}^{n_1 \times n_2 \times n_3}$ and then we:
\begin{enumerate}
 \item compress it with ST-HOSVD so that $(U_1\otimes U_2\otimes U_3)(\tensor{C}) = \tensor{T}$, where the core tensor $\tensor{C}$ has shape $r \times r \times \min\{r, n_3\}$, and $U_i \in \mathbb{R}^{n_i \times r}$, $i=1,2$, and $U_3 \in \mathbb{R}^{n_3\times\min\{r,n_3\}}$ are orthonormal factor matrices---see \cite{VVM2012};
 \item make two random linear combinations of $\tensor{C}$, where the coefficients are sampled independently from a standard normal distribution, resulting in two matrices $X, Y \in \mathbb{R}^{r\times r}$;
 \item compute the generalized eigendecomposition of the pencil $(X^T,Y^T)$ and extract the transposed inverse factor matrix $A^{-T}$ from its eigenvectors;
 \item Compute $\tensor{D} = (A^{-1}\otimes I\otimes I)(\tensor{C}) = \sum_{i=1}^r \vect{e}_i \otimes \vect{b}_i \otimes \vect{c}_i$, where $\vect{e}_i$ is the $i$th standard basis vector of $\mathbb{R}^r$;
 \item extract the vectors $\vect{b}_i, \vect{c}_i$ from a rank-$1$ truncated singular value decomposition of the matrix $\tensor{D}_{i,:,:} \in \mathbb{R}^{r\times \min\{r,n_3\}}$, for all $i=1,\dots,r$.
 \item return the factor matrices $(U_1 A, U_2 B, U_3 C)$, where $B$ and $C$ respectively contain the $\vect{b}_i$'s and $\vect{c}_i$'s as columns.
\end{enumerate}}

\begin{figure}[tb]
\centering
\begin{minipage}{\textwidth} \centering
\includegraphics[width=.95\textwidth]{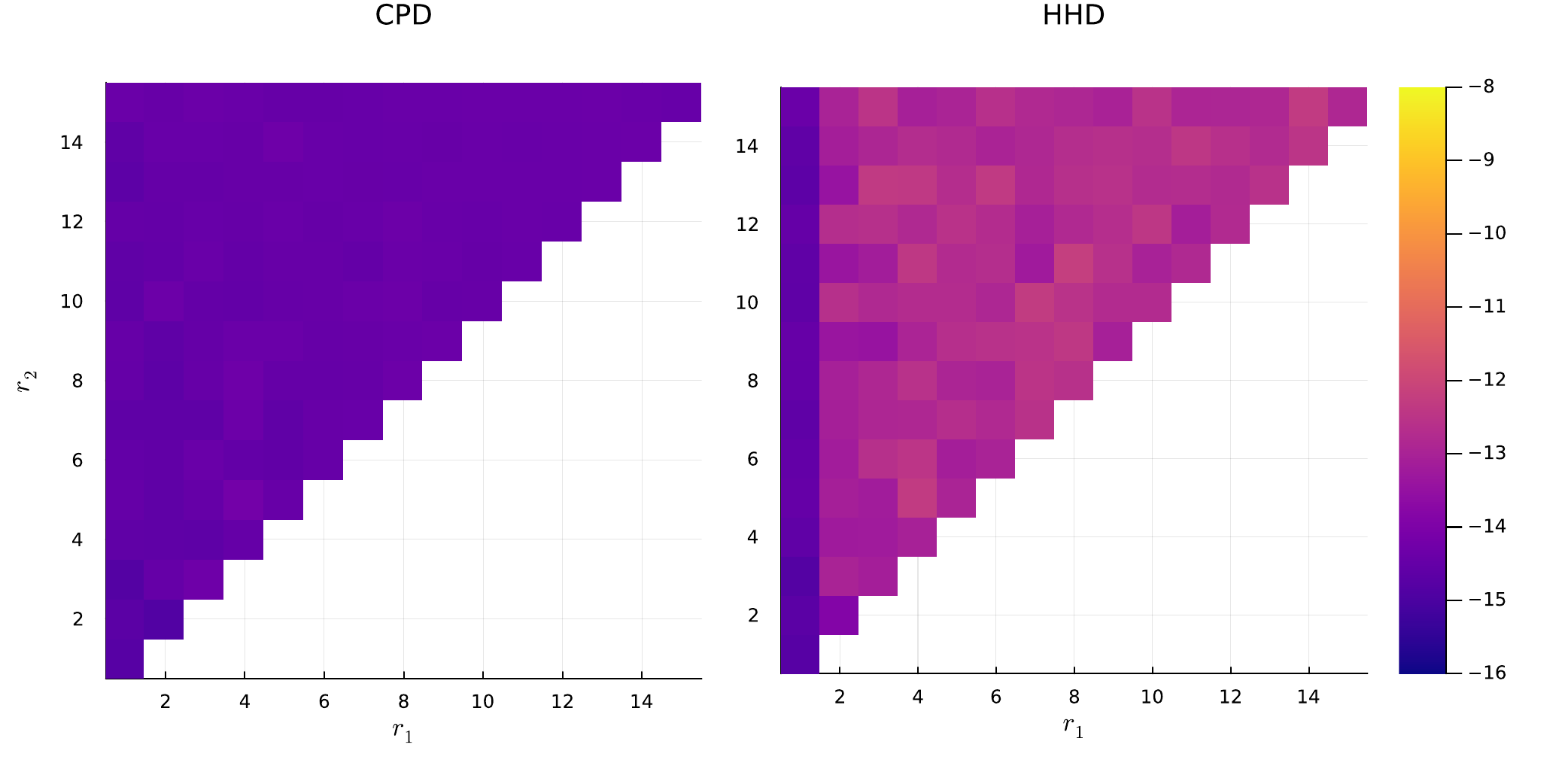}

(a) \texttt{cpd\_hnf}
\end{minipage}

\begin{minipage}{\textwidth} \centering
\includegraphics[width=.95\textwidth]{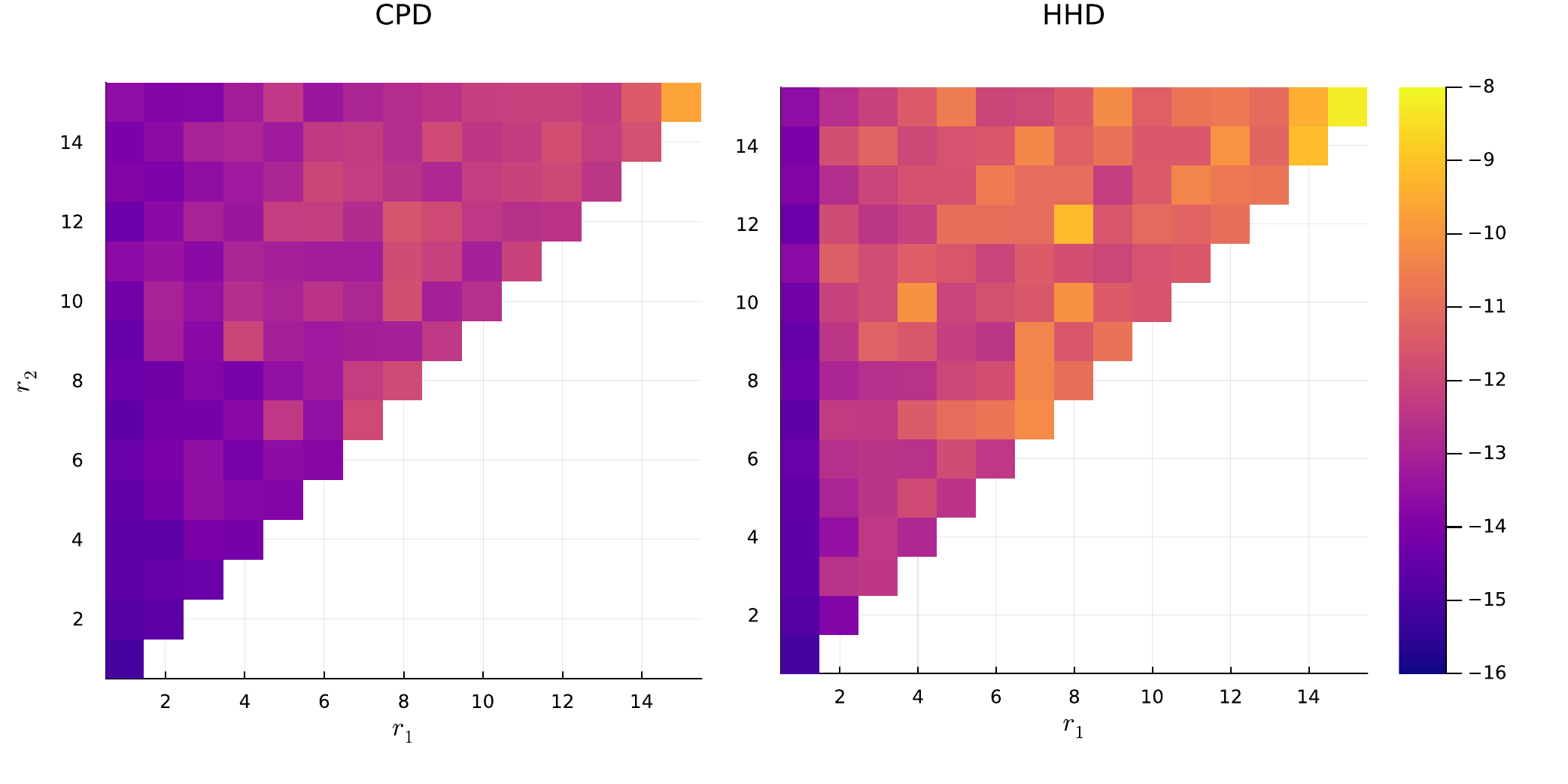}

(b) \texttt{cpd\_pba3}
\end{minipage}
\caption{Relative backward error of the CPD and HHD for random decomposition problems into $(r_1,r_2)$-HHDs from \cref{sec_sub_comparisonpba}, where the CPD is computed with two different algorithms. The scale is the same in all plots.}
\label{fig_accuracy_comparison}
\end{figure}

{The observed relative backward errors are shown in \cref{fig_accuracy_comparison}. The left panels corroborate the theory of \cite{BBV2019}, showing that pencil-based algorithms are not stable. By contrast, \texttt{cpd\_hnf} admits relative backward errors of approximately $10^{-15}$ uniformly across all selected tensors. The right panels show that the part of the algorithm that takes a CPD and further decomposes it into an \ref{eqn_hhd} then causes a further median increase of the relative backward error by a factor of approximately $40.6$ in the case of \texttt{cpd\_hnf} and a factor of about $23.7$ for \texttt{cpd\_pba3}.}

\begin{figure}[tb]
\includegraphics[width=.95\textwidth]{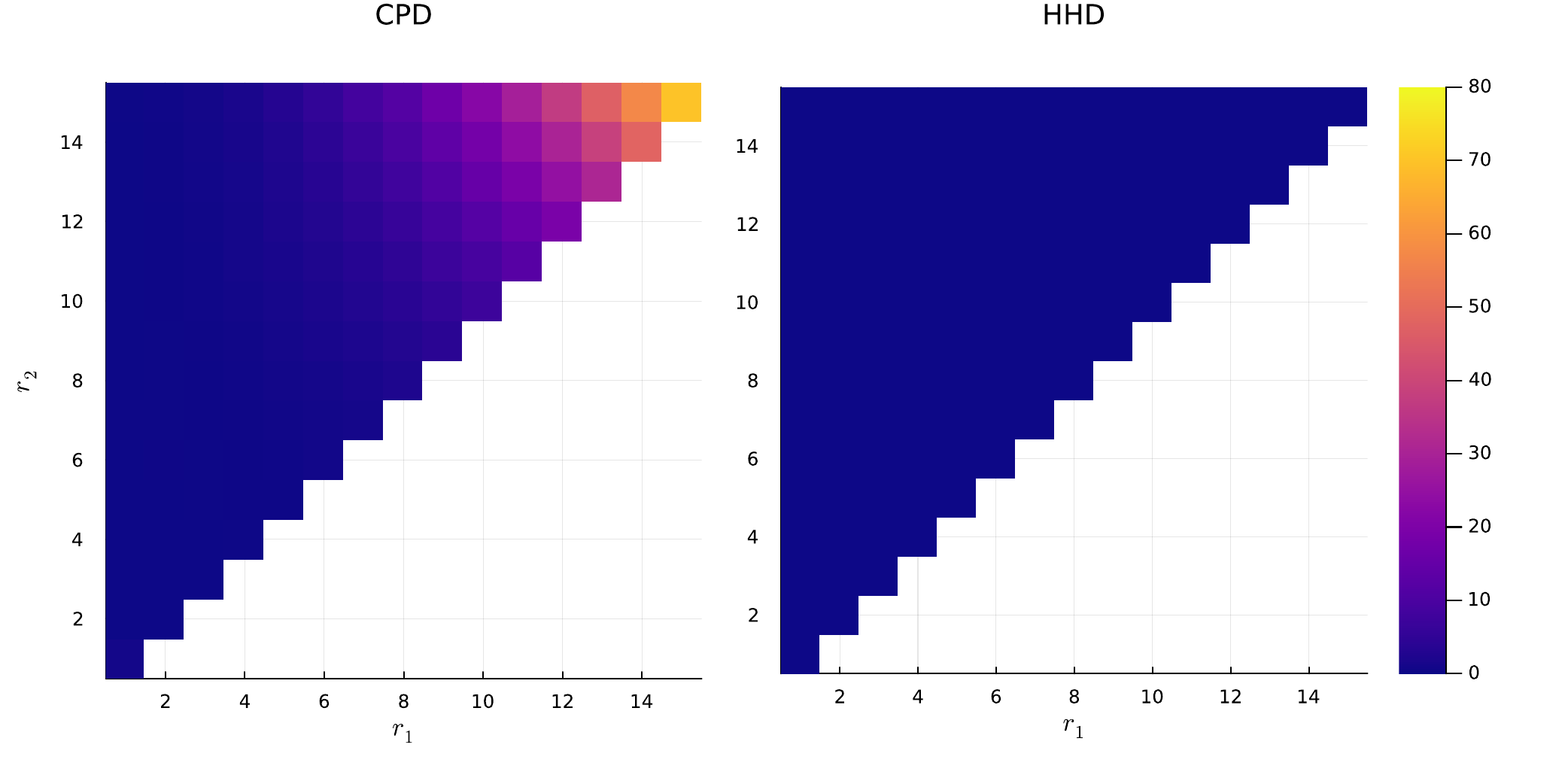}
\caption{The execution time (s) of \cref{alg_hhd} spent in the CPD phase (line \ref{line_cpd}) and the HHD phase (all other lines) for random decomposition problems into rank-$(r_1,r_2)$ HHDs from \cref{sec_sub_comparisonpba}.}
\label{fig_times_comparison}
\end{figure}

{The pencil-based algorithm is significantly faster, however, taking no more than $0.2$ seconds for each of the problem instances. By contrast, the execution times of \texttt{cpd\_hnf} quickly start to rise for $(r_1,r_2) \ge (10,10)$, as can be seen in \cref{fig_times_comparison}.}

{Based on these experiments, we suggest to refine the output of a pencil-based algorithm using a few steps of an optimization-based method for practical setups. This will probably offer a better trade-off between accuracy and computational complexity than using only \texttt{cpd\_hnf} or \texttt{cpd\_pba}.}

\subsection{Noise tolerance}\label{sec_sub_noise_tolerance}
{The next experiment investigates the robustness of \cref{alg_hhd} to mild violations of the model in \cref{eqn_hhd}, as suggested by a reviewer. That is, we consider slightly perturbed \ref{eqn_hhd} tensors of the form
\[
 \widehat{\tensor{T}} = \tensor{T} + \eta \frac{\|\tensor{T}\|_F}{\|\tensor{N}\|_F} \tensor{N},
\]
where $\tensor{T}$ admits a random \ref{eqn_hhd} as described in the general experimental setup in \cref{sec_sub_generalexperiments}, $\eta \approx 0$ is the relative magnitude of the perturbation, and $\tensor{N}$ is a noise tensor. In our experiments, $\tensor{N}$ is white Gaussian noise, i.e., all entries $\tensor{N}_{i_1,\dots,i_d}$ are independent and identically distributed samples from the standard normal distribution $N(0,1)$. For $\tensor{T}$ we consider random HHD order-$d$ tensors ($d=3,4,5$) of shape $10 \times\dots\times 10$ with $2$ Hadamard factors and ranks $\vect{r} = (3, 2)$.}

\begin{figure}[tb]
%
\includegraphics[width=.95\textwidth]{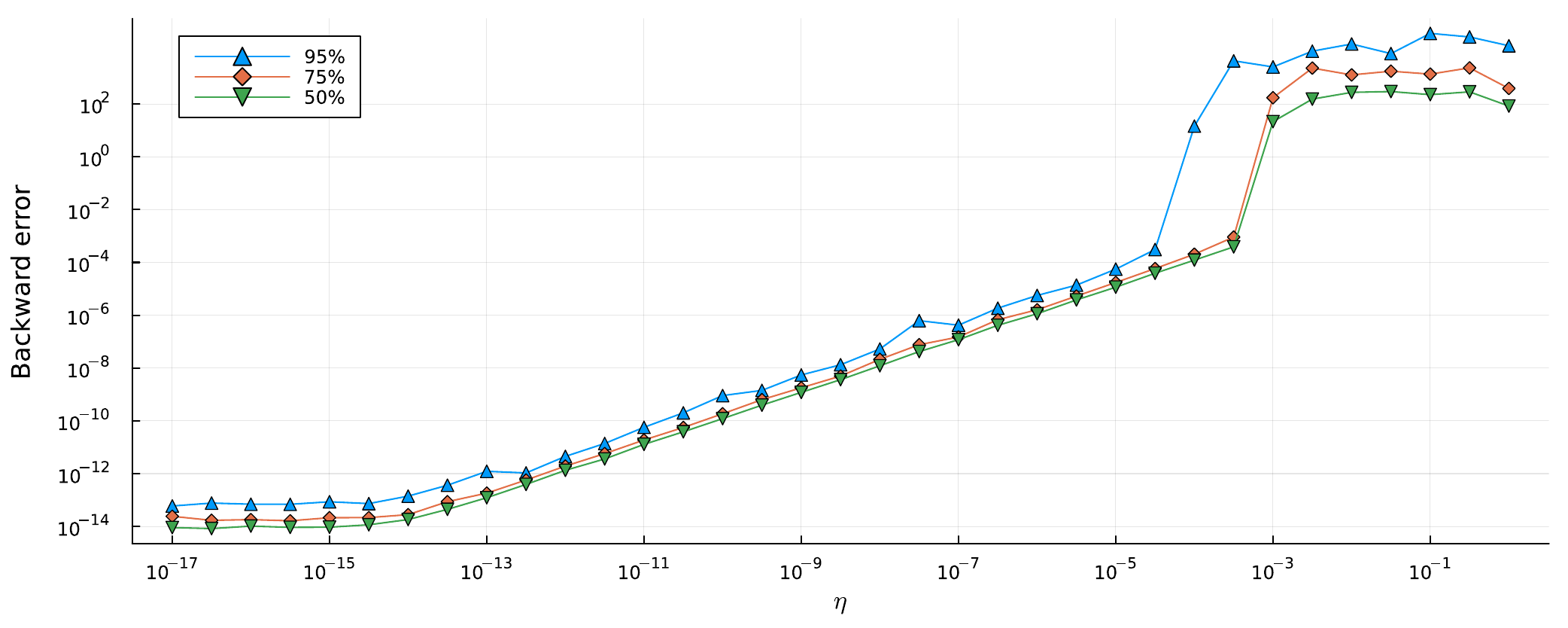}

\caption{Percentiles of the relative backward error for computing the HHD of random noisy problem instances of shape $10\times 10\times 10\times 10 \times 10$ and rank $(3,2)$, as described in \cref{sec_sub_noise_tolerance}.}
\label{fig_noise}
\end{figure}

{For each relative perturbation magnitude $\eta = 10^{-17}, 10^{-16.5}, 10^{-16}, \dots, 10^{0}$, we generate $100$ random problem instances, as described above. The thusly sampled tensors $\widehat{\tensor{T}}$ are then decomposed by \cref{alg_hhd}, recording the relative backward errors.
We summarize the results of these experiments by plotting the $95$th, $75$th, and $50$th (i.e., the median) percentile of these relative backward errors in \cref{fig_noise}. As the results for $d=3,4,5$ all look similar, we chose to include only the case $d=5$.}

{\Cref{alg_hhd} was derived with an exact Hadamard--Hitchcock decomposition as in \cref{eqn_hhd} in mind. It is not an optimization-based algorithm designed to yield the closest \ref{eqn_hhd} to a noisy tensor like $\widehat{\tensor{T}}$. Nevertheless, \cref{fig_noise} shows that \cref{alg_hhd} is reasonably robust to mild perturbations of magnitude less than $10^{-4}$. This robustness strongly relies on the robust minor identification using \cref{alg_robust_minors}. Lowering $\ell$, the number of admissible vectors used, negatively impacts the robustness to noise. In particular, based on our experiments, $\ell \le 10$ is not recommended, unless no noise is present ($\eta=0$).}

{We observed in several additional  experiments that the sudden jump of the median relative backward error near $\eta \approx 5 \cdot 10^{-3}$ in \cref{fig_noise}, in general seems to occur near the chosen relative tolerance $\epsilon$ for the identification of minors in \cref{alg_robust_minors}. The jump can be further delayed by choosing $\epsilon$ up to $10^{-1}$, though we observed that this makes the identification of the vanishing minors much less robust for $\eta \ll \epsilon$. As \cref{eqn_hhd} is intended for the model \cref{eqn_hhd} without noise, we suggest $\epsilon=5\cdot 10^{-3}$ and $\ell=50$ as a good compromise.}

\subsection{A synthetic RBM example} \label{sec_sub_rbm}
{One of the reviewers asked if \cref{alg_hhd} can be applied to an example of RBMs, which was the original inspiration for investigating \ref{eqn_hhd}s. While \cref{alg_hhd} assumes that the \ref{eqn_hhd} model holds exactly, the experiments with noise from the previous subsection suggest that sufficiently accurate approximations of the target tensor may be decomposed successfully.}

{For this experiment, we assume that the phenomenon of interest follows a RBM model exactly. We have chosen five observable variables $X_i \sim D_{16-i}$, $i=1,\dots,5$, and two hidden variables $Y_i \in D_{5-i}$, $i=1,2$. Herein, $D_{k}$ is a distribution on a discrete space with $k$ elements. To obtain a variety of distributions with a nice spread of probabilities, we sampled the distributions $D_k$ randomly from the probability $(k-1)$-simplex by independently drawing $k$ numbers from a $\chi^2_2$ distribution with $2$ degrees of freedom, and then normalizing them.}

{With this random RBM in place, we can realistically assume that we can sample from the population joint pdf $\tensor{P} \in \RR^{15\times14\times13\times12\times11}$ represented by the RBM. We generated $S$ samples from this pdf by \emph{inverse transform sampling} (after imposing an arbitrary order on the lattice $[15]\times[14]\times[13]\times[12]\times[11]$). The samples from the uniform distribution on $[0,1]$ were obtained through \emph{quasi-Monte Carlo sampling} by a \emph{rank-$1$ lattice rule} using $e-2 \approx 0.71828$ as generator.}

{We performed the above process with $S = 10^9$ samples, obtaining a sample estimate $\widehat{\tensor{P}}$ whose relative backward error to $\tensor{P}$ was $9.820738806086711 \cdot 10^{-5}$. Such a high number of samples is required in light of the results of \cref{sec_sub_noise_tolerance} (namely, successful decomposition with high probability requires a relative error below about $10^{-4}$) and because the absolute error between the population pdf $\tensor{P}$ and the empirical estimate $\widehat{\tensor{P}}$ will decay at most like $\sqrt{1/S}$ with high probability as a consequence of \cite[Theorem 2]{VC1971}.
A billion samples is probably too high to be realistic in most applications. However, with additional assumptions on the pdfs in the RBM, more advanced sampling techniques, such as \emph{discrete kernel density estimation}, might generate equally accurate approximations using significantly fewer samples.}

\begin{figure}
 \includegraphics[width=0.99\textwidth]{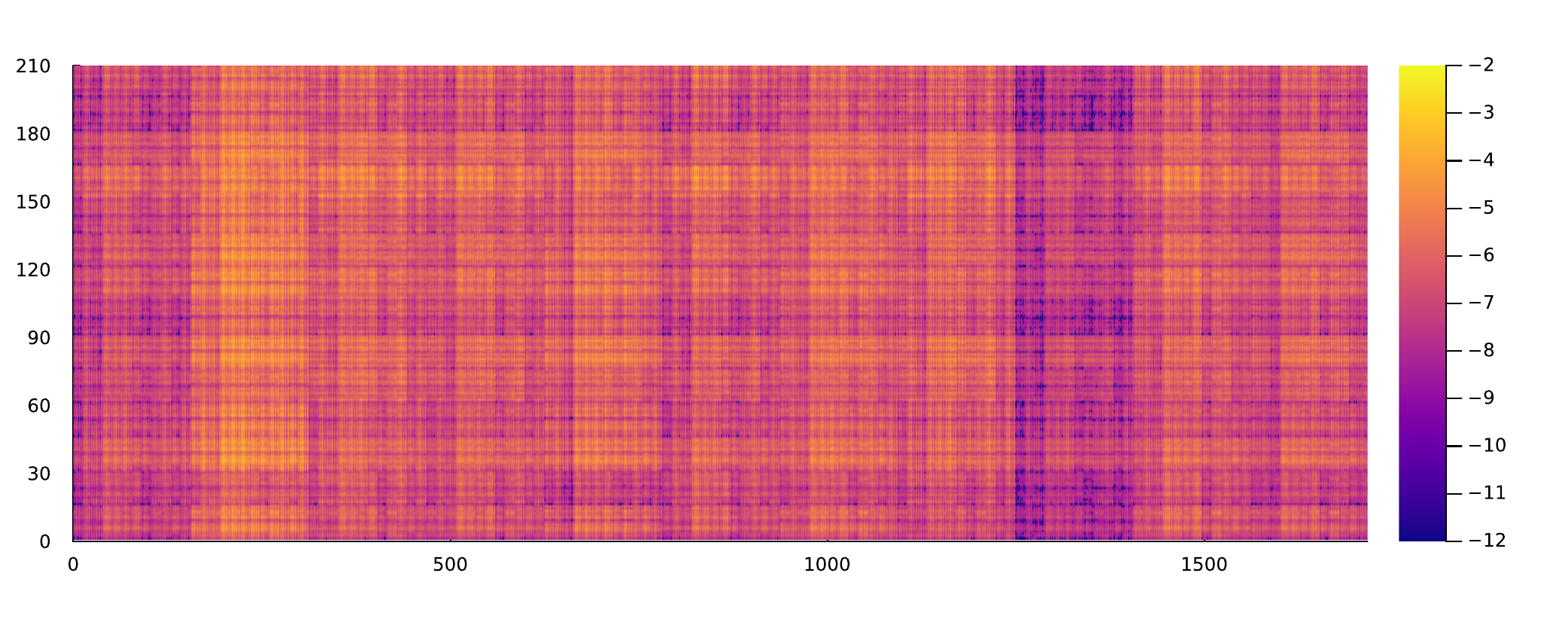}
 \caption{The numerical values of the sample pdf on $[15]\times[14]\times[13]\times[12]\times[11]$ described in \cref{sec_sub_rbm}, flattened as a $210 \times 1716$ matrix and displayed in a base-$10$ logarithmic scale. The visually discontinuous values shaded as if they were $10^{-12}$ all correspond to values that are exactly $0$, i.e., they were never encountered in the sampling process.}\label{fig_rbm_apdf}
\end{figure}

{By reshaping $\widehat{\tensor{P}}$ into an $15\cdot14 \times 13\cdot12\cdot11$ matrix, we can conveniently visualize it. This is shown in \cref{fig_rbm_apdf}. The probabilities lie in a wide range from $0$ to about $6.1 \cdot 10^{-4}$ with the geometric mean (excluding the $2415$ values equal to $0$) being approximately $3.28 \cdot 10^{-7}$.}

{Applying \cref{alg_hhd} to $\widehat{\tensor{P}}$, we observed the following abbreviated output in one of the successful cases:}
{\footnotesize\begin{verbatim}
Performing CPD decomposition:
>> Computed a rank-12 CPD of a (15, 14, 13, 12, 11) tensor in 0.50s.
Performing CPD to HHD decomposition:
>> Decomposed rank-12 CPD into a rank-[4, 3] HHD in 1.74323e-4s.

Information:
- CPD relative backward error: 9.84001177750318e-5
- HHD relative backward error: 1.17948154773631e-4
- Extra lop: 0.07869553276617268
- Total computation time: 1.713093237s.
\end{verbatim}}%
{The above setup will also give rise to unsuccessful runs, as was already observed for high noise levels in \cref{sec_sub_noise_tolerance}. With the above parameters, in $25$ repetitions of the experiment, $17$ of them were ``successful,'' meaning an HHD relative backward error below $2\cdot 10^{-3}$. In all of the failed cases, the CPD relative backward error was less than $10^{-4}$, indicating that \cref{alg_robust_minors} failed to identify the correct minors.}

{One potential reason why the correct rank-$1$ permutation fails to be identified could be due to the presence of many small and clustered numerical values. However, if the true pdf $\tensor{P}$ is supplied to \cref{alg_hhd}, it always succeeds with HHD relative backward error below $10^{-12}$ as in \cref{sec_sub_hhdexper}. For this reason, we do not believe that the small values themselves are responsible for this failure.}

{A contributing factor to why the robust minor identification fails even when the tensor is approximated to a CPD relative backward error of about $10^{-4}$, is because the individual rank-$1$ tensors in the CPD are input to the minor identification. Their relative error is governed by the relative backward error multiplied with the relative \emph{condition number} of the tensor rank decomposition problem \cite{BV2018}. We found in these experiments that the relative condition number was small, between $3$ and $10$. While this will certainly be a contributing factor, we are not convinced that the condition number of the intermediate CPD computation is the sole reason why the correct minors are not identified in the unsuccessful runs. Further research is necessary to pinpoint the exact issue. It could either be an issue of numerical stability, or due to the ill-conditioning of the original HHD decomposition problem.}

\section{Conclusions}\label{sec_conclusions}
Discrete restricted Boltzmann machines are latent variable models whose pdf marginalized over the hidden variables admits an interesting decomposition structure.
This Hadamard--Hitchcock decomposition consists of a Hadamard product of component Hitchcock decompositions.
We investigated the identifiability of \ref{eqn_hhd}s based on their relation to their induced Hitchcock decomposition.
Our main identifiability result, \cref{thm:main_identifiabilityHHD}, shows that tensors admitting a generic \ref{eqn_hhd} whose induced Hitchcock decomposition is $R$-identifiable, is also $\br$-identifiable. Moreover, we showed in \cref{prop_hhdkruskal} that generic \ref{eqn_hhd}s for which $R$ falls in the range of the reshaped Kruskal criterion \cite{AMR2009,COV17} are $R$-identifiable, and, hence, $\br$-identifiable. We then leveraged these identifiability results in \cref{alg_hhd} to decompose a given tensor into an \ref{eqn_hhd} based on an initial rank-$R$ CPD of the tensor under a few genericity assumptions. To establish these results, an essential ingredient was the novel concept of rank-$1$ permutations. Their properties were studied in \cref{sec_rk1permpres} where we related them to rank-$1$ preservers \cite{Westwick1967} and tensor products of permutations \cite{BC2021}. Additionally, we proposed \cref{alg_rank1_permutation,alg_robust_minors} for computing a rank-$1$ permutation of generic inputs, {having at most a cubic} time complexity. Finally, all algorithms were implemented in Julia using floating-point arithmetic. Our numerical experiments in \cref{sec:experiments} show the feasibility of the proposed approach in terms of computation time and accuracy for almost noise-free data.

This paper focused on identifiability and algorithms in the theoretical setting where the tensor exactly admits an \ref{eqn_hhd}. However, in most applications, only an approximation by a low-rank HHD will be meaningful. Much future work lies ahead in this area. For example, to compute approximating HHDs one could develop Riemannian optimization algorithms over the product Segre manifold, similar to those in \cite{BV2018b,SVV2022} for CPDs.

{Our experiments indicate that \cref{alg_robust_minors} works well when the tensor to decompose almost exactly admits an \ref{eqn_hhd}. However, with a more significant model violation, its performance degrades.} More robust ways of identifying the minors, extracting the cross, or alternative algorithms for computing rank-$1$ permutations will therefore be necessary.

The identifiability results we established open the door for an investigation into the condition number of (the computational problem of computing) \ref{eqn_hhd}s. This number quantifies the worst-case local sensitivity of the rank-$1$ components (or the component CPDs, depending on how the decomposition problem is formulated) to perturbations of the tensor. While a study of condition numbers is complicated by action of $(\Var{S},\hadamard)$, one {can apply the techniques} from \cite{DV2023}.

\bibliographystyle{alpha}
\bibliography{HHD}

\newcommand{\etalchar}[1]{$^{#1}$}
\begin{thebibliography}{SVdVV22}

\bibitem[ABC{\etalchar{+}}24]{AB+22}
B.~Atar, K.~Bhaskara, A.~Cook, S.~{Da Silva}, M.~Harada, J.~Rajchgot, A.~{Van
  Tuyl}, R.~Wang, and J.~Yang.
\newblock {Hadamard} products and binomial ideals.
\newblock {\em Journal of Pure and Applied Algebra}, 228(6):107568, 2024.

\bibitem[AMR09]{AMR2009}
E.~S. Allman, C.~Matias, and J.~A. Rhodes.
\newblock Identifiability of parameters in latent structure models with many
  observed variables.
\newblock {\em The Annals of Statistics}, 37(6A), 2009.

\bibitem[AOP09]{abo2009induction}
H.~Abo, G.~Ottaviani, and C.~Peterson.
\newblock Induction for secant varieties of segre varieties.
\newblock {\em Transactions of the American Mathematical Society},
  361(2):767--792, 2009.

\bibitem[BBV19]{BBV2019}
C.~Beltr\'an, P.~Breiding, and N.~Vannieuwenhoven.
\newblock Pencil-based algorithms for tensor rank decomposition are not stable.
\newblock {\em SIAM Journal on Matrix Analysis and Applications},
  40(2):739--773, 2019.

\bibitem[BBV22]{BBV2022}
C.~Beltr\'an, P.~Breiding, and N.~Vannieuwenhoven.
\newblock The average condition number of most tensor rank decomposition
  problems is infinite.
\newblock {\em Foundations of Computational Mathematics}, 23(2):433--491, 2022.

\bibitem[BC19]{BC2019}
C.~Bocci and L.~Chiantini.
\newblock {\em {An Introduction to Algebraic Statistics with Tensors}}.
\newblock Number 118 in UNITEXT. Springer Cham, 2019.

\bibitem[BC21]{BC2021}
F.~Brenti and R.~Conti.
\newblock Permutations, tensor products, and {Cuntz} algebra automorphisms.
\newblock {\em Advances in Mathematics}, 381:107590, 2021.

\bibitem[BC22]{BC22}
C.~Bocci and E.~Carlini.
\newblock {Hadamard} products of hypersurfaces.
\newblock {\em Journal of Pure and Applied Algebra}, 226(11):107078, 2022.

\bibitem[BCC{\etalchar{+}}18]{bernardi2018hitchhiker}
A.~Bernardi, E.~Carlini, M.~V. Catalisano, A.~Gimigliano, and A.~Oneto.
\newblock The hitchhiker guide to: {Secant} varieties and tensor decomposition.
\newblock {\em Mathematics}, 6(12):314, 2018.

\bibitem[BCK16]{BCK16}
C.~Bocci, E.~Carlini, and J.~Kileel.
\newblock {Hadamard} products of linear spaces.
\newblock {\em Journal of Algebra}, 448:595--617, 2016.

\bibitem[Bis06]{Bishop2006}
C.~M. Bishop.
\newblock {\em {Pattern Recognition and Machine Learning}}.
\newblock Springer, 2006.

\bibitem[BPA{\etalchar{+}}21]{JSSv098i16}
M.~Besançon, T.~Papamarkou, D.~Anthoff, A.~Arslan, S.~Byrne, D.~Lin, and
  J.~Pearson.
\newblock Distributions.jl: Definition and modeling of probability
  distributions in the {JuliaStats} ecosystem.
\newblock {\em Journal of Statistical Software}, 98(16):1--30, 2021.

\bibitem[BV18a]{BV2018}
P.~Breiding and N.~Vannieuwenhoven.
\newblock The condition number of join decompositions.
\newblock {\em SIAM Journal on Matrix Analysis and Applications},
  39(1):287--309, 2018.

\bibitem[BV18b]{BV2018b}
P.~Breiding and N.~Vannieuwenhoven.
\newblock A {Riemannian} trust region method for the canonical tensor rank
  approximation problem.
\newblock {\em SIAM Journal on Optimization}, 28(3):2435--2465, 2018.

\bibitem[CCFL20]{CCFL20}
G.~Calussi, E.~Carlini, G.~Fatabbi, and A.~Lorenzini.
\newblock On the {Hadamard} product of degenerate subvarieties.
\newblock {\em Portugaliae Mathematica}, 76(2):123--141, 2020.

\bibitem[CGM24]{CGM2024}
M.~Ciaperoni, A.~Gionis, and H.~Mannila.
\newblock The hadamard decomposition problem.
\newblock {\em Data Mining and Knowledge Discovery}, 38(4):2306--2347, 2024.

\bibitem[CMS10]{CMS10}
M.~A. Cueto, J.~Morton, and B.~Sturmfels.
\newblock Geometry of the restricted {Boltzmann} machine.
\newblock {\em Algebraic Methods in Statistics and Probability}, 516:135--153,
  2010.

\bibitem[CO12]{CO2012}
L.~Chiantini and G.~Ottaviani.
\newblock On generic identifiability of 3-tensors of small rank.
\newblock {\em SIAM Journal on Matrix Analysis and Applications},
  33(3):1018--1037, 2012.

\bibitem[COV14]{COV2014}
L.~Chiantini, G.~Ottaviani, and N.~Vannieuwenhoven.
\newblock On generic and low-rank specific identifiability of complex tensors.
\newblock {\em SIAM Journal on Matrix Analysis and Applications},
  35(4):1265--1287, 2014.

\bibitem[COV17]{COV17}
L.~Chiantini, G.~Ottaviani, and N.~Vannieuwenhoven.
\newblock Effective criteria for specific identifiability of tensors and forms.
\newblock {\em SIAM Journal on Matrix Analysis and Applications},
  38(2):656--681, 2017.

\bibitem[CTY10]{CTY10}
M.~A. Cueto, E.~A. Tobis, and J.~Yu.
\newblock An implicitization challenge for binary factor analysis.
\newblock {\em Journal of Symbolic Computation}, 45(12):1296--1315, 2010.

\bibitem[CZPA09]{NMTFBook2009}
A.~Cichocki, R.~Zdunek, A.-H. Phan, and S.-I. Amari.
\newblock {\em {Nonnegative Matrix and Tensor Factorizations: Applications to
  Exploratory Multi‐Way Data Analysis and Blind Source Separation}}.
\newblock John Wiley \& Sons, Ltd., 2009.

\bibitem[DJL25]{DJL2025}
Z.~Ding, T.~Ji, and M.~Li.
\newblock Improved {Hadamard} decomposition and its application in data
  compression in new-type power systems.
\newblock {\em Mathematics}, 13(4):671, 2025.

\bibitem[dL08]{dSL2008}
V.~{de Silva} and L.-H. Lim.
\newblock Tensor rank and the ill-posedness of the best low-rank approximation
  problem.
\newblock {\em SIAM Journal on Matrix Analysis and Applications},
  30(3):1084--1127, 2008.

\bibitem[DV24]{DV2023}
N.~Dewaele and N.~Vannieuwenhoven.
\newblock Which constraints of a numerical problem cause ill-conditioning?
\newblock {\em Numerische Mathematik}, 156(4):1427--1453, 2024.

\bibitem[EGH09]{EGH2009}
M.~Espig, L.~Grasedyck, and W.~Hackbusch.
\newblock Black box low tensor-rank approximation using fiber-crosses.
\newblock {\em Constructive Approximation}, 30:557--597, 2009.

\bibitem[FI12]{fischer2012introduction}
A.~Fischer and C.~Igel.
\newblock An introduction to restricted {Boltzmann} machines.
\newblock In {\em Progress in Pattern Recognition, Image Analysis, Computer
  Vision, and Applications: 17th Iberoamerican Congress, CIARP 2012, Buenos
  Aires, Argentina, September 3-6, 2012. Proceedings 17}, pages 14--36.
  Springer, 2012.

\bibitem[FOW17]{FOW17}
N.~Friedenberg, A.~Oneto, and R.~L. Williams.
\newblock {Minkowski} sums and {Hadamard} products of algebraic varieties.
\newblock In {\em Combinatorial algebraic geometry}, pages 133--157. Springer,
  2017.

\bibitem[Ges13]{Gesmundo2013}
F.~Gesmundo.
\newblock An asymptotic bound for secant varieties of {Segre} varieties.
\newblock {\em Annali dell' Universit\'a di Ferrara}, 59(2):285--302, 2013.

\bibitem[Gil20]{Gillis2020}
N.~Gillis.
\newblock {\em {Nonnegative Matrix Factorization}}.
\newblock SIAM, Philadelphia, PA, 2020.

\bibitem[Gre78]{Greub1978}
W.~Greub.
\newblock {\em {Multilinear Algebra}}.
\newblock Springer-Verlag, 2 edition, 1978.

\bibitem[Gv08]{Goto2008}
K.~Goto and R.~A. {van de Geijn}.
\newblock Anatomy of high-performance matrix multiplication.
\newblock {\em ACM Transactions on Mathematical Software}, 34(3):1--24, 2008.
\newblock Article 12.

\bibitem[GV13]{GvL4}
G.~H. Golub and C.~F. {Van Loan}.
\newblock {\em Matrix Computations, 4th Edition}.
\newblock The Johns Hopkins University Press, Baltimore, MD, USA, 2013.

\bibitem[Har70]{Harshman1970}
R.~A. Harshman.
\newblock Foundations of the {PARAFAC} procedure: {Models} and conditions for
  an ``explanatory'' multi-modal factor analysis.
\newblock {\em UCLA Working Papers in Phonetics}, 16:1--84, 1970.

\bibitem[Har92]{Harris1992}
J.~Harris.
\newblock {\em {Algebraic Geometry, A First Course}}, volume 133 of {\em
  Graduate Text in Mathematics}.
\newblock Springer-Verlag, 1992.

\bibitem[Hig96]{Higham}
N.~J. Higham.
\newblock {\em {Accuracy and Stability of Numerical Algorithms}}.
\newblock Society for Industrial and Applied Mathematics, 2 edition, 1996.

\bibitem[Hit27]{Hit27}
F.~L. Hitchcock.
\newblock The expression of a tensor or a polyadic as a sum of products.
\newblock {\em Journal of Mathematics and Physics}, 6(1-4):164--189, 1927.

\bibitem[Ho95]{RDF1995}
T.~K. Ho.
\newblock Random decision forests.
\newblock In {\em Proceedings of 3rd International Conference on Document
  Analysis and Recognition}, volume~1 of {\em ICDAR-95}, pages 278--282. IEEE
  Comput. Soc. Press, 1995.

\bibitem[HOT06]{hinton2006fast}
G.~E. Hinton, S.~Osindero, and Y.-W. Teh.
\newblock A fast learning algorithm for deep belief nets.
\newblock {\em Neural computation}, 18(7):1527--1554, 2006.

\bibitem[KB09]{KB09}
T.~G. Kolda and B.~W. Bader.
\newblock Tensor decompositions and applications.
\newblock {\em SIAM Review}, 51(3):455--500, 2009.

\bibitem[Kru77]{kruskal1977three}
J.~B. Kruskal.
\newblock Three-way arrays: rank and uniqueness of trilinear decompositions,
  with application to arithmetic complexity and statistics.
\newblock {\em Linear Algebra and its Applications}, 18(2):95--138, 1977.

\bibitem[Lan12]{Lan12}
J.~M. Landsberg.
\newblock {Tensors: Geometry and Applications}.
\newblock {\em Representation Theory}, 381(402):3, 2012.

\bibitem[Lau96]{lauritzen1996graphical}
S.~L. Lauritzen.
\newblock {\em {Graphical Models}}, volume~17.
\newblock Clarendon Press, 1996.

\bibitem[Lim21]{Lim2021}
L.-H. Lim.
\newblock Tensors in computations.
\newblock {\em Acta Numerica}, 30:555--764, 2021.

\bibitem[LRA93]{LRA1993}
S.~E. Leurgans, R.~T. Ross, and R.~B. Abel.
\newblock A decomposition for three-way arrays.
\newblock {\em SIAM Journal on Matrix Analysis and Applications},
  14(4):1064--1083, 1993.

\bibitem[LWB{\etalchar{+}}19]{Distributions.jl-2019}
D.~Lin, J.~M. White, S.~Byrne, D.~Bates, A.~Noack, J.~Pearson, A.~Arslan,
  K.~Squire, D.~Anthoff, T.~Papamarkou, M.~Besan\c{}on, J.~Drugowitsch,
  M.~Schauer, and other contributors.
\newblock {JuliaStats/Distributions.jl: a Julia package for probability
  distributions and associated functions}, July 2019.

\bibitem[MM15]{MM15}
G.~Mont{\'u}far and J.~Morton.
\newblock Discrete restricted {Boltzmann} machines.
\newblock {\em Journal of Machine Learning Research}, 16(1):653--672, 2015.

\bibitem[MM17]{MM17}
G.~Mont{\'u}far and J.~Morton.
\newblock Dimension of marginals of {Kronecker} product models.
\newblock {\em SIAM Journal on Applied Algebra and Geometry}, 1(1):126--151,
  2017.

\bibitem[MM24]{MM2022}
A.~Massarenti and M.~Mella.
\newblock {Bronowski}'s conjecture and the identifiability of projective
  varieties.
\newblock {\em Duke Mathematical Journal}, 173(17), 2024.

\bibitem[Mon16]{Mon16}
G.~Mont{\'u}far.
\newblock Restricted {Boltzmann} machines: Introduction and review.
\newblock In {\em Information Geometry and Its Applications IV}, pages 75--115.
  Springer, 2016.

\bibitem[PFS16]{PFS2016}
E.~E. Papalexakis, C.~Faloutsos, and N.~D. Sidiropoulos.
\newblock Tensors for data mining and data fusion: models, applications, and
  scalable algorithms.
\newblock {\em ACM Transactions on Intelligent Systems and Technology},
  8(16):1--44, 2016.

\bibitem[Rho10]{Rhodes2010}
J.~A. Rhodes.
\newblock A concise proof of {Kruskal's} theorem on tensor decomposition.
\newblock {\em Linear Algebra and its Applications}, 432(7):1818--1824, 2010.

\bibitem[SM18]{SM18}
A.~Seigal and G.~Mont{\'u}far.
\newblock Mixtures and products in two graphical models.
\newblock {\em Journal of Algebraic Statistics}, 9(1), 2018.

\bibitem[Smo86]{smolensky1986information}
P.~Smolensky.
\newblock Information processing in dynamical systems: Foundations of harmony
  theory.
\newblock Technical report, Department of Computer Science, University of
  Colorado, Boulder, 1986.

\bibitem[Stu97]{sturmfels1996equations}
B.~Sturmfels.
\newblock Equations defining toric varieties.
\newblock In {\em Algebraic geometry---Santa Cruz 1995}, volume 62, Part 2 of
  {\em Proc. Sympos. Pure Math.}, pages 437--449, Providence, RI, 1997.
  American Mathematical Society.

\bibitem[Sul18]{sullivant2018algebraic}
S.~Sullivant.
\newblock {\em {Algebraic Statistics}}, volume 194.
\newblock American Mathematical Society, 2018.

\bibitem[SVdVV22]{SVV2022}
L.~Swijsen, J.~Van~der Veken, and N.~Vannieuwenhoven.
\newblock Tensor completion using geodesics on {Segre} manifolds.
\newblock {\em Numerical Linear Algebra with Applications}, 29(6), April 2022.

\bibitem[TBC23]{blomenhofer2023nondefectivity}
A.~Taveira~Blomenhofer and A.~Casarotti.
\newblock Nondefectivity of invariant secant varieties.
\newblock {\em arXiv:2312.12335v2}, 2023.

\bibitem[TV22]{TV2022}
S.~Telen and N.~Vannieuwenhoven.
\newblock A normal form algorithm for tensor rank decomposition.
\newblock {\em {ACM} Transactions on Mathematical Software}, 48(4):1--35, 2022.

\bibitem[VC71]{VC1971}
V.~N. Vapnik and A.~Ya. Chervonenkis.
\newblock On the uniform convergence of relative frequencies of events to their
  probabilities.
\newblock {\em Theory of Probability \& Its Applications}, 16(2):264--280,
  1971.

\bibitem[VVM12]{VVM2012}
N.~Vannieuwenhoven, R.~Vandebril, and K.~Meerbergen.
\newblock A new truncation strategy for the higher-order singular value
  decomposition.
\newblock {\em SIAM Journal on Scientific Computing}, 34(2):A1027--A1052, 2012.

\bibitem[Wes67]{Westwick1967}
R.~Westwick.
\newblock Transformations on tensor spaces.
\newblock {\em Pacific Journal of Mathematics}, 23(3):613--620, 1967.

\bibitem[ZP08]{ZP2008}
M.~M. Zavlanos and G.~J. Pappas.
\newblock A dynamical systems approach to weighted graph matching.
\newblock {\em Automatica}, 44(11):2817--2824, 2008.

\end{thebibliography}

\end{document}